\documentclass[a4paper,12pt]{article} 
\usepackage[catalan,spanish,english]{babel} 
\usepackage[psamsfonts]{amssymb} 
\usepackage{cmmib57} 
\usepackage{exscale} 
\usepackage{amsmath} 
\usepackage{amscd} 
\usepackage{latexsym} 
\usepackage{graphicx} 
\newtheorem{stat}{Statement}[section]
\newtheorem{thm}[stat]{Theorem}
\newtheorem{cor}[stat]{Corollary}
\newtheorem{prop}[stat]{Proposition}
\newtheorem{lemma}[stat]{Lemma}
\newtheorem{remark}[stat]{Remark}
\newtheorem{example}[stat]{Example}
\newtheorem{assump}[stat]{Assumption}

\numberwithin{equation}{section}
\newcommand{\qed}{\quad{$\square$}} 
\newcommand {\IN}{\mathbb{N}}
\newcommand {\IR}{\mathbb{R}}

\newcommand {\cB}{\mathcal{B}}
\newcommand {\cC}{\mathcal{C}}
\newcommand {\D}{\mbox{${\mathcal D}$}}

\newcommand {\F}{\mbox{${\mathcal F}$}}

\newcommand {\cO}{\mathcal{O}}

\newcommand{\ep}{\varepsilon} 
\newcommand {\IS}{\mbox{${\mathcal S}$}}
\newcommand{\qbar}{\bar{q}}
\newcommand{\qnbar}{q}
\newcommand{\snbar}{s}
\newcommand{\tbar}{\bar{t}}
\newcommand{\wbar}{\bar{w}}
\newcommand{\ro}{\rho}
\newcommand{\bro}{\bar{\rho}}
\newcommand{\half}{\frac{1}{2}}
\newcommand{\kn}{k_n}

\begin{document}

%%%% TITLE PAGE
\begin{titlepage}
\null
\begin{center}
{ \large\bf H\"older-Sobolev regularity of the solution\\[2mm]
 to the stochastic wave equation\\[2mm]
  in dimension 3}\\[2mm]
\bigskip

by\\
\vspace{7mm}

\begin{tabular}{l@{\hspace{10mm}}l@{\hspace{10mm}}l}
{\sc Robert C. Dalang}$\,^{(\ast)}$  &and&{\sc Marta Sanz-Sol\'e}$\,^{(\ast\ast)}$\\
{\small Institut de Math\'ematiques}     &&{\small Facultat de Matem\`atiques}\\
{\small Ecole Polytechnique F\'ed\'erale}          &&{\small Universitat de Barcelona}\\
{\small Station 8 }            &&{\small Gran Via 585}\\         
{\small CH-1015 Lausanne, Switzerland}                
 &&{\small E-08007 Barcelona, Spain}\\
{\small e-mail: robert.dalang@epfl.ch}      &&{\small e-mail: marta.sanz@ub.edu}\\
\null         
\end{tabular}
\end{center}
\vspace{1cm}

\noindent{\bf Abstract:} We study the sample path regularity of the solution of a stochastic wave equation in spatial dimension $d=3$.
The driving noise is white in time and with a spatially homogeneous covariance defined as a product of a Riesz kernel and a smooth function. We prove that at any fixed time, a.s. the sample paths in the spatial variable belong to certain fractional Sobolev spaces. In addition, for any fixed $x\in\IR^3$, the sample paths in time are H\"older continuous functions. Hence, we obtain joint H\"older continuity in the time and space variables.
Our results rely on a detailed analysis of  properties of the stochastic integral used in the rigourous formulation of the spde, as introduced by Dalang and Mueller (2003). Sharp results on one and two dimensional space and time increments
 of generalized Riesz potentials are a crucial ingredient in the analysis of the problem. For spatial covariances given by Riesz kernels, we show that the H\"older exponents that we obtain are optimal. 
\bigskip

\medskip


\footnotesize
{\begin{itemize}
\item[$^{(\ast)}$] Partially supported by the Swiss National Foundation for Scientific Research.
\item[$^{(\ast\ast)}$] Partially supported by the grant BFM2003-01345 from \textit{Direcci\'on
General de Investigaci\'on, Ministerio de Ciencia y Tecnolog\'{\i}a}, Spain.
\end{itemize}}
\end{titlepage}

%%%%END TITLE PAGE

%%%%%%%%%%INTRODUCTION

\section {\bf Introduction}

This paper studies the stochastic wave equation in spatial dimension $d=3$
\begin{align} \label{0.1}
&\left(\frac{\partial^2}{\partial t^2} - \Delta \right) u(t,x) = \sigma\big(u(t,x)\big) \dot F(t,x) 
+ b \big(u(t,x)\big),\\ \nonumber
&\qquad\quad u(0,x) = v_0(x),
\qquad \frac{\partial}{\partial t}u(0,x) = \tilde v_0(x),  
\end{align} 
where $t\in\, ]0,T]$ for some fixed $T>0$, $x\in\IR^3$ and 
$\Delta$ denotes the  Laplacian on $\IR^3$. The coefficients $\sigma$ and $b$ are Lipschitz continuous functions,
the  process $\dot F$ is the {\it formal} derivative of a Gaussian random field, 
white in time and correlated in space. More precisely, for any $d\ge 1$,
let $\mathcal{D}(\mathbb{R}^{d+1})$ be the space of Schwartz test functions 
(see \cite{schwartz}) and let $\Gamma$ be a non-negative and 
non-negative definite tempered measure on $\IR^d$. 
Then, on some probability space, there exists a Gaussian process  $F = \left(F(\varphi),\ \varphi\in \mathcal{D}(\mathbb{R}^{d+1})\right)$
 with mean zero and covariance functional given by 
\begin{equation}\label{cov1}
   E\big(F(\varphi) F(\psi)\big) = \int_{\mathbb{R}_+} ds \int_{\IR^d} \Gamma(dx) (\varphi(s)*\tilde\psi(s))(x), 
\end{equation}
where $\tilde\psi(s)(x)=\psi(s)(-x)$.

   We are interested in solutions which are  {\it random fields}, that is, real-valued processes $(u(t,x)$, $(t,x)\in[0,T]\times\IR^3)$, that are well
defined for every fixed $(t,x)\in[0,T]\times\IR^3$. We want to study their sample path regularity, both in time and space, and check the optimality
of the results.

   There are different possible approaches to giving a rigourous formulation of the Cauchy problem (\ref{0.1}). However, in all of them  
the fundamental solution associated to the wave operator $\mathcal{L} = \frac{\partial^2}{\partial t^2} - \Delta$ naturally 
 plays an important
role. Since its singularity increases with the spatial dimension $d$, the difficulties in studying regularity of the solutions of the stochastic wave equation increase
accordingly. Moreover, keeping the requirement of obtaining random field solutions amounts to adjusting the roughness of the noise to the
degeneracy of the differential operator which defines the equation. It is only for $d=1$ that it is possible to take 
a space-time white noise as 
random input to (\ref{0.1}), while in higher dimensions a non-degenerate spatial correlation is necessary \cite{dalang,KZ}. 

For $d=1, 2$, the stochastic wave equation driven by space-time white noise, and noise that is white in time but spatially correlated, respectively, is now
quite well understood. We refer the reader to \cite{cabana}, \cite{carmonanualart}, \cite{dalangfrangos}, \cite{KZ}, \cite{mms}, \cite{milletss}, \cite{milletss2},
\cite{pz}, \cite{p}, for a sample of articles on the subject.

For $d=3$, the fundamental solution of the wave equation is the measure defined by 
\begin{equation}
\label{G}
   G(t) = \frac{1}{4 \pi t} \sigma_t,
\end{equation}
for any $t>0$, where $\sigma_t$ denotes the uniform surface measure (with total mass $4\pi t^2$) on the sphere of radius $t \in [0, T]$. Hence,
in the mild formulation of equation (\ref{0.1}),  Walsh's  classical stochastic integration theory developed in \cite{walsh} does
not apply. In fact, this question motivated two different extensions of Walsh's integral, given in \cite {dalang} and \cite{DM}, respectively.

The stochastic integral of  \cite {dalang}, written 
$$
   \int_0^t G(t-s,y)\,Z(s,y)\,M(ds,dy),
$$
requires a non-negative distribution $G$, a second-order stationary
process $Z$ in the integrand, and the integrability condition
\begin{equation}
\label{0.2}
\int_0^T ds \int_{\IR^d} \mu(d\xi)\left\vert\mathcal{F}G(t)(\xi)\right\vert^2  < \infty,
\end{equation}
where $\mu = \mathcal{F}^{-1}\Gamma$, 
among other technical properties. As is shown in Section 5 of \cite {dalang}, with this integral one can obtain existence and uniqueness of a random field solution to (\ref{0.1}),
interpreted in the mild form  (\ref{e3.2'}), in the case where the initial conditions vanish. In this framework, results on the regularity of the law of the solution to the 
stochastic wave equation have been proved in \cite{qss1} and \cite{qss2} (see also \cite{sss}).

 In \cite{DM},  a new extension of
Walsh's stochastic integral based on a functional approach is introduced. Neither the positivity of $G$ nor the stationarity of $Z$ are required (see \cite{DM}, Theorem 6). With this integral,
the authors give a precise meaning to the problem (\ref{0.1}) with non vanishing initial conditions  and coefficient $b\equiv 0$ and obtain existence and uniqueness of a solution
$\left(u(t),\ t\in[0,T]\right)$ which is an $L^2(\IR^3)$--valued stochastic process (Theorem 9 in \cite{DM}). This  is the choice of stochastic integral  that we will use  in this paper to study the stochastic
wave equation (\ref{0.1}).

We consider the particular case of a covariance measure that is absolutely continuous with respect to Lebesgue measure, with density given by
\begin{equation}\label{rde1}
f(x) = \varphi(x)\,k_{\beta}(x),
\end{equation}
where $\varphi$ is a smooth positive function and $k_{\beta}$ denotes the Riesz kernel $k_{\beta}(x)= |x|^{-\beta}$, with $\beta\in\, ]0,2[$ (see Assumption \ref{assumpC}).
Riesz kernels are a natural class of  correlation functions and are already present in previous work on the stochastic heat and wave equations,
for instance in \cite{dalangfrangos}, \cite{dalang}, \cite{KZ}, \cite{milletss}. They provide examples where condition (\ref{0.2}) is satisfied: for these covariances, (\ref{0.2}) is equivalent to the condition $0 < \beta < 2$ (see Example \ref{example2.5}).

Related questions for an equation that is second order in time but with fractional Laplacian in any spatial dimension $d$ and  general covariance measure $\Gamma$ have been considered in  \cite{dalangss},
in the setting of an $L^2$--theory (see \cite{rozovsky}). The results there are shown to be optimal in time. 
We adopt here a similar strategy, but we work in an  $L^q$--framework (see \cite{krylov}), for any $q\ge 2$.
Indeed, the particular structure of the wave equation in dimension $d=3$ makes it possible to go beyond the Hilbert space setting
and to obtain sharp results, both in time and space. 

The main result of the paper is Theorem \ref{holdercont}, stating joint H\"older-continuity in $(t,x)$ of the solution to (\ref{0.1}), together with the analysis of the optimality of the exponents
studied in Section \ref{s3}. The optimal H\"older exponent is the same for the time and space variables: this is an intrinsic property of the d'Alembert operator.  Moreover, this result shows how the driving noise $\dot F$ contributes to the roughness of the sample paths, since it expresses the optimal H\"older exponent in terms of the parameters $\beta$ and $\delta$ appearing in Assumption \ref{assumpC} on $\dot F$ (see Section \ref{A2}). 

Notice that for the stochastic heat equation with Lipschitz coefficients in any spatial dimension $d\ge 1$, joint H\"older-continuity in $(t,x)$ of the sample paths of the solution has been established in \cite{sssarra}. Unlike the stochastic wave equation, the H\"older exponent in the time variable is half that for the spatial variable. This is also an intrinsic property of the heat operator. However, it turns out that effect of the driving noise $\dot F$ on the regularity in the spatial variable is the same for both equations (see Theorem \ref{holdercont} and Remark \ref{heatcomparison}). Similar problems for non-Lipschitz coefficients have been recently tackled in \cite{mps}.

   We should point out that despite the similarities just mentioned, establishing regularity results for the solution of the stochastic wave equation requires fundamentally different methods than those for the stochastic heat equation. Indeed, taking for simplicity $b \equiv 0$ and vanishing initial conditions, equation (\ref{0.1}), written in integral form, becomes
$$
   u(t,x) = \int_0^t \int_{\IR^d} G(t-s,x-v) Z(s,v) \dot F(ds,dv),
$$
where $Z(s,v) = \sigma(u(s,v))$. A spatial increment of the solution is
$$
   u(t,x) - u(t,y) = \int_0^t \int_{\IR^d} (G(t-s,x-v) - G(t-s,y-v))  Z(s,v) \dot F(ds,dv).
$$
When the fundamental solution $G$ is smooth, as in the case of the heat equation, one uses Burholder's inequality to see that
\begin{align} \label{Gincr}
 & E(\vert u(t,x) - u(t,y) \vert^p)\\ \nonumber
  &\qquad \le C\, E\Big(\int_0^t ds \int_{\IR^d} du \int_{\IR^d} dv (G(t-s,x-u) - G(t-s,y-u))  Z(s,u) \\  \nonumber
  &\qquad \qquad\qquad \times f(u-v) Z(s,v) (G(t-s,x-v) - G(t-s,y-v)\Big)^{p/2}.
\end{align}
Then, the {\em smoothness} of $G$, together with {\em integrability} of $Z$, leads to regularity of $u(t,\cdot)$. For the wave equation, $G(t)$ is singular with respect to Lebesgue measure (see (\ref{G})), so this kind of approach is not feasible.

   A different idea is to pass the increments on $G$ in (\ref{Gincr}) onto the factor $Z(s,u)f(u-v) Z(s,v)$, using a change of variables; the right-hand side of (\ref{Gincr}) becomes the sum of 
\begin{equation}\label{term2}
   - \int_0^t ds \int_{\IR^3} G(s,du) \int_{\IR^3} G(s,dv) D^2f(v-u,x-y) E(Z(t-s,x-u)Z(t-s,x-v)),
\end{equation}
where $D^2 f(u,x) = f(u+x) - 2f(u) + f(u-x)$, and of three other terms of similar form (see the proof of Lemma \ref{l2.2} for details). Focussing on the term (\ref{term2}), one checks that in the case where $f(x) = \vert x \vert^{-\beta}$ is a Riesz kernel, 
\begin{equation}\label{D2f}
  |D^2 f(u,x)| \leq c |f''(u)|\, \vert x \vert^2 \leq c \vert u \vert^{-(\beta + 2)} \vert x \vert^2,
\end{equation}
where $f''$ denotes the second order diferential of $f$.
This would lead to the following bound for (\ref{term2}):
\begin{equation}\label{infinity}
   \vert x - y\vert^2 \int_0^t ds \int_{\IR^3} G(s,du) \int_{\IR^3} G(s,dv)\,  \vert v-u \vert^{-(\beta + 2)}.
\end{equation}
The factor $\vert x - y\vert^2$ looks too good to be true, and it is! Indeed, the triple integral is equal to the left-hand side of (\ref{0.2}), and we have already pointed out that this is finite if and only if the exponent $\beta + 2$ is less than $2$. However, this is not the case since $\beta \in\, ]0,2[$!

   Even though the bound (\ref{infinity}) equals $+\infty$, this approach contains the premises of our argument. Indeed, instead of differentiating $f$ twice as in (\ref{D2f}), we shall estimate $D^2 f$ by using a {\em fractional derivative} of order $\gamma$, where $\gamma< 2 - \beta$. It turns out that for $f(u) =  \vert u \vert^{-\beta}$, the fractional derivative $f^{(\gamma)}(u) \simeq \vert u \vert^{-(\beta + \gamma)}$. This leads to the following bound for (\ref{term2}): 
$$
   \vert x - y\vert^\gamma \int_0^t ds \int_{\IR^3} G(s,du) \int_{\IR^3} G(s,dv)\, \vert v-u \vert^{-(\beta + \gamma)}.
$$
The triple integral is now finite since $\beta + \gamma < 2$ and this gives the correct order of regularity for $u(t,\cdot)$. The precise properties of Riesz kernels and rigorous use of their fractional derivatives (or rather, their fractional Laplacians) are given in Lemma \ref{lemB.1} (with $a$ there replaced by $3 - \beta - \gamma$, $b$ replaced by $\gamma$, and $d=3$).   

In short, there are mainly three ideas which have been central to obtaining the results of this paper. First, the smoothing in space of the fundamental solution by means of a regularisation procedure based on time-scaled approximations of the Dirac delta--function (see (\ref{A.1})), and the study of the corresponding smoothed equation, to which we can successfully apply standard techniques of stochastic calculus. Secondly, at the level of the smoothed equation, increments of stochastic integrals, whether in space or in time, initially expressed in terms of increments of the fundamental solution, can be reexpressed in terms of increments of the covariance density of the noise. Using the semigroup property of Riesz potentials (see for instance  \cite{stein}), we implement  the ideas concerning fractional derivatives described above (see Lemma \ref{lemB.1}). With these results, we are able to obtain bounds on one and two dimensional increments, in space and in time, of  certain generalized Riesz potentials  of a smoothed version of the fundamental solution of the wave equation. The sharp character of these estimates leads to the optimality of our results.

The paper is organized as follows.  
In Section \ref{F}, we define the smoothing $G_n$ of the fundamental solution $G$ and prove some of their basic properties. Then we describe precisely the type of stochastic noise we are considering in the paper and prove the above mentioned fractional derivative properties.

In Section \ref{s1}, we study the path properties of the indefinite stochastic integral introduced in \cite{DM}. Briefly stated,
we prove that if the sample paths of the stochastic integrand belong to some fractional Sobolev space
with a fixed order of differentiability,  then the stochastic integral inherits the same property with a related order of differentiability (Theorem \ref{thm2.1}).
This fact, together with Sobolev's embeddings and $L^p$-estimates of increments in time of the integral (Theorem \ref{thm2.3}), complete the analysis.

Section \ref{s2} is devoted to the study of  equation (\ref{0.1}) itself. The idea is to transfer the  properties of the stochastic integral obtained in Section  \ref{s1} to the solution of the equation.
First, in Section \ref{ss2.1}, we give a more general version of existence and uniqueness of a solution and its properties than in \cite{DM}, allowing non vanishing initial conditions and an additive non-linearity $b$. We also show how  the $L^q$--moments of the solution depend on properties of the initial conditions (see Theorem \ref{thm3.1}). 

Next, in Section \ref{ss2.2}, we go beyond the $L^q$--norm in the space variable. We see in Theorem \ref{thm3.6} 
how the assumptions on the initial conditions and on $\dot F$ imply that the fractional Sobolev norm in the space variable of the solution of equation (\ref{0.1}) is finite.
The analysis is carried out at the level of the s.p.d.e.~driven by the smoothed kernel $G_n$, and then transferred to the solution of 
equation (\ref{0.1}) by means of the properties on the fractional Sobolev norm of the contribution of the initial conditions (Lemma \ref{lem3.4}), 
an approximation result proved in Proposition \ref{prop3.3} and Fatou's lemma. With the Sobolev embedings, we obtain the H\"older continuity
property in the space variable of the sample paths.  

In Section \ref{ss2.3}, we prove regularity in time using the classical approach based on Kolmogorov's continuity criterion. We fix a bounded
domain $D\subset \IR^3$ and we first 
study the H\"older continuity in time, uniformly in $x\in D$, of the contribution of the initial conditions (Lemma \ref{lem3.9}). Secondly, we
find upper bounds for the $L^q(\Omega)$ norm of increments in time of the solution of (\ref{0.1}), uniformly in $x\in D$ (Theorem \ref{tcontinuity}).
We end up with the joint H\"older continuity in space and in time stated in Theorem \ref{holdercont}.
In particular, these results are non-trivial even for the deterministic inhomogeneous three-dimensional wave equation (see Remark \ref{detwave}). 

In Section \ref{s3}, we check the sharpness of the results proved in Section \ref{s2} by considering the most simple example consisting of an equation with vanishing initial conditions and coefficient $b$, and constant coefficient $\sigma$. In this case, the solution is a  stationary Gaussian process and all the information about the sample paths is contained in the covariance function. From the results of Section \ref{s2}, we already have upper bounds on $L^q$--norms of increments of the solution. We complete the analysis by obtaining
sharp lower bounds for these increments; this requires precise estimates of integrals related to (\ref{0.2}).  

The last section of the article, Section \ref{B},  gathers the somewhat technical but crucial sharp estimates on integrals of increments of  a class of  generalized Riesz potentials that also involve the smoothed version of the fundamental solution of the wave equation.

%%%%%%END INTRODUCTION

%%%%%%BEGINNING SECTION 2
\section{The fundamental solution of the wave equation and the covariance function}
\label{F}

The first part of this section is devoted to introducing the smoothing of the fundamental solution of the wave equation used throughout the paper.
We prove some of its properties as well as some of the properties of the fundamental solution itself. In the second part, 
we obtain expressions for first and second order increments of the covariance function.  Informally, these express the covariance function as a fractional integral of its fractional derivative; they are proved by applying the semigroup property of the Riesz kernels.  

\subsection{Some properties of the fundamental solution and its regularisations}
\label{A1}

  Let $d\ge 1$ and $\psi : \IR^d \to \IR_+$ be a function in $C^\infty(\IR^d)$ with support included in $B_1(0)$ and such that $\int_{\IR^d}
\psi(x)dx = 1$ ($B_r(x)$ denotes the open ball centered at $x\in \IR^d$ with radius $r \geq 0$). For any $t \in\, ]0, 1]$ and $n \geq 1$, we define
$$
   \psi_n(t,x) = \left( \frac{n}{t}\right)^d \psi \left(\frac{n}{t} x \right)
$$
and
\begin{equation}\label{A.1}
   G_n(t,x) = (\psi_n(t,\cdot) \ast G(t))(x),
\end{equation}
where $``\ast''$ denotes the convolution operation in the spatial variable. Observe that $\int_{\IR^d} \psi_n(t,x)dx = 1$ and
$$
 \mbox{supp}~G_n(t, \cdot) \subset B_{t(1+\frac{1}{n})} (0).
$$

   The following elementary scaling property plays an important role in the study of regularity properties in time of the
stochastic integral. Its proof is included for convenience of the reader.

\begin{lemma} Let $d=3$. For any $s, t \in [0,T]$ and $v_0 \in C(\IR^3)$,
\begin{equation}\label{A.2a}
   \int_{\IR^3} G(s,du)\, v_0(u) = \frac{s}{t} \int_{\IR^3} G(t,du)\, v_0(\frac{s}{t}u),
\end{equation}
and for any $x \in \IR^3,$
\begin{equation}\label{A.2}  
   G_n(t, \frac{t}{\snbar} x) = \left(\frac{\snbar}{t}\right)^2 G_n(\snbar, x).
\end{equation}
\label{lemA.1}
\end{lemma}

\noindent{\em Proof.} The first equality follows from the fact that the transformation $u \mapsto \frac{s}{t} u$ maps $G(t,\cdot)$ onto $\frac{t}{s}\,G(s,\cdot)$.

   The change of variables $y \mapsto \frac{t}{\snbar} \ y$ yields
\begin{eqnarray*}
 G_n(t, \frac{t}{\snbar} x) &=& \int_{\IR^3} G(t,dy) \psi_n\left(t, \frac{t}{\snbar} x-y\right)\\
    &=& \int_{\IR^3} G(\snbar, dy) \psi_n\left(t, \frac{t}{\snbar}(x-y)\right) \frac{t}{\snbar}\,.
\end{eqnarray*}
Since
$$
  \psi_n\left(t, \frac{t}{\snbar}(x-y)\right) = \left(\frac{n}{t}\right)^3 \psi \left(\frac{n}{\snbar}(x-y)\right)% \frac{t}{\snbar} 
  = \psi_n(\snbar, x-y)) \left(\frac{\snbar}{t}\right)^3,
$$
it follows that
\begin{eqnarray*}
G_n(t, \frac{t}{\snbar} x) &=& \int_{\IR^3} G(\snbar, dy) \psi_n(\snbar, x-y)) \left(\frac{\snbar}{t}\right)^2 \\
   &=& \left(\frac{\snbar}{t}\right)^2 G_n(\snbar, x).
\end{eqnarray*}    
%\begin{eqnarray*}
%\displaystyle\int_{\IR^3} G(\snbar, dy) \psi_n(t, \frac{t}{\snbar}(x-y)) \frac{t}{\snbar} &=&  \int_{\IR^3} G(\snbar,
%dy)\left(\frac{n}{t}\right)^3 \psi \left(\frac{n}{\snbar}(x-y)\right) \frac{t}{\snbar}\\
%   \\
%   &=& \displaystyle \int_{\IR^3} G(\snbar, dy) \left(\frac{n}{\snbar}\right)^3 \psi \left(\frac{n}{\snbar} (x-y)\right)
%\left(\frac{\snbar}{t}\right)^2 \\
%   &=& \int_{\IR^3} G(\snbar, dy) \psi_n(\snbar, x-y)) \left(\frac{\snbar}{t}\right)^2\\
%   \\
%   &=& \left(\frac{\snbar}{t}\right)^2 G_n(\snbar, x).
%\end{eqnarray*} 
%\qed
This proves the lemma.
\hfill $\Box$
\vskip 16pt

   We recall the following integrability condition of the fundamental solution of the wave equation, valid for any $\beta \in\,
]0,2[$:
\begin{equation}\label{A.3}
   \sup_{t \in[0,T]} \int_{\IR^d} \frac{\vert{\cal{F}}G(t)(\xi)\vert^2}{\vert \xi \vert^{d-\beta}}\, d\xi \leq C(1+T^2).
\end{equation}
Indeed, 
\begin{equation}\label{fourierG}
   {\cal{F}}G(t)(\xi) = \vert\xi\vert^{-1}\sin(t \vert \xi \vert)
\end{equation} 
(see \cite{treves}), and therefore
$$
\int_{\IR^d} \frac{\vert {\cal{F}} G(t)(\xi)\vert^2}{\vert\xi\vert^{d-\beta}} d \xi \leq I_1(t)+I_2(t),
$$
where
$$
   I_1(t) = \int_{\vert \xi\vert \leq 1} \frac{t^2}{\vert \xi \vert^{d-\beta}}\, d \xi \leq C_1 T^2,\, \qquad
   I_2(t) = \int_{\vert \xi \vert > 1} \frac{d \xi}{\vert \xi \vert^{d+2-\beta}}\, d\xi \leq C_2.
$$
A similar property holds for $G_n$. In fact, since $\vert{\cal{F}} \psi_n(t)(\xi)\vert \leq 1$, 
\begin{equation}\label{A.4}
  \vert {\cal{F}} G_n(t) (\xi)\vert = \vert {\cal{F}}  
   \psi_n(t)(\xi)\vert\, \vert {\cal{F}} G(t)(\xi) \vert \leq
    \vert {\cal{F}} G(t) (\xi) \vert
\end{equation}
and therefore
\begin{equation}\label{A.5}
\sup_{n \geq 1} \sup_{t \in[0,T]} \int_{\IR^d} \frac{\vert {\cal{F}} G_n(t)(\xi)\vert^2}{\vert \xi \vert^{d-\beta}}\, d \xi
\leq C(1+ T^2).
\end{equation}
The next statement gives a more precise result than (\ref{A.3}).

\begin{lemma}
For any $t \in [0,T]$ and $\beta \in\,]0,2[$,
\begin{equation}\label{A.6}
   \int_0^t ds \int_{\IR^d} \frac{\vert {\cal{F}} G(s)(\xi)\vert^2}{\vert \xi \vert^{d-\beta}} d\xi \leq C t^{3-\beta}.
\end{equation}
\label{lemA.2}
\end{lemma}

\noindent{\it Proof.} By Fubini's theorem, and using the change of variables $w= t \, \xi$, we see that
\begin{align*}
 \int^t_0 ds \int_{\IR^d} \frac{\vert {\cal{F}} G(s)(\xi)\vert^2}{\vert \xi \vert^{d-\beta}} d\xi 
 &= \int_{\IR^d} \frac{d \xi}{\vert \xi \vert^{d+2-\beta}} \int_0^t \frac{1-\cos(2s \vert \xi \vert)}{2} ds\\
 \\
 &= \int_{\IR^d} \frac{d\xi}{\vert \xi \vert^{d+2-\beta}} \left(\frac{t}{2} - \frac{\sin(2t \vert \xi \vert}{4
\vert \xi \vert}\right)\\
\\
 &=  t^{3-\beta} J,
%\int_{\IR^3} \frac{dw}{\vert w \vert^{5-\beta}} \left( \frac{1}{2} - \frac{\sin(2 \vert w\vert)}{4\vert w \vert} \right),
\end{align*}
where  
$$
   J = \int_{\IR^d} \frac{dw}{\vert w \vert^{d+2-\beta}} \left(\frac{1}{2} - \frac{\sin(2 \vert w \vert)}{4 \vert w
\vert}\right).
$$
\smallskip

\noindent Note that $J < \infty$. Indeed, $J \leq J_1 + J_2,$ where
$$
   J_1 = \int_{\vert w \vert > 1} \frac{dw}{\vert w \vert^{d+2-\beta}},
   \qquad
   J_2 = \int_{\vert w \vert \leq 1} \frac{dw}{\vert w \vert^{d+2-\beta}} \left( \frac{1}{2}- \frac{\sin(2 \vert w
\vert)}{4\vert w\vert} \right).
$$
Clearly, $J_1 < \infty$. For $J_2$, since $\sin(2\vert w \vert) = 2 \vert w \vert - \frac{2^3}{3!} \cos (\zeta) \vert w
\vert^3,$ with $\zeta \in\, ]0, \vert w \vert[$,
$$
   J_2 \leq C \int_{\vert w \vert \leq 1} \frac{dw}{\vert w \vert^{d-\beta}} < \infty.
$$
This establishes (\ref{A.6}).
\hfill $\Box$
\vskip 16pt
%\qed

\begin{lemma} For any $b > 0$ and $\beta \in\, ]0,2[$ such that $\beta + b \in\, ]0,3[$,
\begin{equation}\label{A.7}
   \sup_{t\in[0,T]}\,\int_0^t \frac{ds}{s^b} \int_{\IR^d} \frac{\vert{\cal{F}} G(s)(\xi)\vert^2}{\vert \xi\vert^{d-\beta}}\, d\xi < \infty .
\end{equation}
\label{lemA.3}
\end{lemma}

\noindent{\it Proof.} The change of variables $\xi \mapsto s \xi$ shows that the integral in (\ref{A.7}) is equal to
$$
%\int_o^t \frac{ds}{x^b} \int_{\IR^3} \frac{\vert %{\cal{F}}G(s)(\xi)\vert^2}{\vert \xi \vert^{3-\beta }}d\xi = 
\int_0^t ds \ s^{2-(\beta+b)}\int_{\IR^d} \frac{\sin^2 \vert \xi \vert}{\vert \xi \vert^{d+2-\beta}} d \xi.
$$
The inner integral 
%$\int_{\IR^3} \frac{\sin^2 \vert \xi \vert}{\vert \xi \vert^{5-\beta}} d\xi$ 
is finite for $\beta \in\, ]0,2[$. Therefore (\ref{A.7}) holds when $\beta + b < 3$.
\hfill $\Box$
\vskip 16pt
%\qed

%SUBSECTION 2.2	
\subsection{The covariance function and Riesz kernels}
\label{A2}

   We assume that the covariance measure of the noise is absolutely continuous with respect to Lebesgue measure, that is,
$\Gamma(dx) = f(x) dx.$ In addition, we suppose that $f$ satisfies the following assumption.
\vskip 12pt
\begin{assump} There is $\beta \in \, ]0,2[$ and $\delta \in \,]0,1]$ such that 
\begin{equation*} f(x) = \varphi(x) k_\beta(x),
\end{equation*}
 where
$k_\beta(x) = \vert x \vert^{-\beta}$, $x \in \IR^d\setminus \{0\}$, and $\varphi$ is bounded and positive, $\varphi \in C^1(\IR^d)$ and
$\nabla \varphi \in C^\delta_{\mbox{\scriptsize b}}(\IR^d)$ (the space of bounded and H\"older continuous functions with exponent $\delta$).
\label{assumpC}
\end{assump}

\begin{example}\label{example2.5} (a) The basic case is when $\varphi \equiv 1$. In this case, $f \equiv k_\beta$ is termed a {\em Riesz kernel}.  We recall that $k_\beta = c_{d,\beta} \F k_{d-\beta}$ \cite[Chapter V]{stein}.

   (b) Another possibility is $\varphi(x) = \exp(- \sigma^2 \vert x \vert^2/2)$. In this case, $f(x)$ is indeed a covariance function, since $f = \F(k_{d-\beta}\ast \psi)$, where $\psi(\xi) = (2\pi \sigma^2)^{-3/2}\exp(-\vert \xi \vert^2/(2 \sigma^2))$. The parameter $\delta$ in Assumption \ref{assumpC} can be set to $1$. Condition (\ref{0.2}) is satisfied since $\beta \in \, ]0,2]$. Indeed,
\begin{align*}
   \int_{\IR^d} d\xi\, (k_{d-\beta}\ast \psi)(\xi)\, \vert \F G(t)(\xi)\vert^2 &=
      \int_{\IR^d} dz \, \psi(z) \int_{\IR^d} d\xi\, k_{d-\beta}(\xi)\, \vert \F G(t)(\xi + \eta)\vert^2 \\
      &\leq \sup_{\eta \in \IR^d} \int_{\IR^d} d\xi\, k_{d-\beta}(\xi)\, \vert \F G(t)(\xi + \eta)\vert^2 ,
\end{align*}
and the right-hand side is finite when $\beta \in \, ]0,2]$: see \cite[Lemma 8]{DM}.
\end{example}

   The  {\em Riesz potentials} $I_a$ associated with the function $k_\beta(x)$ are defined by
$$
   (I_a \varphi)(x) = \frac{1}{\gamma(a)} \int_{\IR^d} \vert x-y \vert^{-d+a}\varphi(y)dy,
$$
for $\varphi \in \IS(\IR^d)$, $a \in\, ]0,d[$ and $\gamma(a) = \pi^{d/2} 2^{a}
\Gamma(\frac{a}{2})/\Gamma(\frac{d-a}{2}).$ Riesz potentials can be interpreted as fractional integrals and have the semigroup property
$$
I_{a+b}\, \varphi = I_a(I_b \varphi),\ \varphi \in \IS(\IR^d),\ a+b \in\, ]0,d[
$$
(see \cite[p.118]{stein}). This property implies in particular that
\begin{equation}\label{B.1}
   \vert x-y \vert^{-d+(a+b)} = \int_{\IR^d} dz\, \vert x-z \vert^{-d+a}\vert z-y \vert^{-d+b},
\end{equation}
provided $a+b \in\,]0,d[$. This equality can be informally interpreted by saying that $\vert \cdot \vert^{-d+(a+b)}$ is the fractional integral of order $a$ of $\vert \cdot \vert^{-d+b}$, which is natural since $(-\Delta)^{a/2}(\vert \cdot \vert^{-d+(a+b)}) = \vert \cdot \vert^{-d+b}$ as can be checked by taking Fourier transforms.
\smallskip

   We will make heavy use of properties of first and second order increments of Riesz kernels. For a function $f : \IR^d \to \IR$, we
set 
\begin{align}
 Df(u,x) & = f(u+x) - f(u),\nonumber\\
 D^2f(u,x) & = f(u-x) -2 f(u) + f(u+x),\nonumber\\
 \bar D^2f(u,x,y) & = f(u+x+y)-f(u+x)-f(u+y)+f(u).
 \label{B.1'}
\end{align}
Notice that $D^2f(u,x) =  \bar D^2f(u-x,x,x)$.
\smallskip

\begin{lemma} Fix $u, x, y \in \IR^d, a+b \in\, ]0,d[$. 
%and set $e=\frac{x}{|x|}$.
The following properties hold:
\begin{description}
\item{(a)} For any
$c \in\IR$, 
$$
   Dk_{d-a-b}(u,cx) = |c|^b \int_{\IR^d}dw\, k_{d-a}(u- c\,w) Dk_{d-b}(w,x).
$$
\item{(b)} For any $b \in\, ]0,1[$ and any vector $e\in\IR^d$ with $|e|=1$,
$$
   \int_{\IR^d} dw \, \vert Dk_{d-b}(w,e)\vert < \infty.
$$ 
\item{(c)} Set $e=\frac{x}{|x|}$, $x\in\IR^d$. Then
$$
   \vert D^2k_{d-a-b}(u,x) \vert \leq \vert x \vert^b \int_{\IR^d}dw \, k_{d-a}(u- \vert x \vert w) \vert D^2k_{d-b}(w,e)
\vert.
$$
\item{(d)} For any $b \in\, ]0,2[$ and each $e\in\IR^d$ with $|e|=1$,
$$
   \int_{\IR^d} dw\, \vert D^2k_{d-b}(w,e) \vert < \infty.
$$ 
\item{(e)}
$$
\bar D^2k_{d-a-b}(u,cx,cy)\le |c|^b \int_{\IR^d} dw\, k_{d-a}(u - c\,w)|\bar D^2k_{d-b}(w,x,y)|.
$$
\end{description}
\label{lemB.1}
\end{lemma}

\noindent{\it Proof.}  (a) From (\ref{B.1}), we obtain
$$
   Dk_{d-a-b}(u,cx) = \int_{\IR^d} dz\, k_{d-a}(u-z) Dk_{d-b}(z,cx).
$$
Set $z = c w$ to see that this is equal to
\begin{align*}
&|c|^d \int_{\IR^d} dw\, k_{d-a}(u-c\,w) Dk_{d-b}(c\,w, cx) \\
&\quad =  |c|^b \int_{\IR^d} dw\, k_{d-a}(u- c\, w) Dk_{d-b}(w, x),
\end{align*}
which proves (a).

   We now check (b). Set 
$$
   I = \int_{\IR^d} dw\, \vert Dk_{d-b}(w,e) \vert
$$ 
and consider the decomposition $I = I_1 + I_2$, where $I_1$ (resp. $I_2$) is the integral of the same expression but over
$B_2(0)$ (resp. $B_2(0)^c$) instead of $\IR^d$. Then
$$
   I_1 \leq 2 \int_{B_3(0)} dw\, \vert w \vert^{-d+b} \leq C \int_0^3 \rho^{b-1} d \rho < \infty,
$$
if $b > 0.$ As for $I_2$, we write
\begin{align*}
I_2 &=\int_{\vert w \vert > 2} dw\, \vert \int_0^1 d \lambda \frac{d}{d \lambda} ( \vert w-\lambda e\vert^{-d+b} )\\
&\leq \int^1_0 d \lambda \int_{\vert w \vert > 2} dw \vert w-\lambda e \vert^{-d+b-1}\\
& \leq C \int^1_0 d \lambda \int_1^\infty \rho^{b-2} d \rho,
\end{align*}
which is finite if $b < 1$.

   The proof of (c) is analogous to that of (a). It suffices to apply the identity (\ref{B.1}) twice, since $D^2k_{d-a-b}(u,x) =
Dk_{d-a-b}(u,x) + Dk_{d-a-b}(u,-x)$, and to use the change of variables $z = \vert x \vert w.$ 

   Let us now prove (d). Observe that $D^2k_{d-b}(\omega,e) = k_{d-b}(w\pm e) - 2k_{d-b}(w)$ and that the integrals
\begin{equation*}  
\int_{\vert w \vert \leq 2} dw\, \vert w \vert^{-d+b}, \quad \int_{\vert w \vert \leq 2}dw\,
\vert e\pm w \vert^{-d+b},
\end{equation*}
converge for each $b >0$.

We next study
$\int_{\vert w\vert > 2} dw\, \vert D^2k_{d-b}(w,e)\vert$. Set
%$$
%\varphi_-(\lambda) = \vert w-\lambda e \vert^{-d+b},\qquad
%\varphi_+(\lambda) = \vert w + \lambda e \vert^{-d+b}, \qquad \lambda \in [0,1].
%$$
\begin{equation*}
\Phi (\lambda,\mu) = k_{d-b}\left(w-(\lambda - \mu)e\right), \, \lambda, \mu \in [0,1].
\end{equation*}
Then 
$$
   D^2k_{d-b}(w,e) =  \int_0^1 \,d\lambda \, \int_0^1 \,d\mu \, \frac{\partial^2\Phi}{\partial\lambda\partial\mu}(\lambda,\mu).
$$ 
Elementary computations lead to
%$$
%\varphi_\pm^\prime(\lambda) = \pm(b-d) \vert w \pm \lambda e \vert^{-d+b-2} \displaystyle\sum^d_{i=1} e_i (w_i \pm
%\lambda e_i).
%$$
\begin{equation*}
\left\vert \frac{\partial^2\Phi}{\partial\lambda\partial\mu}(\lambda,\mu)\right\vert \le C\,k_{d-b+2}\left(w-(\lambda-\mu)e\right).
\end{equation*}
%Thus,
%\begin{align*}
%   D^2k_{d-b}(w,e) &= -\int^1_0 d \lambda \sum^d_{i=1} (d-b) \ e_i \Big((w_i-\lambda e_i) \vert w-\lambda
%e\vert^{-d+b-2}\\
%&\quad -(w_i+\lambda e_i) \vert w + \lambda e\vert^{-d+b-2}).
%\end{align*}
%Set $F_i(\mu) = (w_i - \mu \lambda e_i) \vert w-\mu \lambda e \vert^{-d+b-2}$, $\mu \in [-1,1].$ It is easy to check that
%$$
 %  \vert F_i^\prime (\mu)\vert \le C \vert w- \mu \lambda e \vert^{-d+b-2}.
%$$
%Therefore,
%
%\begin{equation*}
 %  D^2k_{d-b}(w,e) \le  C \int_0^1 d \lambda \int^1_{-1} d\mu\, \vert w - \mu  \lambda e   \vert^{-d+b-2}.
%\end{equation*}

Therefore, by Fubini's theorem,
$$
   \int_{\vert w \vert > 2} dw\, \vert D^2 k_{d-b}(w,e)\vert \leq C \int^1_0 d\lambda \int^1_{0} d \mu \int_{\vert w
\vert > 2} dw\, \vert w-(\lambda-\mu)e \vert^{-d+b-2}.
$$
The integral $\int_{\vert w\vert > 2} dw\, \vert w-(\lambda-\mu)e\vert^{-d+b-2}$ converges for any $b < 2$.  Consequently, (d) is proved.

The proof of (e) is analogous to that of (a). It suffices to apply the identity (\ref{B.1}) four times,
 and to use the change of variables $z = c w.$ 
\hfill $\Box$
\vskip 16pt
%\qed

   Lemma \ref{lemB.1} is the basis for the main technical estimates of this paper, whose statements and proofs are deferred to Section \ref{B}.

%%%%%%%%%%%%%%%%%%%%%%%%%%%%%%%%
%%%%%% BEGINNING SECTION 3
\section{{H\"{o}lder-Sobolev regularity of the stochastic integral}}
\label{s1}

In this section, we consider the extension of the stochastic integral given in Section 2 of \cite{DM} in the particular
case where $G$ is the fundamental solution of the wave equation in spatial dimension $d=3$,  defined in (\ref{G}). More precisely, let $\{Z(s),\ s
\in [0,T]\}$ be an ${\cal{F}}_s$-measurable, $L^2(\IR^3)$-valued stochastic process that is mean-square continuous. Let $F$ be the Gaussian process with covariance measure $\Gamma$ and covariance functional as in (\ref{cov1}). Assume that its spectral measure $\mu=\mathcal{F}^{-1}\Gamma$ satisfies condition (\ref{0.2}).

 As it has been proved in \cite{dalangfrangos}, the process $F$ can be extended to a worthy martingale measure in the sense given in \cite{walsh}. We shall denote by $M = (M_t(A),\ t\ge 0,\ A \in \cB_b(\IR^3))$ this extension. The relationship between $F$ and $M$ is 
$$
   F(\psi) = \int_{\IR_+ \times \IR^3} \psi(t,x) M(dt,dx),
$$
for all $\psi \in \D(\IR^{d+1})$. The covariance measure $Q$ and dominating measure $K$ of $M$ are 
$$
   Q(A \times B \times \, ]s,t]) = (t-s) \int_{\IR^3} dx \int_{\IR^3} dy\, 1_A(x) f(x-y) 1_B(y)
$$
and $K \equiv Q$. 

   Then, following Theorem 6 in \cite{DM}, the stochastic integral
$$
   v^t_{G,Z} = \int^t_0 \int_{\IR^3} G(t-s, \cdot -y) Z(s, y) M(ds, dy)
$$
is well-defined as a random vector in $L^2(\Omega; L^2(\IR^3))$ and has the isometry property
\begin{equation}\label{isometry1}
   E\left(\Vert v_{G,Z}^t\Vert^2_{L^2(\IR^3)}\right) = \int^t_0 ds \int_{\IR^3} d\xi E\left(\vert {\cal{F}} Z(s)(\xi)\vert^2\right)
\int_{\IR^3} \mu(d\eta)\, \vert {\cal{F}} G(t-s) (\xi-\eta) \vert^2.
\end{equation}
It also satisfies the bound
$$
   E\left(\Vert v_{G,Z}^t\Vert^2_{L^2(\IR^3)}\right) \leq
   \int^t_0 ds\, E\left(\Vert Z(s)\Vert^2_{L^2(\IR^3)} \right) \sup_{\xi \in \IR^3} \int_{\IR^3} \mu(d\eta)\, \vert {\cal{F}} G(t-s) (\xi-\eta) \vert^2.
$$

   Let $\cO$ be a bounded or unbounded open subset of $\IR^3$, $p \in [1, \infty[$, $\gamma \in \, ]0,1[$. Define
$$
    \Vert g \Vert_{\gamma,p,\cO} = \left( \int_{\cO} dx \int_{\cO} dy\, \frac{\vert g(x)-g(y)\vert^p}{\vert x- y
\vert^{3+\gamma p}} \right)^{1/p}.
$$
When $\cO = \IR^3$, we write $\Vert g \Vert_{\gamma,p}$ instead of $\Vert g \Vert_{\gamma,p,\IR^3}$.

   We denote by $W^{\gamma, p}(\cO)$ the Banach space consisting of functions $g: \IR^3 \to \IR$ such that
\begin{equation}\label{eq2.1}
   \Vert g \Vert_{W^{\gamma, p}(\cO)} : = \Vert g \Vert_{L^p(\cO)} +  \Vert g \Vert_{\gamma,p,\cO}
\end{equation}
is finite. The spaces $W^{\gamma, p}(\cO)$ are the fractional Sobolev spaces  
(see for instance \cite{adams}, \cite{shima}). 

   Given a bounded set $K\subset {\IR}^3$ and $\ep>0$, we let $K^{\ep}$ be the open set
$$
   K^{\ep}=\{x\in {\IR}^3: \exists z \in K \mbox{ such that } \vert x - z \vert < \ep\}.
   $$ 

\subsection{Regularity in space}

The following result concerns the Sobolev regularity of the stochastic integral in the space variable. We assume that $\Gamma(dx)$ is of the form $\Gamma(dx) = f(x) dx$, where $f(x)$ is as in (\ref{rde1}), and there is $\beta \in \, ]0,2[$ and $\delta \in \, ]0,1[$ such that Assumption \ref{assumpC} is satisfied.

\begin{thm} 
\label{thm2.1}
Fix $T>0$, $q \in \, ]3,\infty[$ and a bounded domain $\cO\subset \IR^3$. Suppose that $\tau_q(\beta,\delta):= (\frac{2-\beta}{2} \wedge \frac{1+\delta}{2}) - \frac{3}{q} > 0$ and fix $\gamma \in \, ]0,1[$, $\rho \in \, ]0,\tau_q(\beta,\delta) \wedge\gamma[$.  Consider a stochastic process $Z$  such that for some fixed $t \in [0,T]$,
\begin{equation}\label{2.2}
   \int_0^t ds\, E\left(\Vert Z(s)\Vert^q_{W^{\gamma,q}(\cO^{t-s})}\right) < \infty.
\end{equation}
There is $C < \infty$ (depending on $\rho$ but not on $Z$) such that
\begin{equation}\label{2.2.1}
E\left(\Vert v^t_{G,Z} \Vert^{\qnbar}_{W^{\rho, \qnbar}(\cO)}\right) \le C 
\int_0^t ds\, E\left(\Vert Z(s)\Vert^q_{W^{\rho,q}(\cO^{t-s})}\right).
\end{equation}
\end{thm}
\medskip

The main ingredient in the proof of this theorem is the next lemma, which will also be used in the proof of Theorem \ref{thm3.6} in the next section.

%%%%%%%%%%%%%%%%%%Crucial lemma
\begin{lemma}
\label{l2.2}
Fix $q\in\, ]3,\infty[$ such that $\tau_q(\beta,\delta)> 0$. Fix $\ro\in\,]0,\tau_q(\beta,\delta)[$, a compact set $K \subset \IR^3$ 
and a bounded domain $\cO\subset K$.
Let $G_n$ be given in (\ref{A.1}). 
For any $\alpha\in\, ]0,(2-\beta)\wedge 1[$ satisfying $\frac{\alpha}{2}-\frac{3}{q} < \ro$,  there is a finite constant $C=C(T,\ro,\alpha,K)$ such that for every $t\in[0,T]$
and $n\ge 1$,
\begin{align}
&\,E\left(\Vert v^t_{G_n,Z}\Vert^q_{\ro,q,\cO}\right)\le C\int_0^t ds \Big( E\big(\Vert Z(s)\Vert^q_{W^{\ro,q}(\cO^{(t-s)(1+\frac{1}{n})})}\big)\nonumber\\
&\quad +\big(E\big(\Vert Z(s)\Vert^q_{L^q(\cO^{(t-s)(1+\frac{1}{n})})}\big) E\big(\Vert Z(s)\Vert^q_{2\ro-\alpha+\frac{3}{q},q,\cO^{(t-s)(1+\frac{1}{n})}}\big)\big)^{\frac{1}{2}}\Big).\label{222}
\end{align}
In addition,
\begin{equation}
\label{222'}
E\left(\Vert v^t_{G_n,Z}\Vert^q_{\ro,q,\cO}\right)\le C\int_0^t ds E\big(\Vert Z(s)\Vert^q_{W^{\ro,q}(\cO^{(t-s)(1+\frac{1}{n})})}\big).
\end{equation}
\end{lemma}

\noindent{\it Proof}:  Set $p=q/2$ and fix $t\in\, ]0,T]$. Note that $v^t_{G_n,Z}(x)$ is a Walsh stochastic integral and
\begin{align*}
   & v^t_{G_n,Z}(x) - v^t_{G_n,Z}(y)\\
   &\qquad\qquad = \int_0^t ds \int_{\IR^3} (G_n(t-s, x-u) - G_n(t-s, y-u))
         Z(s,u) M(ds,du).
\end{align*}
Burkholder's inequality yields
\begin{align*}
&E\big(\vert v^t_{G_n,Z}(x) - v^t_{G_n,Z}(y)\vert^{q}\big) \\
&\quad \leq C E\Big(\vert \int_0^t ds \int_{\IR^3} du \int_{\IR^3} dv\, Z(s, u)
f(u-v) Z(s,v)\\
& \quad\quad\times (G_n(t-s, x-u) - G_n(t-s, y-u)) \\
& \quad\quad\times (G_n(t-s, x-v) - G_n(t-s, y-v))\vert^p\Big)\\
& \quad = C E\Big(\vert I^t_{n}(x,x) - I_{n}^t(x,y)-I_{n}^t(y,x) + I^t_{n}(y,y)\vert\Big)^p,
\end{align*}
with
$$
I^t_{n}(x,y) = \int_0^t ds \int_{\IR^3} du \int_{\IR^3} dv\, g_n(t,s,x,y,u,v) Z(s, u) f(u-v) Z(s,v)
$$
and
$$
   g_n(t,s,x,y,u,v) = G_n(t-s, x-u) G_n(t-s, y-v).
$$
After the change of variables $s \mapsto t-s$ and the spatial transformation %$(u,v) \mapsto (x-u, x-v)$, 
$(u,v) \mapsto (x-u, y-v)$, %$(u,v) \mapsto (y-u, x-v)$, $(u,v) \mapsto (y-u, y-v)$, respectively,
 we obtain
$$
   I^t_{n}(x,y) = \int^t_0 ds \int_{\IR^3} du \int_{\IR^3} dv\, G_n(s,u)G_n(s,v) h(t,s,x,y,u,v),
$$
where
$$
    h(t,s,x,y,u,v) = Z(t-s, x-u) f(x-y-u+v) Z(t-s, y-v).
$$
Using these expressions and rearranging terms, it is straightforward to check that
$$
   I_{n}^t(x,x) - I_{n}^t(x,y) - I_{n}^t(y,x) + I_{n}^t(y,y) = \sum^4_{i=1}  J_{i,n}^t (x,y) ,
$$
where, for $i=1,\dots,4$,
$$
    J_{i,n}^t(x,y) = \int^t_0 ds \int_{\IR^3} du \int_{\IR^3}dv \, G_n(s,u)G_n(s,v) h_i(t,s,x,y,u,v),
$$
and, using the notation in (\ref{B.1'}), 
\begin{align} \nonumber
   h_1(t,s,x,y,u,v) &=f(y-x+v-u)(Z(t-s, x-u) - Z(t-s, y-u)) \\ \nonumber
   &\qquad\times  (Z(t-s, x-v) - Z(t-s,y-v)),\\ \nonumber
   \\ \nonumber
   h_{2}(t,s,x,y,u,v) &= Df(v-u,x-y) Z(t-s, x-u)\\ \nonumber
   & \qquad\times  (Z(t-s, x-v)-Z(t-s, y-v)),\\ \nonumber
   \\ \nonumber
   h_{3}(t,s,x,y,u,v) &= Df(v-u, y-x) Z(t-s,y-v) \\ \nonumber
   & \qquad\times  (Z(t-s, x-u) - Z(t-s, y-u)),\\ \nonumber
   \\
     h_{4}(t,s,x,y,u,v) &= -D^2 f(v-u,x-y) Z(t-s,x-u) Z(t-s, x-v).
     \label{starp13}
\end{align}

%%%%%%%%%%%
Consequently, setting $\bro=\ro+\frac{3}{q}$ and 
$$T_n(t,\cO)=E\left(\Vert v^t_{G_n,Z}\Vert^q_{\ro,q,\cO}\right),$$
we have
$$T_n(t,\cO) \leq C \sum^4_{i=1} T^i_n(t,\cO),$$
where 
$$
   T^i_n(t,\cO) = \int_{\cO} dx \int_{\cO} dy\, \frac{E(\vert J^t_{i,n}(x,y) \vert^p)}{\vert x-y \vert^{q\bro}}.
$$

   Set
$$
   \mu^1_n(x,y) = \sup_{s \in [0,T]} \int_{\IR^3} du \int_{\IR^3} dv\, G_n(s,u)G_n(s,v)f(y-x+v-u).
$$
Notice that, since $f(x) = \varphi(x) k_\beta(x)$ with $\varphi$ bounded (see Assumption \ref{assumpC}) and $G_n \geq 0$, 
\begin{equation*}
   \sup_{n, x, y} \mu_n^1(x,y) \leq \int_{\IR^3} \frac{ \vert {\cal{F}} G(s)(\xi)\vert^2}{\vert \xi \vert^{3-\beta}} d \xi <
\infty,
\end{equation*}
for any $\beta \in\,]0,2[$ (see (\ref{A.3})).

 Therefore, H\"{o}lder's inequality 
%with respect to the measure $dudv G_n(s,u)G_n(s,v)f(y-x+v-u)$ 
implies that
\begin{align*}
 E\left(\vert J^t_{1,n}(x,y)\vert^p\right) &\leq \left(T\mu^1_n(x,y)\right)^{p-1}\\
 &\quad\times  E\Big(\int_0^t ds \int_{\IR^3} du  \int_{\IR^3} dv\,
G_n(s,u)G_n(s,v)f(y-x+v-u)\\
&\quad\times \vert Z(t-s,x-u)-Z(t-s,y-u)\vert^p\, \\
&\quad\times\vert Z(t-s,x-v)-Z(t-s,y-v)\vert^p\Big).
\end{align*} 

The function $(s,x,y,u,v)\mapsto G_n(s,u)G_n(s,v)f(y-x+v-u)$ is the density of a finite measure on $[0,T]\times\cO\times \cO\times{\IR}^3\times {\IR}^3$. Indeed,
\begin{align*}
&\sup_{s\in[0,T]}\int_{\cO} dx \int_{\cO} dy \int_{{\IR}^3}du  \int_{{\IR}^3}dv\, G_n(s,u)G_n(s,v)f(y-x+v-u)\\
&\quad \le \vert\cO\vert^2  \sup_{n, x, y} \mu_n^1(x,y) \le C,
\end{align*}
where $\vert\cO\vert$ denotes the volume of $\cO$. 

Applying the Cauchy-Schwarz inequality with respect to this measure, we see that $T_n^1(t,\cO)\le C \left(T_n^{1,1}(t,\cO)\, T_n^{1,2}(t,\cO)\right)^{\half}$, with
\begin{align*}
T_n^{1,1}(t,\cO)&=\int_0^t ds \int_{\cO} dx \int_{\cO} dy \int_{{\IR}^3}du  \int_{{\IR}^3}dv\, G_n(s,u)G_n(s,v)f(y-x+v-u)\\
\\
&\quad \times \frac{E\left(\vert Z(t-s,x-u)-Z(t-s,y-u)\vert^{2p}\right)}{\vert x-y\vert^{q\bro}},\\
\\
T_n^{1,2}(t,\cO)&=\int_0^tds \int_{\cO} dx \int_{\cO} dy \int_{{\IR}^3}du  \int_{{\IR}^3}dv\,  G_n(s,u)G_n(s,v)f(y-x+v-u)\\
\\
&\quad \times \frac{E\left(\vert Z(t-s,x-v)-Z(t-s,y-v)\vert^{2p}\right)}{\vert x-y\vert^{q\bro}}.
\end{align*}
Consider the change of variables defined by $\bar x=x-u$, $\bar y=y-u$. We notice that for any $s\in[0,T]$ and $n\ge 1$, the support of the function $G_n(s)$ is included in the
ball centered at zero with radius $s \kn$, where $\kn = 1+\frac{1}{n}$. Therefore, the domain of the new variables $\bar x$, $\bar y$ is included in $\cO^{s \kn}$. Hence,
\begin{align*}
T_n^{1,1}(t,\cO)&\le \int_0^t ds \int_{\cO^{s\kn}} d\bar x \int_{\cO^{s\kn}} d\bar y\, E\left(\left\vert\frac{Z(t-s,\bar x)-Z(t-s,\bar y)}{\vert\bar x-\bar y\vert^{\bro}}\right\vert^q\right)\\
\\
&\quad\times  \int_{{\IR}^3}du  \int_{{\IR}^3}dv\,  G_n(s,u)G_n(s,v)f(\bar y-\bar x+v-u)\\
\\
&\le \sup_{n,\bar x,\bar y} \mu_n^1(\bar x,\bar y)\\
\\
&\quad\times\int_0^t ds \int_{\cO^{s\kn}} d\bar x \int_{\cO^{s\kn}} d\bar y\,E\left(\left\vert\frac{Z(t-s,\bar x)-Z(t-s,\bar y)}{\vert\bar x-\bar y\vert^{\bro}}\right\vert^q\right)\\
\\
&\le C \int_0^t ds\, E\left(\Vert Z(s)\Vert^q_{\ro,q,\cO^{(t-s)\kn}}\right).
\end{align*}

   We deal with the term $T_n^{1,2}(t,\cO)$ similarly, using the change of variables $\bar x=x-v$, $\bar y=y-v$. We obtain the same upper bound as for $T_n^{1,1}(t,\cO)$. Summarising, 
\begin{equation}
T_n^1(t,\cO)\le C \int_0^t ds\, E\left(\Vert Z(s)\Vert^q_{\ro,q,\cO^{(t-s)(1+\frac{1}{n})}}\right).
\label{223}
\end{equation}

Set 
$$
   \mu^2_n(x,y) =  \sup_{s \in [0,T]}\int_{\IR^3}du \int_{\IR^3}dv\, G_n(s,u) G_n(s,v) \vert Df(v-u,x-y)\vert.
$$
By Lemma \ref{lemB.2},
$$
   \sup_{n \geq 1}  \mu_n^2(x,y) \leq C \vert x-y \vert^\alpha,
$$
for any $\alpha\in\, ]0, (2-\beta) \wedge 1[$.

   By H\"{o}lder's inequality, for any $\alpha \in\, ]0, (2-\beta) \wedge 1[$,
\begin{align*}
& \frac{E( \vert J^t_{2,n}(x,y) \vert^p)}{\vert x-y \vert^{q\bro}}
 \leq \left(\sup_{n, x, y}  \frac{T \mu^2_n(x,y)}{ \vert x-y \vert^\alpha}\right)^{p-1}\nonumber\\ 
 \nonumber \\
  & \quad \quad \quad\times \int_0^t ds \int_{\IR^3} du \int_{\IR^3} dv\, G_n(s,u)G_n(s,v)\, \frac{\vert Df(v-u,x-y) \vert}{\vert x-y\vert^{\alpha}} \nonumber\\
 \nonumber \\
& \quad \quad \quad\times E\left( \vert Z(t-s, x-u)\vert^p \frac{\vert Z(t-s, x-v) - Z(t-s, y-v) \vert^p}{\vert x-y \vert^{p(2
\bro-\alpha)}} \right). 
\end{align*}

%Hence,
%\begin{align*}
%T_n^2(t,\cO)&\le C \int_0^t ds \int_{\cO} dx \int_{\cO} dy \int_{\IR^3} du \int_{\IR^3} dv\, G_n(s,u)G_n(s,v)\\
%\\
%&\quad \times\frac{\vert Df(v-u,x-y) \vert}{\vert x-y\vert^{\alpha}}\\
%\\
%&\quad\quad\times E\left(\vert Z(t-s,x-u)\vert^p \frac{\vert Z(t-s,x-v)-Z(t-s,y-v)\vert^p}{\vert x-y\vert^{p(2\bro-\alpha)}}\right).
%\end{align*}
For each $s\in[0,T]$, the function $(x,y,u,v)\mapsto G_n(s,u)G_n(s,v)\,\frac{\vert Df(v-u,x-y) \vert}{\vert x-y\vert^{\alpha}}$ is the density of a finite measure on
the set $\cO\times \cO\times\IR^3\times\IR^3$. In fact,
\begin{align*}
&\int_{\cO} dx \int_{\cO} dy \int_{\IR^3} du \int_{\IR^3} dv\, G_n(s,u)G_n(s,v)\,\frac{\vert Df(v-u,x-y) \vert}{\vert x-y\vert^{\alpha}}\\
\\
&\quad \le \vert\cO\vert^2 \sup_{n, x, y}  \frac{\mu^2_n(x,y)}{ \vert x-y \vert^\alpha}\le C.
\end{align*}
By the Cauchy-Schwartz inequality,
\begin{equation*}
T_n^2(t,\cO)\le C \int_0^t ds \left(T_n^{2,1}(s,\cO)\, T_n^{2,2}(s,\cO)\right)^{\half},
\end{equation*}
with the terms $T_n^{2,1}(s,\cO)$ and $T_n^{2,2}(s,\cO)$ defined as follows:
\begin{align*}
T_n^{2,1}(s,\cO)&=\int_{\cO} dx \int_{\cO} dy \int_{\IR^3} du \int_{\IR^3} dv\, G_n(s,u)G_n(s,v)\,\frac{\vert Df(v-u,x-y) \vert}{\vert x-y\vert^{\alpha}}\\
&\quad\times E\left(\vert Z(t-s,x-u)\vert^{q}\right),\\
\\
T_n^{2,2}(s,\cO)&=\int_{\cO} dx \int_{\cO} dy \int_{\IR^3} du \int_{\IR^3} dv\, G_n(s,u)G_n(s,v)\,\frac{\vert Df(v-u,x-y) \vert}{\vert x-y\vert^{\alpha}}\\
\\
&\quad \times E\left(\left\vert\frac{ Z(t-s,x-v)-Z(t-s,y-v)}{\vert x-y\vert^{2\bro-\alpha}}\right\vert^q\right).
\end{align*}
Introducing the change of variables $\bar x=x-u$, $\bar y= y-u$ yields
\begin{align*}
T_n^{2,1}(s,\cO)&\le \int_{\cO^{s\kn}}d\bar x \int_{\cO^{s\kn}}d\bar y\int_{\IR^3} du \int_{\IR^3} dv\, G_n(s,u)G_n(s,v)\\
\\
&\quad \times \frac{\vert Df(v-u,\bar x-\bar y) \vert}{\vert \bar x-\bar y\vert^{\alpha}}\,E\left(\vert Z(t-s,\bar x)\vert^{q}\right)\\
\\
&\le C\sup_{n,\bar x,\bar y} \left(\frac{\mu_n^2(x,y)}{\vert\bar x-\bar y\vert^{\alpha}}\right)\, E\left(\Vert Z(t-s)\Vert^q_{L^q(\cO^{s\kn})}\right).
\end{align*}

For the term $T_n^{2,2}(s,\cO)$, we consider the new variables $\bar x=x-v$, $\bar y= y-v$. We obtain
\begin{align*}
T_n^{2,2}(s,\cO)&\le  \int_{\cO^{s\kn}}d\bar x \int_{\cO^{s\kn}} d\bar y\int_{\IR^3} du \int_{\IR^3} dv\, G_n(s,u)G_n(s,v)\\
\\
&\quad\times \frac{\vert Df(v-u,\bar x-\bar y) \vert}{\vert \bar x-\bar y\vert^{\alpha}}\,E\left(\left\vert\frac{ Z(t-s,\bar x)-Z(t-s,\bar y)}{\vert \bar x-\bar y\vert^{2\bro-\alpha}}\right\vert^q\right)\\
\\
&\le C\sup_{n,\bar x,\bar y} \left(\frac{\mu_n^2(x,y)}{\vert\bar x-\bar y\vert^{\alpha}}\right)\, E\left(\Vert Z(t-s)\Vert^q_{2\bro-\alpha-\frac{3}{q},q,\cO^{s\kn}}\right).
\end{align*}
Hence, we have obtained
\begin{align}
\label{224}
T_n^2(t,\cO)&\le C \int_0^t ds \Big(E\big(\Vert Z(s)\Vert^q_ {L^q(\cO^{(t-s)(1+\frac{1}{n})})}\big)\nonumber\\
& \qquad\qquad \times E\big(\Vert Z(s)\Vert^q_{2\bro-\alpha-\frac{3}{q},q,\cO^{(t-s)(1+\frac{1}{n})}}\big)\Big)^{\half}.
\end{align}
With the same arguments, one can find an identical upper bound for the term $T_n^3(t,\cO)$. 

   We continue with the analysis of $T^4_n(t,\cO).$ Set
$$
   \mu^4_n(x,y) = \sup_{s \in [0,T]}\int_{\IR^3} du \int_{\IR^3} dv\, G_n(s,u)G_n(s,v)\, \vert
   D^2 f(v-u,x-y)\vert,
$$
where $D^2 f$ is defined in (\ref{B.1'}). By Lemma \ref{lemB3}, if $\bro\in\, ]0, \frac{2-\beta}{2}\wedge\frac{1+\delta}{2}[$,
\begin{equation*}
\sup_{n\ge 1}\mu^4_n(x,y)\le C\vert x-y\vert^{2\bro}.
\end{equation*}
Then, H\"{o}lder's inequality yields
\begin{align*}
 \frac{E(\vert J^t_{4,n}(x,y) \vert^p)}{\vert x-y \vert^{2p\bro}} &\leq
  \left(\sup_{n,x,y} \frac{T\mu^4_n(x,y)}{\vert x-y \vert^{2\bro}}\right)^{p-1} \nonumber\\
  \nonumber\\
&\quad \times \int_0^t ds \int_{\IR^3} du
\int_{\IR^3} dv\, G_n(s,u)G_n(s,v)\nonumber\\ 
\nonumber\\
 & \quad\times \frac{\vert D^2 f(v-u,x-y) \vert}{\vert x-y \vert^{2\bro}}\nonumber\\
\nonumber \\
 &\quad\times E( \vert Z(t-s, x-u) \vert^p\, \vert
Z(t-s,x-v)\vert^p),
\label{2.9a}
\end{align*}
Therefore,
\begin{align*}
T^4_n(t,\cO) &\leq C  \int_0^t ds  \int_{\cO} dx \int_{\cO} dy \int_{\IR^3} du \int_{\IR^3} dv\, G_n(s,u)G_n(s,v)\\
\\
&\quad \times \frac{\vert D^2 f(v-u,x-y) \vert}{\vert x-y \vert^{2\bro}}\,E( \vert Z(t-s, x-u) \vert^p\, \vert Z(t-s,x-v)\vert^p).
\end{align*}
Since the function $(s,x,y,u,v)\mapsto G_n(s,u)G_n(s,v)\frac{\vert D^2 f(v-u,x-y) \vert}{\vert x-y \vert^{2\bro}}$ is the density of a finite measure on the set
$[0,T]\times \cO\times\cO\times\IR^3\times\IR^3$, we can apply the Cauchy-Schwarz inequality with respect to this measure to obtain
\begin{equation*}
T^4_n(t,\cO)\le C \left(T^{4,1}_n(t,\cO) T^{4,2}_n(t,\cO)\right)^{\half},
\end{equation*}
with
\begin{align*}
T^{4,1}_n(t,\cO)&= \int_0^t ds \int_{\cO} dx \int_{\cO} dy \int_{\IR^3} du \int_{\IR^3} dv\, G_n(s,u)G_n(s,v)\\
\\
&\quad\times \frac{\vert D^2 f(v-u,x-y) \vert}{\vert x-y \vert^{2\bro}}\,E( \vert Z(t-s, x-u) \vert^{q}),\\
\\
T^{4,2}_n(t,\cO)&= \int_0^t ds \int_{\cO} dx \int_{\cO} dy  \int_{\IR^3} du \int_{\IR^3} dv\, G_n(s,u)G_n(s,v)\,\\
&\quad\times\frac{\vert D^2 f(v-u,x-y) \vert}{\vert x-y \vert^{2\bro}}
\, E( \vert Z(t-s, x-v) \vert^{q}).
\end{align*}
Using the new variables $\bar x=x-u$, $\bar y=y-u$, one can handle these terms as follows:
\begin{align*}
T^{4,1}_n(t,\cO)&\le \int_0^t ds  \int_{\cO^{s\kn}} d\bar x \int_{\cO^{s\kn}} d\bar y \int_{\IR^3} du \int_{\IR^3} dv\, G_n(s,u)G_n(s,v)\\
\\
&\qquad\qquad \times\frac{\vert D^2 f(v-u,\bar x-\bar y) \vert}{\vert \bar x-\bar y \vert^{2\bro}}
 E( \vert Z(t-s, \bar x) \vert^{q})\\
 \\
&\le C \vert \cO^{2T}\vert 
\sup_{n,x,y} \frac{\mu^4_n(x,y)}{\vert x-y \vert^{2\bro}} \int_0^t ds \int_{\cO^{s\kn}} d\bar x\, E( \vert Z(t-s, \bar x) \vert^{q})\\
\\
&\le C  \int_0^t ds \,E\left(\Vert Z(s)\Vert^q_{L^q(\cO^{(t-s)\kn})}\right).
\end{align*}
In the same way, with the change of variables $\bar x=x-v$, $\bar y=y-v$ we obtain a similar upper bound for the term $T^{4,2}_n(t,\cO)$. Consequently,
\begin{equation}\label{225}
   T^4_n(t,\cO) \leq C\int_0^t ds E\left(\Vert Z(s)\Vert^q_{L^q(\cO^{(t-s)(1+\frac{1}{n})})}\right).
\end{equation}
With (\ref{223})--(\ref{225}), we obtain (\ref{222}).

   The inequality (\ref{222'}) is a consequence of (\ref{222}). Indeed, take $\alpha=\rho+\frac{3}{q}$, which clearly satisfies the requirements of the statement.
Hence, the proof of the lemma is complete.
\hfill$\Box$
\vskip 16pt

   With these ingredients, we can now proceed to the proof of Theorem \ref{thm2.1}.
\medskip

\noindent{\it Proof of Theorem \ref{thm2.1}}:  Set $p = \qnbar/2$ and fix $t\in\,]0,T]$. Fix $m \in \IN$ and let $\cO_m = \{x \in \cO: \overline{B_{1/m}(x)} \subset \cO\}$, where $\overline{B_{1/m}(x)}$ is the closed Euclidean ball centered at $x$ with radius $1/m$. Then $\{\cO_m,\ m\in\IN\}$ is an increasing sequence of open sets and $\cup_{m \in \IN}\, \cO_m = \cO$. For $n \in \IN$, let $k_n = 1 + \frac{1}{n}$. Observe that for $n \in \IN$ such that $T/n < 1/m$, $\cO_m^{(t-s) \kn} \subset \cO^{t-s}$ if $0 \leq s \leq t \leq T$.

   We start by showing that for such $m,n \in \IN$,
\begin{align}
E\left(\Vert v^t_{G_n,Z} \Vert^{\qnbar}_{L^{\qnbar}(\cO_m)}\right) &\le C \int_0^t \,ds\,J(t-s)\,E\left(\Vert Z(s)\Vert^q_{L^q(\cO_m^{(t-s) \kn})}\right)\label{2.3.7}\\
&\le C \int_0^t \,ds\,E\left(\Vert Z(s)\Vert^q_{L^q(\cO_m^{(t-s) \kn})}\right), \label{2.3}
\end{align}
where, for $0\le s < t\le T$, 
\begin{equation}
\label{2.3.1}
J(t-s)= \sup_{\xi\in\IR^3} \int_{\IR^3} \mu(d\eta)\,\vert\mathcal{F} G(t-s)(\xi-\eta)\vert^2.
\end{equation}

Consider the (Walsh) stochastic integral $v^t_{G_n, Z}$ with $G_n$ given in (\ref{A.1}). Set
$$
   \mu_n(t,x) = \int_0^t ds \int_{\IR^3} dy\, f(y) \left(G_n(s, x - \cdot) \ast \tilde{G}_n(s, x - \cdot)\right) (y).
$$
By (\ref{fourierG}), (\ref{A.4}), (\ref{A.5}) and Assumption \ref{assumpC},
$$
   \sup_{n,x,\tbar\leq T} \mu_n(\tbar,x) \leq \sup_{\tbar \leq T} \int_0^{\tbar} ds \int_{\IR^3} \mu(d \xi)\, \vert {\cal{F}} G(s) (\xi)\vert^2 \le C.
$$
Burkholder's inequality yields
\begin{align*}
 E\left(\Vert v^t_{G_n, Z} \Vert^{\qnbar}_{L^{\qnbar}(\cO_m)} \right) \leq&  \int_{\cO_m} dx\,
 \Big(\int_0^t ds \int_{\IR^3} du \int_{\IR^3} dv\, f(u-v) \\
 & \qquad\qquad \times G_n(t-s, x-u)G_n(t-s, x-v)\\
 & \qquad\qquad \times 1_{\cO_m^{(t-s)k_n}}(u) Z(s,u) 1_{\cO_m^{(t-s)k_n}}(v) Z(s,v) \Big)^p 
\end{align*}
(note that the presence of the indicator functions does not change the value of the integral, since $G_n(t-s, x-u) = 0$ if $u \not\in  \cO_m^{(t-s)k_n}$).
H\"{o}lder's inequality yields
\begin{align}
   E\left(\Vert v^t_{G_n, Z} \Vert^{\qnbar}_{L^{\qnbar}(\cO_m)} \right) \leq& \int_{\cO_m} dx\,
     (\mu_n(T,x))^{p-1} E\Big( \int_0^t ds \int_{\IR^3} dy\, f(y) \int_{\IR^3} dz \nonumber\\
    &\nonumber \\ \nonumber
   & \qquad \times G_n(s, x-z)G_n(s, x-y-z) \\
   & \qquad \times \vert Z_{m,n}(t-s, z) \vert^p \vert Z_{m,n}(t-s, z+y)\vert^p \Big),\label{A}
\end{align} 
where $Z_{m,n}(t-s,z) = 1_{\cO_m^{(t-s)k_n}}(z)\, Z(t-s,z)$.

  Assume first that for any $s\in[0,T]$, $Z_{m,n}(t-s)\in\mathcal{C}_0^\infty(\cO_m^{(t-s)k_n})$. In this case, $G_n(s,x-\cdot)\vert Z_{m,n}(s,\cdot)\vert\in \mathcal{C}_0^\infty(\cO_m^{(t-s)k_n})$, so the last integral is bounded by
\begin{equation}\label{Aa}
  C E \left( \int_{\cO_m} dx\, \int_0^t ds \int_{\IR^3} \mu(d \eta)\, \vert {\cal{F}} [G_n(s, x-\cdot) \vert Z_{m,n}(t-s,\cdot) \vert^p](\eta)\vert^2 \right).
\end{equation}   
We now apply the arguments of the proof of Lemma 1 of \cite{DM} (in particular (2.5) there), as follows. Since the Fourier transform of a product is the convolution of the Fourier transform of the factors, we have
\begin{align*}
 & {\cal{F}} [G_n(s, x-\cdot) \vert Z_{m,n}(t-s,\cdot) \vert^p](\eta) \\
 & \qquad\qquad =\int_{\IR^3} d{\xi}'\, e^{ix\cdot(\eta-{\xi}')} {\cal{F}} G_n(s,\cdot)(\eta-{\xi}') {\cal{F}}[\vert Z_{m,n}(t-s,\cdot) \vert^p](\xi'),
\end{align*}
where $x\cdot{\xi}'$ denotes the Euclidean inner product in $\IR^3$ of the vectors $x$ and ${\xi}'$. By Plancherel's identity, 
 \begin{align*}
 &\int_{\cO_m} dx\, \vert{\cal{F}} [G_n(s, x-\cdot) \vert Z_{m,n}(t-s,\cdot) \vert^p](\eta)\vert^2\\
 &\quad \le \int_{\IR^3} dx\, \vert{\cal{F}} [G_n(s, x-\cdot) \vert Z_{m,n}(t-s,\cdot) \vert^p](\eta)\vert^2\\
  &\quad =  \int_{\IR^3} d\xi\, \vert {\cal{F}} G_n(s,\cdot)(\eta-\xi)\,{\cal{F}}\left(\vert Z_{m,n}(t-s,\cdot)\vert^p\right) (\xi)\vert^2.
  \end{align*}
Consequently, (\ref{Aa}) is bounded by
\begin{align}   
   & C \int_0^t ds\, E\left(\int_{\IR^3} d \xi\, \vert {\cal{F}} ( \vert Z_{m,n}(t-s, \cdot)\vert^p)(\xi) \vert^2\right) \int_{\IR^3} \mu
(d\eta)\, \vert{\cal{F}} G_n(s)(\xi-\eta)\vert^2 \label{2.4.1}\\
  &\leq C \int^t_0 ds\, E\left(\int_{\IR^3} dz\,  \vert Z_{m,n}(t-s, z) \vert^{2p}\right) \sup_{\xi \in \IR^3} \int_{\IR^3} \mu (d\eta)\,
\vert {\cal{F}} G(s) (\xi-\eta) \vert^2 \nonumber\\
&\leq C \int_0^t ds\, E\left(\Vert Z_{m,n}(s) \Vert^{\qnbar}_{L^{\qnbar}(\IR^3)}\right) \sup_{\xi \in \IR^3} \int_{\IR^3} \mu(d \eta)\,
\vert {\cal{F}} G(t-s)(\xi-\eta) \vert^2 .
\label{2.4}
\end{align}
Since supp~$Z_{m,n}(s) \subset \cO_m^{(t-s)k_n}$, we can replace $\Vert \cdot \Vert_{L^{\qnbar}(\IR^3)}$ by $\Vert \cdot \Vert_{L^{\qnbar}(\cO_m^{(t-s)k_n})}$ in (\ref{2.4}).

   Without the smoothness and compact support restrictions on $Z_{m,n}(t-s)$, we can check, by regularizing $\vert Z\vert 1_{\cO_m^{(t-s)k_n}}$ by means of convolution and using Fatou's Lemma, that the inequality
\begin{align}
&E\left(\Vert v^t_{G_n, Z} \Vert^{\qnbar}_{L^{\qnbar}(\cO_m)} \right)\nonumber\\
&\quad \le C \int_0^t ds\, E\left(\Vert Z(s) \Vert^{\qnbar}_{L^{\qnbar}(\cO_m^{(t-s)k_n})}\right) \sup_{\xi \in \IR^3} \int_{\IR^3} \mu(d \eta)\,
\vert {\cal{F}} G(t-s)(\xi-\eta) \vert^2\nonumber\\
&\quad = C \int_0^t ds\, J(t-s) E\left(\Vert Z(s) \Vert^{\qnbar}_{L^{\qnbar}(\cO_m^{(t-s)k_n})}\right)\label{2.4.2}
\end{align}
holds for any process $Z$ satisfying the assumptions of the Theorem. This proves (\ref{2.3.7}).

Lemma 8 in \cite{DM} and (\ref{A.3}) yield 
%$$ \int_0^T ds \sup_{\xi \in \IR^3} \int_{\IR^3} \mu(d \eta)\, \vert {\cal{F}}G(s)(\xi-\eta)\vert^2 < \infty.$$ 
\begin{equation}
\label{J}
   \sup_{0<s<t\le T}J(t-s)\le \int_{\IR^3} \mu(d\eta) \vert{\mathcal{F}} G(t-s)(\eta)\vert^2 \le C.
\end{equation}
This establishes  (\ref{2.3}).

As in \cite[Lemma 5]{DM}, one checks that, for any fixed $m\in\mathbb{N}$,
\begin{align}\label{eq2.4'}
  & \lim_{n \to \infty} \sup_{t \in [0,T]} E\left(\Vert v^t_{G_n,Z}- v^t_{G,Z} \Vert^2_{L^2(\cO_m)}\right)\\ \nonumber
  &\qquad\qquad \leq \lim_{n \to \infty} \sup_{t \in [0,T]} E\left(\Vert v^t_{G_n,Z}- v^t_{G,Z} \Vert^2_{L^2(\cO_m)}\right)\\ \nonumber
  &\qquad\qquad =  0. 
\end{align}
By Fatou's Lemma, this yields
$$
E\left(\Vert v^t_{G,Z}\Vert^{\qnbar}_{L^{\qnbar}(\cO_m)}\right) \leq \liminf_{k \to \infty} E\left(\Vert v^t_{G_{n_k},Z}
\Vert^{\qnbar}_{L^{\qnbar}(\cO_m)}\right)
%\leq \sup_{n \geq 1} E\left(\Vert v^t_{G_n, Z}\Vert^{\qnbar}_{L^{\qnbar}(\cO_m)}\right),
$$
for some subsequence $(n_k)_{k \geq 1}$. %With (\ref{2.4.2}) this proves the estimate  (\ref{2.3.7}).

% and even, using (\ref{2.4}),
%\begin{equation}\label{2.5a}
%   E\left(\Vert v^t_{G,Z} \Vert^{\qnbar}_{L^{\qnbar}(\IR^3)}\right) \leq C\int_0^t ds\, E\left(\Vert Z(s)
%\Vert^{\qnbar}_{L^{\qnbar}(\IR^3)}\right) \sup_{\xi \in\IR^3} \int_{\IR^3} \mu(d\eta)\, \vert {\cal{F}} G(t-s) (\xi-\eta)\vert^2 .
%\end{equation}

Since $\cO_m^{(t-s)k_n} \subset \cO^{t-s}$, for any $n, m \in \mathbb{N}$ satisfying $\frac{T}{n} < \frac{1}{m}$, we deduce from (\ref{2.4.2}) 
and (\ref{eq2.4'}) that
$$
  E\left(\Vert v^t_{G,Z} \Vert^{q}_{L^{q}(\cO_m)}\right) \leq C\int_0^t ds\, E\left(\Vert Z(s)\Vert^{q}_{L^{q}(\cO^{t-s})}\right).
$$
Let $m \to \infty$. Using the monotone (increasing) convergence theorem, we obtain
\begin{equation}\label{p16star}
E\left(\Vert v^t_{G,Z} \Vert^{q}_{L^{q}(\cO)}\right) \leq C\int_0^t ds\, E\left(\Vert Z(s)\Vert^{q}_{L^{q}(\cO^{t-s})}\right).
\end{equation}

   Fix $\ro\in\,]0,\tau_q(\beta,\delta)\wedge\gamma[$.  By Fatou's Lemma and (\ref{222'}),
\begin{align}
 E\left(\Vert v^t_{G,Z} \Vert^{q}_{\ro,q,\cO_m}\right)&\leq 
 \liminf_{\ell \to \infty} E\left(\Vert v^t_{G_{n_\ell},Z} \Vert^{q}_{\ro,q,\cO_m}\right)\nonumber\\
 &\le C \liminf_{\ell \to \infty} \int_0^t ds\, E\left(\Vert Z(s)\Vert^q_{W^{\rho,q}(\cO_m^{(t-s)(1+\frac{1}{n_\ell})})}\right)\label{228},
 \end{align}
for some subsequence $(n_\ell)_{\ell \geq 1}$. %By virtue of (E.31) in \cite{shima} and the assumption (\ref{2.2}), we obtain
Since $ \cO_m^{(t-s)(1+\frac{1}{n_\ell})}\subset\cO^{t-s}$ for $\frac{T}{n_l}<\frac{1}{m}$, we obtain
\begin{equation*}
   E\left(\Vert v^t_{G,Z} \Vert^{q}_{\ro,q,\cO_m}\right)\le C \int_0^t ds\, E\left(\Vert Z(s)\Vert^q_{W^{\rho,q}(\cO^{t-s})}\right).
\end{equation*}
Let $m \to \infty$. Using monotone convergence, we obtain
\begin{equation}\label{p16star2}
   E\left(\Vert v^t_{G,Z} \Vert^{q}_{\ro,q,\cO}\right)\le C \int_0^t ds\, E\left(\Vert Z(s)\Vert^q_{W^{\rho,q}(\cO^{t-s})}\right).
\end{equation} 
%by monotone convergence.
%Then $\alpha:=\ro+\frac{3}{q}$ satisfies the restrictions of
%Lemma \ref{l2.2}. Thus, for this value of $\alpha$, for any compact set $\cO$, the inequality (\ref{222}) becomes
%Owing to (\ref{222'}),
%\begin{equation}
%\label{228}
%E\left(\Vert v^t_{G_n,Z}\Vert^q_{\ro,q,\cO}\right)\le C \int_0^t ds E\left(\Vert Z(s)\Vert^q_{W^{\ro,q}(\cO^{(t-s)(1+\frac{1}{n})})}\right)
%\end{equation}

%By Fatou's Lemma, (\ref{228}) and (\ref{2.2}), %for some subsequence $(n_k)_{k \geq 1},$
%\begin{align*}
%& E\left(\Vert v^t_{G,Z} \Vert^{q}_{\ro,q,\cO}\right)\\
%&\quad \leq \liminf_{k \to \infty} E\left(\Vert v^t_{G_{n_k},Z} \Vert^{q}_{\ro,q,\cO}\right)\\
% &\quad \leq  C \int_0^t ds E\left(\Vert Z(s)\Vert^q_{W^{\rho,q}(\cO^{t-s})}\right).
%\end{align*}

%   We have shown that
%$$
%   E\left(\Vert v^t_{G,Z} \Vert^{q}_{W^{\rho,q}(\cO_m)}\right) \leq C \int_0^t ds\, E\left(\Vert Z(s)\Vert^q_{W^{\rho,q}(\cO_m^{(t-s)k_n})}\right).
%$$
The inequalities (\ref{p16star}) and (\ref{p16star2}) together establish (\ref{2.2.1}).
This completes the proof of the theorem. 
\hfill $\Box$
\vskip 16pt
\smallskip

\begin{remark}
\label{rem2.2}
{\rm 
Fix a bounded domain $\cO\subset \IR^3$ and assume that the stochastic process $Z$ in Theorem \ref{thm2.1} is such that the right-hand side of 
(\ref{2.2.1}) is finite for any $t\in\, ]0,T]$. 
%This holds for instance if $Z$ satisfies 
%$$
%\int_0^t ds\,\left[E\left(\Vert Z(s)\Vert^q_{W^{\gamma,q}(\bar{\cO}^{(t-s)(1+\varepsilon)})}\right) + E\left(\Vert Z(s)\Vert^q_{L^q(\IR^3)}\right)\right] < \infty,
%$$
%for any  $t\in]0,T]$.
By the Sobolev embedding theorem (see for instance \cite[Theorem E.12 p.257]{shima}), Theorem \ref{thm2.1} yields that, for each $t\in[0,T]$, a.s., $x
\mapsto v^t_{G,Z}(x)$ is $a$-H\"{o}lder continuous, 
%uniformly in $t \in [0, T]$,
with $a \in\,]0, (\gamma \wedge \tau_q(\beta,\delta)) - \frac{3}{\qnbar}[$.
 Indeed, for any bounded or unbounded domain $\cO \subset
\IR^d$, $W^{\ro, q}(\cO) \subset \cC^\beta(\cO),$ for each $\beta < \ro - \frac{d}{q}.$
}
\end{remark}
\smallskip

%\begin{remark}
%\label{rem2.3} \ref{?} REMOVE ?
%{\rm Consider the particular case of integrands  $Z=\{Z(t,x), (t,x)\in[0,T]\times \IR^3\}$ 
%which are real-valued stochastic processes with
%spatial stationary finite dimensional distributions. That is, for any $m\ge 1$, $x_1,\cdots,x_m \in \IR^3$, $t\in[0,T]$,
%$x\in\IR^3$, the joint distribution of $(Z(t,x_1+x), \cdots,Z(t,x_m+x))$ coincides with that of $(Z(t,x_1), \cdots,Z(t,x_m))$.
%In \cite{dalang},  examples of such processes are given; they are obtained as solutions of a class of spde's.

%The next result is a version of Theorem \ref{thm2.1} for this class of integrands. 

%\begin{thm}
%\label{thm2.2} \ref{?} REMOVE ?
%Fix $q\in]3,\infty[$ and a bounded domain $\mathcal{O}\subset\IR^3$. Suppose that $\tau_q(\beta,\delta):= (\frac{2-\beta}{2} \wedge \frac{1+\delta}{2}) - \frac{3}{q} > 0$. 

%Then, for any $t\in[0,T]$ and every $\rho \in \, ]0,\tau_q(\beta,\delta) \wedge\gamma[$,
%\begin{align*}
%E\left(\Vert v^t_{G,Z} \Vert^{\qnbar}_{W^{\rho,q}(\mathcal{O})}\right)&\le C \int_0^t ds \Big[E\left(\Vert Z(s)\Vert^q_{L^q(\IR^3)}\right)\\
%&+ E\left(\Vert Z(s)\Vert^q_{\rho,q,\cO}\right)\Big].
%\end{align*}
%\end{thm}
%}
%\end{remark}

%%%%%%%%%%%%%%%%%%%%%%%%%%%%%%%%%%%%%%%%
\subsection{Regularity in time of the stochastic integral}
\label{s1.3}
 This section is devoted to the analysis of the H\"{o}lder continuity in time 
of the stochastic integral process $\{v^t_{G,Z}(x),\ t \in [0,T]\}$, when $x$ is fixed. Throughout this section, $\cO$ denotes a bounded domain in $\IR^3$ and we shall make the following assumption on the integrand process $Z$:
\medskip

\begin{assump} For some fixed $q\in\, ]3,\infty[$ and $\gamma\in\, ]0,1[$, 
\begin{equation*}
   \sup_{t \in [0,T]} E\left(\Vert Z(t)\Vert^q_{W^{\gamma,q}(\cO^{T-t})} \right)  
< \infty,
\end{equation*}
\label{assumpI}
\end{assump}

   Notice that Assumption \ref{assumpI} implies (\ref{2.2}). Therefore, by Remark \ref{rem2.2}, it makes sense to fix the argument $x\in\cO$ in the stochastic integral process. 
In addition, by the above mentioned Sobolev embedding, there is $C < \infty$ such that the H\"older norm $\Vert \cdot \Vert_{\cC^\rho(\cO)}$ is bounded by a constant times the Sobolev norm $\Vert \cdot \Vert_{W^{\gamma,q}(\cO)}$, provided $\rho \in \, ]0,\gamma - \frac{3}{q}[$. In particular, for any stochastic process $Z$ 
%satisfying the assumption (b) of Theorem \ref{thm2.1} for some %$\ro\in]\frac{3}{q},1[$
satisfying Assumption \ref{assumpI}, one has 
$$
   \sup_{t \in [0, T]} E(\vert Z(t,x) - Z(t, y) \vert^{\qnbar}) \leq C \vert x-y \vert^{\rho \qnbar}.
$$
By H\"older's inequality, we deduce that for $\qbar \in \, ]0,q]$,
\begin{equation}\label{2.11}
   \sup_{t \in [0, T]} E(\vert Z(t,x) - Z(t, y) \vert^{\qbar}) \leq C \vert x-y \vert^{\rho \qbar},
\end{equation}
for any $x,y \in \cO^{T-t}$, $\rho \in\, ]0, \gamma - \frac{3}{q}[$, $\qbar \in\, ]0,q]$, with a positive constant $C$ not depending on $x$, $y$ or $\qbar$. Moreover,
\begin{equation}\label{2.12}
\sup_{t \in [0, T]}\sup_{x \in \cO^{T-t}} E(\vert Z(t,x)\vert^q) < \infty.
\end{equation}

\begin{thm}
\label{thm2.3}
Let $\tau(\beta,\delta) = \frac{2-\beta}{2}\wedge\frac{1+\delta}{2}$ and let $Z$ be a stochastic process satisfying Assumption \ref{assumpI} for a fixed $q \in\,]3, \infty[$ and $\gamma \in\,]\frac{3}{q},1[$. Then the stochastic process $\{v^t_{G,Z}(x),\ t \in [0,T]\}$ ($x \in \IR^3$) satisfies
\begin{equation}\label{2.13}
\sup_{x \in \cO} E\left(\vert v^t_{G,Z} (x) - v^{\tbar}_{G,Z}(x) \vert^{\qbar}\right) \leq C \vert t- \tbar \vert^{\rho \qbar},
\end{equation}
for each $t, \tbar \in [0,T]$, any $\qbar \in\, ]2, \qnbar[$ and $\rho \in\,]0, (\gamma - \frac{3}{q}) \wedge
\tau(\beta,\delta)[$. Consequently, the process $\left(v^t_{G,Z}(x), t \in [0, T]\right)$ is
a.s.~$a$-H\"{o}lder continuous in $t$, for any $a \in\, ]0, ((\gamma - \frac{3}{q}) \wedge
\tau(\beta,\delta)) - \frac{3}{\qnbar}[$.
\end{thm}

%Notice that since {\em{supp}}~$G(t-s) \subset B_{t-s}(0)$ and $Z(t)$ has compact support a.s., the supremum on the left
%hand-side of (\ref{2.13}) is limited to a compact set. Actually, to the set $K^T$, where $K$ is described in assumption (a) of %Theorem \ref{thm2.1}.

\noindent{\it Proof}: Let $\cO_m$ be the open sets defined at the beginning of the proof of Theorem \ref{thm2.1} and set $k_n = 1 + \frac{1}{n}$. Recall that $\cO_m^{(t-s)k_n} \subset \cO^{t-s}$ if $0 \leq s \leq t \leq T$ and $\frac{T}{n} < \frac{1}{m}$.

   Set $p=\qbar/2$, $p \in\,]1, \infty[$ and fix $0 \leq t \leq \tbar \leq T$. The first part of the proof is devoted to
showing that there is $C < \infty$ such that for all $m \in \IN$,
\begin{equation}\label{eq2.14}
\sup_{n \geq mT} \sup_{x \in \cO_m} E \left( \vert v^{\tbar}_{G_n,Z}(x) - v^t_{G_n,Z}(x) \vert^{2p}\right) \leq C \vert t- \tbar
\vert^{\tilde{\rho} p},
\end{equation}
with $\tilde{\rho} \in\,]0, 2 ((\gamma - \frac{3}{2p}) \wedge \tau(\beta,\delta))[$, where the $G_n$, $n \geq 1$, are defined in (\ref{A.1}). Indeed, for $x \in \cO_m$, consider the decomposition
$$
   E\left(\vert v^{\tbar}_{G_n, Z}(x) - v^t_{G_n,Z}(x) \vert^{2p}\right) \leq C \left(T^n_1(t, \tbar, x) + T^n_2(t, \tbar, x)\right),
$$
where
\begin{align*}
T_1^n(t, \tbar, x) &= E \Big(\big( \int^{\tbar}_t \int_{\IR^3} G_n(\tbar-s, x-y) Z(s,y) M(ds, dy)
\big)^{2p}\Big),\\
   T^n_2(t, \tbar, x) &= E\Big(\big(\int_0^t \int_{\IR^3}(G_n(\tbar-s, x-y) - G_n(t-s, x-y))\\ 
   &\qquad\qquad \times Z(s,y) M(ds, dy)
\big)^{2p}\Big).
\end{align*}
Burkholder's inequality yields
\begin{align*}
   T_1^n(t,\tbar,x) &\leq C \ E \Big( \int_t^{\tbar} ds \int_{\IR^3} dy \int_{\IR^3} dz\, G_n (\tbar-s, x-y) G_n(\tbar-s,x-z)\\ 
&\quad \times f(y-z) Z(s,y) Z(s,z) \Big)^p\\
&= C \ E \Big( \int_0^{\tbar-t} ds \int_{\IR^3} dy \int_{\IR^3} dz\, G_n(s,x-y) G_n(s,x-z)\\
&\quad \times  f(y-z) Z(\tbar-s,y)Z(\tbar-s,z)\Big)^p.
\end{align*}
Set
$$
\mu^n(x,t,\tbar) = \int_0^{\tbar-t} ds \int_{\IR^3} dy \int_{\IR^3} dz\, G_n(s,x-y)G_n(s, x-z) f(y-z).
$$
By Assumption \ref{assumpC} and Lemma \ref{lemA.2},
$$
   \sup_{n \geq 1} \sup_{x \in \IR^3} \mu^n(x,t,\tbar) \leq \int_0^{\tbar-t} ds \int_{\IR^3} \frac{\vert {\cal{F}} G(s)(\xi)
\vert^2}{\vert \xi \vert^{3-\beta}} d\xi \leq C \vert \tbar-t \vert^{3-\beta}.
$$
Hence, applying H\"{o}lder's inequality and using (\ref{2.12}), we obtain
\begin{equation}\label{eq2.15}
\sup_{n \geq mT} \sup_{x \in \cO_m} T^n_1 (t, \tbar, x) \leq C \vert \tbar-t \vert^{p(3-\beta)},
\end{equation}
where $C$ does not depend on $m$.

   We now study the contribution of $T^n_2(t, \tbar, x)$. We proceed first in a manner analogous to the proof of Theorem
\ref{thm2.1}. After having applied Burkholder's inequality, changed variables (using (\ref{A.2})) and rearranged terms, we obtain 
$$
   T^n_2(t, \tbar, x) \leq C \sum_{i=1}^4 E\left(\vert R^{i,n}(t, \tbar, x)\vert^p\right),
$$
where
$$
   R^{i,n} (t,\tbar,x) = \int^t_0 ds \int_{\IR^3} du \int_{\IR^3} dv\, G_n(t-s, u) G_n(t-s,v) r_i(t, \tbar,s,x,u,v),
$$
and
\begin{align*}
r_1(t, \tbar,s,x,u,v)=&f\left(\frac{\tbar-s}{t-s} v-u\right) \frac{\tbar-s}{t-s} \left(Z(s, x - \frac{\tbar-s}{t-s} u) - Z(s,x-u)\right)\\
& \qquad \times \left(Z(s,x-\frac{\tbar-s}{t-s} v) - Z(s, x-v)\right),
\end{align*}
\begin{align*}
r_2(t, \tbar,s,x,u,v)&=\left(\left(\frac{\tbar-s}{t-s}\right)^2 f\left(\frac{\tbar-s}{t-s}(v-u)\right)-\frac{\tbar-s}{t-s} f\left(\frac{\tbar-s}{t-s}
v-u\right)\right)\\
&\qquad \times Z\left(s,x - \frac{\tbar-s}{t-s} u\right)\left(Z(s,x - \frac{\tbar-s}{t-s} v) - Z(s,x-v)\right),
\end{align*}
\begin{align*}
r_3(t, \tbar,s,x,u,v)&= \left(\left(\frac{\tbar-s}{t-s}\right)^2 f\left(\frac{\tbar-s}{t-s}(v-u)\right)- \frac{\tbar-s}{t-s} f\left(v -
\frac{\tbar-s}{t-s} u\right)\right) \\
&\qquad \times Z(s,x-v) \left(Z(s,x - \frac{\tbar-s}{t-s}u) - Z(s,x-u)\right),
\end{align*}
\begin{align*}
r_4(t, \tbar,s,x,u,v) &= \Bigg(\left(\frac{\tbar-s}{t-s}\right)^2f\left(\frac{\tbar-s}{t-s} (v-u)\right)- \frac{\tbar-s}{t-s} f\left(\frac{\tbar-s}{t-s}
v-u\right)\\
&\qquad  - \frac{\tbar-s}{t-s} f\left(v- \frac{\tbar-s}{t-s} u\right) + f(v-u)\Bigg)\\
\\
&\qquad\times Z(s,x-u)Z(s,x-v). 
\end{align*}
Set
$$
   \nu_1^n(s,t, \tbar) = \int_{\IR^3} du \int_{\IR^3} dv\, G_n(t-s,u) G_n(t-s,v) f\left(\frac{\tbar-s}{t-s} v-u\right)
\frac{\tbar-s}{t-s}.
$$
By Lemma \ref{lemB.4}, we can apply H\"{o}lder's inequality and then Schwarz's inequality to obtain
\begin{align*}
& E\left(\vert R^{1,n}(t,\tbar,x)\vert^p\right)
 \leq \left(\sup_{n \geq 1} \sup_{0\leq s\leq t\leq \tbar\leq T}\nu_1^n (s,t, \tbar)\right)^{p-1}\\
 \\
  &\qquad \times \int_0^t ds \int_{\IR^3} du \int_{\IR^3} dv\, G_n(t-s,u) G_n(t-s,v)
f\left(\frac{\tbar-s}{t-s} v-u\right) \\
\\
 &\qquad\times  \frac{\tbar-s}{t-s}\left( E\left( \left\vert Z\left(s,x- \frac{\tbar-s}{t-s} u\right) - Z(s,x-u)
\right\vert^{2p}\right) \right)^{1/2}\\
\\
 &\qquad\qquad\times  \left( E\left(\left\vert Z\left(s,x- \frac{\tbar-s}{t-s}v\right) - Z(s,x-v)\right\vert^{2p}\right)\right)^{1/2}.
\end{align*}
By (\ref{2.11}), for $\rho \in\, ]0, \gamma - \frac{3}{q}[$, the product of the last two factors is bounded by
$$
%   (\nu^n_1(t,\tbar))^{p-1} \int^t_0 ds \int_{\IR^3} du \int_{\IR^3}dv\, G_n(t-s,u) G_n(t-s,v) f(\frac{\tbar-s}{t-s} v-u)
\frac{\tbar-s}{t-s}\\
%\\  \times 
   \left(\frac{\tbar-t}{t-s}\right)^{2 p \rho} \vert u \vert^{\rho p} \vert v \vert^{\rho p}.
$$
Since for any $t \in (0,T],$ supp $G_n(t, \cdot) \subset B_{t(1+ \frac{1}{n})}(0),$  $(\vert u\vert \vert v \vert /\vert
t-s \vert^2)^{\rho p}$ in the last integral is bounded by a finite constant. Therefore, by Lemma \ref{lemB.4},
\begin{align}
\sup_{n \geq mT} \sup_{x \in \cO_m} E\left( \vert R^{1,n} (t, \tbar, x) \vert^p\right)
  &\leq C \vert t- \tbar \vert^{2p\rho} \sup_{n \geq 1} \sup_{0\leq s\leq t\leq \tbar\leq T} 
     \left(\nu^n_1 (t, \tbar)\right)^p\nonumber\\
   &\leq C \vert t- \tbar \vert^{2p\rho} \label{marta1},
\end{align}
with $\rho \in\, ]0, \gamma - \frac{3}{q}[$ and $C$ does not depend on $m$.

   Taking into account the result proved in Lemma \ref{lemB.5} and using the quantity $\nu_2^n(t,\tbar)$ defined in that lemma, we apply first H\"{o}lder's inequality, then Schwarz's
inequality, to obtain
\begin{align*}
E\left(\vert R^{2,n} (t,\tbar, x) \vert^p\right) &\leq \left(\nu^n_2(t,\tbar)\right)^{p-1} \int_0^t ds \int_{\IR^3} du \int_{\IR^3} dv\, G_n(t-s,u) G_n(t-s,v)\\
\\
&\quad \times\, \left\vert \left(\frac{\tbar-s}{t-s}\right)^2 f\left(\frac{\tbar-s}{t-s}(v-u)\right) - \frac{\tbar-s}{t-s}
f\left(\frac{\tbar-s}{t-s} v-u\right)\right\vert\\
\\
&\quad\times \left(E\left( \left\vert Z\left(s,x- \frac{\tbar-s}{t-s} u\right)\right \vert^{2p}\right)\right)^{1/2}\\
\\
&\quad\times \left( E\left(\left\vert Z\left(s,x- \frac{\tbar-s}{t-s} v\right)-Z(s,x-v)\right \vert^{2p}\right)\right)^{1/2}.
\end{align*}
Then, by Lemma \ref{lemB.5}, (\ref{2.11}) and (\ref{2.12}),
\begin{equation}\label{eq2.16}
   \sup_{n \geq mT} \sup_{x \in \cO_m} E\left( \vert R^{2,n} (t, \tbar, x) \vert^p\right) \leq C \vert t-\tbar \vert^{p(\rho+\alpha)},
\end{equation}
for any $\alpha \in\, ]0,1[$ with $\alpha + \beta \in\, ]0,2[$, $\rho \in\, ]0, \gamma - \frac{3}{q}[$, and $C$ does not depend on $m$. The same result holds for the term $R^{3,n}(t, \tbar, x).$

 Using Lemma \ref{lemB.6} and the quantity $\nu^n_3(t,\tbar)$ defined in that lemma, we can apply H\"{o}lder's inequality to obtain
\begin{align*}
E\left(\vert R^{4,n}(t,\tbar, x) \vert^p\right) &\leq \left(\nu^n_3(t, \tbar)\right)^{p-1} \int_0^t ds \int_{\IR^3} du \int_{\IR^3} dv\,
G_n(t-s,u) G_n(t-s,v)\\
\\
& \quad \times \left \vert \Delta^2 f(s,t,\tbar,u,v)\right\vert\, E\left(\vert Z(s,x-u) \vert^p \vert Z(s,x-v)\vert^p\right),
\end{align*}
where $\Delta^2 f(s,t,\tbar,u,v)$ is defined in (\ref{B.13.1}).

   Finally, by Schwarz's inequality, property (\ref{2.12}) and (\ref{B.14}), we reach
\begin{eqnarray}
   \sup_{n \geq mT} \sup_{x \in \cO_m} E\left( \vert R^{4,n}(t,\tbar,x) \vert^p\right) &\leq& C \sup_{n \geq 1}
   \left(\nu^n_3 (t, \tbar)\right)^p
   \sup_{t \in [0,T]} \sup_{x \in \cO^{T-t}} E\left( \vert Z(t,x)\vert^{2p}\right) \nonumber
   \\
   &\leq& C \vert t - \tbar \vert^{\alpha p},
\label{eq2.18}
\end{eqnarray}
with $\alpha \in\, ]0,(2-\beta) \wedge (1+\delta)[$ and $C$ does not depend on $m$. Hence, (\ref{eq2.15})-(\ref{eq2.18}) establish (\ref{eq2.14}).

   The second part of the proof consists in deducing (\ref{2.13}) from (\ref{eq2.14}). To this end, we first prove that for any fixed $m \in \IN$ and $t \in [0,T]$, $(v^t_{G_n,Z},\ n \geq mT)$ is a sequence of bounded and equicontinuous functions defined on $\cO_m$ with values in $L^{\qbar}(\Omega)$, for any  $\qbar \in [1, q]$. Indeed, 
%we have proved in \ref{thm2.1} (see (\ref{z2}) that
from (\ref{2.4.2}) and (\ref{222'}) together with Assumption \ref{assumpI} and the inclusion $\cO_m^{(t-s)k_n} \subset \cO^{t-s}$, we see that
\begin{equation}\label{eq2.19}
  \sup_{m \in \IN} \sup_{n \geq mT} \sup_{t\in[0,T]} E\left(\Vert v^t_{G_n,Z} \Vert^{q}_{W^{\rho,q}(\cO_m)}\right) < \infty
\end{equation}
for any $\rho \in\, ]0, \rho_0[$, with $\rho_0=\gamma \wedge\left(\tau(\beta,\delta) - \frac{3}{q}\right)$.
Therefore, the Sobolev embedding yields
$$
   \sup_{m \in \IN} \sup_{n \geq mT}\sup_{t\in[0,T]} E\left(\Vert v^t_{G_n,Z} \Vert^q_{\cC^{\tilde{\rho}}(\cO_m)}\right) < \infty,
$$
for any $\tilde{\rho} \in\, ]0, \rho_0 - \frac{3}{q}[$. Consequently, for
any $x,y \in \cO_m$,
\begin{equation}\label{eq2.20}
\sup_{n \geq mT} E\left( \vert v^t_{G_n,Z}(x) -v^t_{G_n,Z}(y) \vert^{\qbar}\right) \leq C \vert x-y \vert^{\tilde{\rho} \qbar},
\end{equation}
for every $\qbar \in[1, q]$, where $C$ does not depend on $m$. This establishes the desired boundedness and equicontinuity properties.

%   Properties (\ref{eq2.19}) and (\ref{eq2.20}) yield the announced result. 
   We now establish a uniform convergence result. The convergence (\ref{eq2.4'}) implies the
existence of a subsequence $(n_k,\ k \geq 1)$ and a Lebesgue-measurable subset $N \subset \cO_m$ with null Lebesgue measure such that for any $x \in \cO_m \setminus N$,
$$
   \lim_{k \to \infty} E\left( \vert v^t_{G_{n_k},Z}(x) - v^t_{G,Z}(x) \vert^2\right) = 0.
$$
The property  (\ref{eq2.20}) yields $\sup_{n \geq mT} E(\vert v^t_{G_n,Z}(x) \vert^q) < \infty.$ Hence, by uniform integrability, for any $x \in \cO_m \setminus N$ and for any $\qbar \in [1,q[$,
\begin{equation}\label{eq2.21}
   \lim_{k \to \infty} E\left(\vert v^t_{G_{n_k, Z}}(x) - v^t_{G,Z} (x) \vert^{\qbar}\right) = 0.
\end{equation}

   Thus, we have proved that the bounded and equicontinuous sequence $(v^t_{G_{n_k},Z},\ k \geq 1)$ of functions converges
pointwise in $L^{\qbar}(\Omega)$, on a dense subset $\cO_m \setminus N$ of $\cO_m$. Therefore, since $v^t_{G,Z}$ is $L^{\qbar}(\Omega)$--continuous (see Remark \ref{rem2.2}), we have
\begin{equation}\label{eq2.22}
   \lim_{k \to \infty} \sup_{x \in \cO_m} E\left(\vert v^t_{G_{n_k}, Z} (x) - v^t_{G,Z}(x) \vert^{\qbar}\right) = 0.
\end{equation} 
(see for instance \cite[Ch. VII]{dieudonne}), that is, $(v^t_{G_{n_k}, Z},\ k \in \IN) \subset C(\cO_m,L^{\qbar}(\Omega))$ converges uniformly to $v^t_{G,Z}$.

   Clearly, for any $m \in \IN$ and $t, \tbar \in [0, T]$,
\begin{equation}\label{eq2.23}
    \sup_{x \in \cO_m} E\left(\vert v^t_{G,Z}(x) - v^{\tbar}_{G,Z}(x) \vert^{\qbar}\right) \leq C \sum^3_{i=1} A_i(t, \tbar),
\end{equation}
 where $C$ does not depend on $m$ and
 \begin{eqnarray*}
    A_1(t, \tbar) &=& \sup_{x \in \cO_m} E\left(\vert v^t_{G,Z} (x) - v^t_{G_{n_k},Z}(x) \vert^{\qbar}\right),\\
 \\
    A_2(t, \tbar)&=& \sup_{k \geq 1} \sup_{x \in \cO_m} E\left(\vert v^t_{G_{n_k}, Z}(x) - v^{\tbar}_{G_{n_k},Z}(x)
\vert^{\qbar}\right),\\
 \\
   A_3(t, \tbar) &=& \sup_{x \in \cO_m} E\left(\vert v^{\tbar}_{G_{n_k},Z}(x) - v^{\tbar}_{G,Z}(x)\vert^{\qbar}\right).
 \end{eqnarray*}
 By (\ref{eq2.23}), (\ref{eq2.22}) and (\ref{eq2.14}), we obtain
$$
    \sup_{x \in \cO_m} E\left(\vert v^t_{G,Z}(x) - v^{\tbar}_{G,Z}(x) \vert^{\qbar}\right) \leq C \vert t - \tbar \vert^{\rho
\qbar} ,
$$
for any $\rho \in\, ]0, (\gamma - \frac{3}{q}) \wedge \tau(\beta,\delta)[$. Since $C$ does not depend on $m$, $\cO_m$ can de replaced by $\cO$ in the inequality above and this establishes (\ref{2.13}).
\hfill \qed

%%%%%%%%END OF SECTION 3

%%%%%%%% BEGINNING OF SECTION 4

\section{Path properties of the solution of the sto-\break chastic wave equation}
\label{s2}

   This section is devoted to studying the properties of the sample paths of the stochastic wave equation in spatial dimension three. More precisely, consider the s.p.d.e.
\begin{equation}\label{e3.1}
\left(  \frac{\partial^2}{\partial t^2}-\Delta\right) u(t, x) = \sigma\left(u(t,x)\right) \dot{F}(t,x) + b\left(u(t,x)\right),
    \,\, t \in\, ]0, T],\ x \in \IR^3,
\end{equation}
with initial conditions $u(0,x) = v_0(x)$, $\frac{\partial}{\partial t} u(0,x) = \tilde v_0(x)$. We are interested in the solution
$u(t,x)$ of this equation for $(t,x) \in\, ]0,T] \times D,$ where $D$ denotes a bounded domain of $\IR^3$ included in the ball
$B_{m_0}(0)$, for some $m_0 > 0$.
%Throughout this section, we assume the following 
%\vskip 12pt

%\noindent{\bf Condition (S).} NOT NEEDED \ref{?} The supports of the functions $v_0$ and $\tilde v_0$ are contained in the ball $B_{r_0}(0)$, for some $r_0>0$.
%\vskip 12pt

%   For simplicity, we set $D = \overline{B_{m_0}(0)},$ where $m_0 > 0$. 

For any $a\ge1$, set
$$
   K^D_a(t) = \{y \in \IR^3: d(y,D)\leq a(T-t)\}, \,\, t \in [0,T],
$$
where $d$ denotes the Euclidean distance. Notice that $t \mapsto K^D_a(t)$ is decreasing.

\subsection{Existence, uniqueness and moments}
\label{ss2.1}

   Recall that the solution $w$ of the homogeneous wave equation with the same initial conditions as $u$ is $w(t) = \frac{d}{dt} G(t) \ast v_0 + G (t) \ast \tilde v_0$. We are mainly interested in the solution of (\ref{e3.1}) for $(t,x) \in [0,T]\times D$, though we will need to construct the solution in a slightly larger set that contains the ``past light cone" of $\{T\}\times D$.
%The relation between the size of the set $D$ and the support of the initial conditions ensures that outside of $K^D(t),$ the functions $w(t,\cdot)$ %and $u(t,\cdot)$ vanishes. 
Therefore, we term a {\em solution of the s.p.d.e.}~(\ref{e3.1})~{\em ``in $D$"} a stochastic process $\left(u(t) 1_{K^D_a(t)},\ t \in [0,T] \right)$ with values in
$L^2(\IR^3)$, satisfying
\begin{align}
&u(t, \cdot) 1_{K^D_a(t)}(\cdot) = 1_{K^D_a(t)} (\cdot) \left(\frac{d}{dt} G(t) \ast v_0 + G(t) \ast \tilde v_0\right) (\cdot) \nonumber\\
&\qquad + 1_{K^D_a(t)} (\cdot) \int_0^t \int_{\IR^3} G(t-s, \cdot -y) \sigma \left(u(s,y)\right) 1_{K^D_a(s)}(y) M(ds, dy)
\nonumber\\
& \qquad + 1_{K^D_a(t)} (\cdot) \int^t_0 ds\, G(t-s) \ast \left(b(u(s, \cdot))1_{K^D_a(s)}(\cdot)\right) ,
\label{e3.2'}
\end{align}
a.s., for any $t \in [0, T].$ The integrands in (\ref{e3.2'}) have compact support. Further, the support of the measure $G(t-s)$ is the boundary of the ball $B_{t-s}(0)$. Therefore, if $x \in K^D_a(t)$ and $x-y \in \partial B_{t-s}(0)$, then $d(y, D) \leq
d(y,x)+d(x,D) \leq t-s + a(T-t) \leq a(T-s).$ Consequently, $y \in K^D_a(s)$. Therefore, $(K^D_a(t))^{t-s} \subset K^D_a(s)$ and so (\ref{e3.2'}) is coherent.

   The stochastic integral in (\ref{e3.2'}) is to be considered in the sense of \cite[Section 2]{DM}.  We note that we have introduced the indicator functions $1_{K^D_a(s)}(y)$ in order to use this particular stochastic integral, which requires that the integrands be square-integrable over all of $\IR^3$.  

   Concerning the pathwise integral, we now give some details. Let $(Z(s),\ s \in [0, T])$ be a stochastic process taking its values in $L^2(\IR^d)$ satisfying
\begin{equation}\label{e3.3}
   \sup_{t \in [0,T]} E\left(\Vert Z (t) \Vert^2_{L^2(\IR^d)}\right) < \infty.
\end{equation}
Let $G:[0,T] \to \IS^\prime(\IR^d)$ be such that for any $s \in [0,T]$, ${\cal{F}} G(s)$ is a function and
\begin{equation}\label{e3.4}
   \int_0^T ds \sup_{\xi \in \IR^d} \vert {\cal{F}} G(s) (\xi)\vert^2 < \infty
\end{equation}
(notice that the fundamental solution of the wave equation satisfies this condition).
Then for any $t \in [0, T]$,
$$
   J^t_{G,Z} (x): = \int_0^t ds\, \left(G(s) \ast Z(s)\right)(x)
$$
defines a function in $L^2(\IR^d)$, a.s.\,. Moreover,
$$
   \Vert J^t_{G,Z} \Vert^2_{L^2(\IR^d)} \leq t \int_0^t ds \sup_{\xi \in \IR^d} \vert {\cal{F}} G(s)(\xi)\vert^2\, \Vert
Z(s) \Vert^2_{L^2(\IR^d)},
$$
hence
\begin{equation}\label{e3.5}
E\left(\Vert J^t_{G,Z} \Vert^2_{L^2(\IR^d)}\right) \leq t \int_0^t ds \sup_{\xi \in \IR^d} \vert {\cal{F}} G(s)(\xi) \vert^2 E\left(\Vert
Z(s) \Vert^2_{L^2(\IR^d)}\right).
\end{equation}

   For higher moments, assume that $d=3$, $(Z(s),\ s \in [0,T])$ takes values in $L^q(\IR^3)$ for some $q \in [1, \infty[$ and that $G(t)$ is
the fundamental solution of the wave equation in $\IR^3$. Then, 
%by Young's inequality (\cite{adams}, Corollary 2.25),
Fubini's theorem and the facts that $G \geq 0$ and supp~$G(s) = \partial B_s(0)$ imply that
\begin{equation}\label{e3.5a}
   \Vert G(s) \ast Z(s) \Vert_{L^q(D)}^q \leq \left(\int_{\IR^3} G(s, dy)\right)^q
        \Vert Z(s) \Vert^q_{L^q(D^s)} = s^q \Vert Z(s) \Vert^q_{L^q(D^s)}.
\end{equation}
From H\"older's inequality, it follows that
\begin{equation}\label{e3.6}
   E\left(\Vert J^t_{G,Z} \Vert^q_{L^q(D)}\right) \leq C \int_0^t ds\, \left(\int_{\IR^3} G(s, dy)\right)^q E\left(\Vert Z(s)
\Vert^q_{L^q(D^s)}\right).
\end{equation}

   The next theorem gives existence and uniqueness of the solution of equation (\ref{e3.2'}) and states some of its
properties. Notice that in comparison with Theorems 9 and 13 of \cite{DM}, we allow a non vanishing coefficient $b$ and
avoid introducing a weight function.

\begin{thm} Let $\sigma, b: \IR \to \IR$ be Lipschitz continuous functions, $v_0$, $\tilde v_0$ be real-valued functions such that $v_0$ is of class $C^2$ and $\tilde v_0 \in L^q_{\text{loc}}(\IR^3)$, for some $q \in [2, \infty[.$ Then there exists a unique process $\left(u(t) 1_{K^D_a(t)},\ t \in [0,T]\right)$ with values in $L^2(\IR^3)$ satisfying equation (\ref{e3.2'}).\\
Moreover, for any $t\in[0,T]$,
\begin{equation}
\label{e3.7}
E\left(\Vert u(t) \Vert _{L^q(K^D_a(t))}^q\right) \le C(1+ I_a(t)),
\end{equation}
% and
%\begin{equation}\label{e3.7}
%\sup_{t \in[0,T]} E\left(\Vert u(t) \Vert _{L^q(K^D(t))}^q\right) \le C(1+ I_0),
%\end{equation}
with 
\begin{equation}
\label{e3.7.tris}
   I_a(t): = \Vert v_0 \Vert^q_{L^q((K^D_a(t))^t)} + \Vert \Delta v_0 \Vert^q_{L^q((K^D_a(t))^t)} + \Vert \tilde v_0\Vert^q_{L^q((K^D_a(t))^t)}.
\end{equation}
%$$
 %  I_0: = \Vert v_0 \Vert^q_{L^q(K^D(0))} + \Vert \Delta v_0 \Vert^q_{L^q(K^D(0))} + \Vert \tilde v_0\Vert^q_{L^q(K^D(0))}.
%$$
As a consequence, 
\begin{equation}
\label{e3.7.bis}
\sup_{t\in[0,T]}E\left(\Vert u(t) \Vert _{L^q(K^D(t))}^q\right) \le C(1+ J),
\end{equation}
where $K^D(t)= K^D_1(t)$ and
$$
J=\Vert v_0 \Vert^q_{L^q(D^T)} + \Vert \Delta v_0 \Vert^q_{L^q(D^T)} + \Vert \tilde v_0\Vert^q_{L^q(D^T)}.
$$
\label{thm3.1}
\end{thm}

Before proving this theorem, we state some results concerning the initial conditions. First, we notice that under the hypotheses of the previous theorem, $J$ and
$I_a(t)$ are finite for any $a\ge 1$ and $t\in[0,T]$.

%$$
%   I_0: = \Vert v_0 \Vert^q_{L^q(\IR^3)} + \Vert \Delta v_0 \Vert^q_{L^q(\IR^3)} + \Vert \tilde v_0\Vert^q_{L^q(\IR^3)}
%< \infty.
%$$
For $v_0 \in L^2(\IR^3),$ it is well-known that
\begin{equation}\label{e3.8}
\frac{d}{dt} G(t) \ast v_0 = \frac{1}{t} \left(v_0 \ast G(t)\right) + \frac{1}{4 \pi} \int_{\vert y \vert < 1} (\Delta v_0)(\cdot+ty)dy,
\end{equation}
(see for instance \cite{sogge}).

\begin{lemma} Let $v_0, \tilde v_0$ be real-valued functions and $\cO$ be a bounded domain of $\IR^3$.

   (a) Suppose that $v_0 \in C^2(\cO^T).$ Then for any $t\in[0,T]$,
   $$
   \left\Vert \frac{d}{dt} G(t) \ast v_0 \right\Vert^q_{L^q(\cO)} \leq C \left(\Vert v_0 \Vert^q_{L^q(\cO^t)} +
\Vert \Delta v_0 \Vert^q_{L^q(\cO^t)}\right),
$$
%$$
%   \sup_{t \in [0,T]} \left\Vert \frac{d}{dt} G(t) \ast v_0 \right\Vert^q_{L^q(K^D(t))} \leq C \left(\Vert v_0 \Vert^q_{L^q(K^D(0))} +
%\Vert \Delta v_0 \Vert^q_{L^q(K^D(0))}\right),
%$$
for any $q \in [2, \infty[.$

   (b) Assume that $\tilde v_0 \in L^q(\cO^T)$ for some $q \in [2, \infty[.$ Then for any $t\in[0,T]$,
$$
\Vert G(t) \ast \tilde v_0 \Vert_{L^q(\cO)} \leq \Vert \tilde v_0 \Vert_{L^q(\cO^t)}.
$$
%$$
%\sup_{t \in [0, T]} \Vert G(t) \ast \tilde v_0 \Vert_{L^q(K^D(t))} \leq \Vert \tilde v_0 \Vert_{L^q(K^D(0))}.
%$$
\label{lem3.2}
\end{lemma}

\noindent{\it Proof}:  By (\ref{e3.8}),
$$
   \left\Vert \frac{d}{dt} G(t) \ast v_0\right\Vert^q_{L^q(\cO)} \leq C \left(A(t) + B(t)\right),
$$
where
$$
   A(t) = \left\Vert \frac{1}{t}\left(v_0 \ast G(t)\right)\right \Vert^q_{L^q(\cO)}, \quad
   B(t) = \left\Vert \frac{1}{t^3} \left(\Delta v_0 \ast 1_{B_t(0)}\right)\right \Vert^q_{L^q(\cO)}.
% \int_{\vert y \vert < 1} (\Delta v_0)(\cdot + ty)dy \Vert^q_{L^q(\IR^3)}.
$$
Inequality (\ref{e3.5a})  yields
$$
   A(t) \leq \frac{1}{t^q} \left(\int_{\IR^3}G(t,du)\right)^{q} \Vert v_0 \Vert^q_{L^q(\cO^t)} = \Vert v_0 \Vert^q_{L^q(\cO^t)}.
$$
Similarly,
$$
   B(t) \leq C \Vert \Delta v_0\Vert^q_{L^q(\cO^t)},
\quad\mbox{and}\quad 
  \Vert G(t) \ast \tilde v_0 \Vert^q_{L^q(\cO)} \leq \Vert \tilde v_0 \Vert^q_{L^q(\cO^t)}.
$$
This proves the lemma.
\hfill $\Box$
\vskip 16pt

\noindent{\em Proof of Theorem \ref{thm3.1}.} Consider the Picard iteration scheme given by
$$
   u^0(t, \cdot) = \frac{d}{dt} G(t) \ast v_0 + G(t) \ast \tilde v_0,
$$
and, for $n \geq 1 $,
\begin{align*}
u^n(t, \cdot)1_{K^D_a(t)}(\cdot)&= 1_{K^D_a(t)}(\cdot)\Big( u^0(t, \cdot) + %1_{K^D(t)}(\cdot) 
 \int^t_0 \int_{\IR^3} G(t-s, \cdot -y) \sigma\left(u^{n-1}(s, y)\right)\\
 \\
& \qquad \qquad\qquad \qquad \times 1_{K^D_a(s)}(y) M(ds, dy)\big)\\
\\
&\qquad + 1_{K^D_a(t)}(\cdot) \int_0^t ds\, G(t-s) \ast \left(b\left(u^{n-1} (s, \cdot)\right)1_{K^D_a(s)}(\cdot)\right)(\cdot) .
\end{align*}

   For $n \geq 0$, set $v^n(t, \cdot) = u^n(t, \cdot) 1_{K^D_a(t)} (\cdot)$.
  % and
%$$
%  J(t-s) = \sup_{\xi \in \IR^3} \int_{\IR^3} \mu (d\eta)\, \vert {\cal{F}} G(t-s) (\xi-\eta)\vert^2.
%$$  
%Notice that $\sup_{0\le s<t\le T} J(t-s) <\infty$.
Lemma \ref{lem3.2} with $\cO:= K^D_a(t)$, along with the inequalities (\ref{p16star}) and (\ref{e3.6}) applied to $Z(s):= g\left(u^{n-1}(s, \cdot)\right) 1_{K^D_a(s)}(\cdot)$
with $g \equiv \sigma$ and $g \equiv b,$ respectively, and the Lipschitz properties of $\sigma$ and $b$, tell us that for any $n \geq 1$,
\begin{align*}
   E\left(\Vert v^n (t) \Vert^q_{L^q(K^D_a(t))}\right) &\leq C I_a(t) + C_1
\int_0^t ds\, E\Big( \int_{\IR^3} dx\, 1_{K^D_a(s)}(x) (1+ \vert u^{n-1}(s, x)\vert^q)\Big)\\
   &\leq C I_a(t) +C_1 \int_0^t ds\, \left(1+ E\left(\Vert v^{n-1}(s)\Vert^q_{L^q(K^D_a(s))}\right)\right).
\end{align*}
Notice that $v^0(t, \cdot)1_{K^D_a(t)}(\cdot) = u^0(t, \cdot)1_{K^D_a(t)}(\cdot)$ and by Lemma \ref{lem3.2},
$$
   E\left(\Vert v^0(t) \Vert^q_{L^q(K^D_a(t))}\right) = \left(\Vert u^0(t) \Vert^q_{L^q(K^D_a(t))}\right) \leq C
\ I_a(t) < \infty.
$$
%$$
 %  \sup_{t \in [0,T]} E\left(\Vert v^0(t) \Vert^q_{L^q(K^D(t))}\right) = \sup_{t \in [0,T]} \left(\Vert u^0(t) \Vert^q_{L^q(K^D(t))}\right) \leq C
%\ I_0 < \infty.
%$$

Thus, we obtain existence, uniqueness and (\ref{e3.7}) using arguments similar to those in the proof of Theorem 9 in \cite{DM}, based on Gronwall's lemma.
\smallskip

   A priori, the solution should be written $(u_a(t)1_{K^D_a(t)},\ t \in [0,T])$. However, $u_a(t)$ does not depend on $a$. Indeed, for $1\leq a \leq b$, both $(u_a(t)1_{K^D_a(t)},\ t \in [0,T])$ and $((u_b(t)1_{K^D_b(t)}) 1_{K^D_a(t)},\ t \in [0,T])$ satisfy (\ref{e3.2'}), so uniqueness for $a$ implies that $u_a(t)1_{K^D_a(t)} = u_b(t)1_{K^D_a(t)}$.
\smallskip

We now prove (\ref{e3.7.bis}).
For any $m\ge 1$, set $D_m= \{x\in D: \overline{B_{1/m}(x)}\subset D\}$ and $a_m:= 1+\frac{1}{m^2}$. We notice that for any $m\ge T-t$,
and $t\in[0,T]$,
the set $(K^{D_m}_{a_m}(t))^t$ is included in $D^T$. Consequently, $I_{a_m}(t)\le J$. Moreover, the sequence of sets $K^{D_m}_{a_m}(t)$,
$m\ge 2(T-t)$, is increasing and $\cup_m K^{D_m}_{a_m}(t)= K^D(t)$. Hence, (\ref{e3.7.bis}) follows from (\ref{e3.7}) by letting $m$ tend to $\infty$.
\hfill $\Box$
\vskip 16pt

\subsection{Regularity in the space variable}
\label{ss2.2}
   The solution of equation (\ref{e3.2'}) can be approximated by a sequence of solutions of similar equations obtained by regularising the fundamental solution $G$ in the term involving a stochastic integral. More precisely, let $G_n(t,x)$, $n \geq 1$, be as in (\ref{A.1}), fix a bounded domain
$D \subset \IR^3$ satisfying the conditions described at the beginning of the section and set, for any $a > 1$ and $n$ with $1+\frac{1}{n} < a$,
\begin{align}
&u_n(t,x)1_{K^D_a(t)}(x)= 1_{K^D_a(t)}(x)\left(\frac{d}{dt} G(t) \ast v_0 + G(t) \ast \tilde v_0\right) (x)\nonumber\\
&\qquad+ 1_{K^D_a(t)}(x) \int_0^t \int_{\IR^3} G_n(t-s, x-y) \sigma\left(u_n(s, y)\right) 1_{K^D_a(s)}(y) M(ds, dy)\nonumber\\
&\qquad+ 1_{K^D_a(t)}(x) \int^t_0 ds\, \left(G(t-s) \ast \left(b\left(u_n(s, \cdot)\right) 1_{K^D_a(s)}(\cdot)\right)\right)(x).
\label{e3.10}
\end{align}
Since $G_n$ is smooth, the stochastic integral in (\ref{e3.10}) is considered in Walsh's sense \cite{walsh}. In particular, we could remove the $1_{K^D_a(t)}(x)$ and $1_{K^D_a(s)}(y)$ from (\ref{e3.10}).

\begin{prop} Let $\sigma$, $b$, $v_0$ and $\tilde v_0$ be as in Theorem \ref{thm3.1}. %Assume that $v_0$ is
% of class $C^2$ and $\tilde v_0 \in L^q(\IR^3)$, for any $q \in [2,\infty[$. 
Then there exists a unique process $\{u_n(t, \cdot)1_{K^D_a(t)}(\cdot),\ t
\in[0,T]\}$ solution of (\ref{e3.10}), and this process is such that for $q \in [2,\infty[$,
\begin{equation}\label{e3.11}
\sup_{n\geq (a-1)^{-1}}\, \sup_{t\in[0,T]} E\left(\Vert u_n(t)\Vert^q_{L^q(K^D_a(t))}\right) < \infty.
\end{equation}
Moreover
\begin{equation}\label{e3.12}
\lim_{n\to \infty} \sup_{t \in [0,T]}E\left(\Vert u_n(t) - u(t)\Vert^q_{L^q(K^D_a(t))}\right) = 0.
\end{equation}
\label{prop3.3}
\end{prop}

\noindent{\it Proof}:  The proof of existence, uniqueness and (\ref{e3.11}) can be obtained using arguments similar to those applied in Theorem \ref{thm3.1}, taking into account (\ref{A.4}). 

   To prove (\ref{e3.12}), we apply the $L^p$-estimates of
the stochastic and the pathwise integrals given in (\ref{p16star}) and (\ref{e3.6}), respectively (with $G$ replaced by $G_n$). Let $n\ge 1$ be such that $1+\frac{1}{n}\le a$; then  
\begin{equation}\label{KDT}
   (K^D_a(t))^{(t-s)(1+1/n)} \subset K^D_a(s),
\end{equation}
for $0\leq s \leq t\leq T$ and we obtain
\begin{align} \nonumber
& E\left(\Vert u_n(t) -u(t) \Vert^q_{L^q(K^D_a(t))}\right)\\
%&\qquad \leq C \Big(\int_0^t ds\, \left(J(t-s)+ 1\right) E\left(\Vert\left(u_n(s)-u(s)\right)1_{K^D(s)}\Vert^q_{L^q(\IR^3)}\right)\\
%&\qquad +M_n(t)\\
&\qquad \leq C \int_0^t ds\, E\left(\Vert\left(u_n(s)-u(s)\right) \Vert^q_{L^q(K^D_a(s))}\right)+M_n(t),
\label{e3.121}
\end{align}
where
%\begin{align*}
%M_n(t) &= E\Big( \Big\Vert\int_0^t\int_{\IR^3}\big(G_n(t-s,\cdot-t)-G(t-s,\cdot-t)\big)\\
%&\quad\times \sigma\big(u(s,y)\big)1_{K^D(s)}M(ds,dy)\Big\Vert_{L^q(\IR^3)}^q\Big).
%\end{align*}
\begin{align*}
M_n(t) &= E\Big( \Big\Vert\int_0^t\int_{\IR^3}\big(G_n(t-s,\cdot-y)-G(t-s,\cdot-y)\big)\\
&\qquad \qquad\times \sigma\big(u(s,y)\big)1_{K^D_a(s)}(y)M(ds,dy)\Big\Vert_{L^q(K^D_a(t))}^q\Big).\\
&= E(\Vert v_{G_n,Z} - v_{G,Z} \Vert^q_{L^q(K^D_a(t))}),
\end{align*}
and $Z(s,y) = \sigma(u(s,y)) 1_{K^D_a(s)}(y).$

   We shall show that
\begin{equation}\label{G1}
   \lim_{n\to \infty} M_n(t) = 0.
\end{equation}
Since (\ref{e3.12}) follows from this property, (\ref{e3.121}) and Gronwall's lemma, the proposition will be proved.

   To prove (\ref{G1}), we check first the case $q=2$, that is,
\begin{equation}\label{G2}
   \lim_{n \to \infty} E(\Vert v_{G_n,Z} - v_{G,Z} \Vert^2_{L^2(K^D_a(t))}) = 0.
\end{equation}
Indeed,
\begin{eqnarray*}
   E(\Vert v_{G_n,Z} - v_{G,Z} \Vert^2_{L^2(K^D_a(t))}) &=& E(\Vert v_{G_n-G,Z} \Vert^2_{L^2(K^D_a(t))})\\
&\leq& E(\Vert v_{G_n-G,Z} \Vert^2_{L^2(\IR^3)}),
\end{eqnarray*}
and by the isometry property (\ref{isometry1}) of the stochastic integral, this is equal to
$$
   \int_0^t ds \int_{\IR^3} d\xi\, E(\vert \F Z(t-s)(\xi)\vert^2) \int_{\IR^3} \mu(d\eta)\, \vert \F(G_n(t-s) - G(t-s)) (\xi-\eta)\vert^2.
$$
Notice that the integrand converges to 0 pointwise, and therefore the integral converges to 0 by the dominated convergence theorem, since the integrand is bounded by
$$
 4 E(\vert \F Z(t-s)(\xi)\vert^2) \vert \F G(t-s) (\xi)\vert^2
$$
and by the Lipschitz property of $\sigma$, (\ref{J}) and (\ref{e3.7}),
\begin{align*}
   &\int^t_0 ds \int_{\IR^3} d \xi\, E(\vert \F Z(t-s)(\xi)\vert^2) \int_{\IR^3} \mu(d\eta) \, \vert \F G(t-s) (\xi-\eta) \vert^2\\
   &\qquad\leq \int^t_0 ds\, E(\Vert Z(t-s)\Vert^2_{L^2(\IR^3)}) \sup_{\xi\in \IR^3} \int_{\IR^3} \mu(d \eta)\, \vert \F G(t-s)(\xi-\eta)\vert^2\\
   &\qquad\leq C \int_0^t ds\, (1+ E(\Vert u(t-s) 1_{K^D_a(t-s)} \Vert^2_{L^2(\IR))})\\
\\
   &\qquad< + \infty.
\end{align*}
This establishes (\ref{G2}).

   In order to deduce (\ref{G1}) from (\ref{G2}), we use the fact that $L^2$-convergence (on $\Omega \times K^D_a(t))$ together with $L^q$-boundedness (for all $q$), implies $L^q$-convergence. In particular, it now suffices to show that
\begin{equation}\label{G3}
   \sup_{n> (a-1)^{-1}} E(\Vert v^t_{G_n,Z} \Vert^q_{L^q(K^D_a(t))}) < + \infty.
\end{equation}
%and
%\begin{equation}\label{G4}
%   E(\Vert v^t_{G,Z} \Vert^q_{L^q(K^D(t))}) < + \infty.
%\end{equation}

   By (\ref{2.4.2}) with $\cO_m$ replaced by $K^D_a(t)$, we obtain using (\ref{KDT}) that
\begin{equation}\label{G5}
   E(\Vert v^t_{G_n,Z} \Vert^q_{L^q(K^D_a(t))}) \leq C \int_0^t ds\, E(\Vert Z(s)\Vert^q_{L^q(K^D_a(s))}).
\end{equation}
%Similarly, (\ref{p16star}) gives the same bound for the expectation in (\ref{G4}).
 The right-hand side of (\ref{G5}) is finite by the Lipschitz property of $\sigma$ and (\ref{e3.7}). The proposition is proved. 
\hfill $\Box$
\vskip 16pt

%By dominated convergence, $\lim_{n\to\infty} M_n(t)=0$. Indeed, let $q=2p$; by the arguments in (\ref{2.4.1}) we have 
%\begin{align*}
%M_n(t) &= C \int_0^t ds E\left(\int_{\IR^3} d\xi\, \vert{\cal{F}}(\vert\sigma(u(t-s,\cdot))1_{K^D(t-s)}\vert^p)(\xi)\vert^2\right)\\
%&\quad\times \int_{\IR^3} \mu(d\eta) \vert {\cal{F}}(G(s)-G_n(s))(\xi- \eta)\vert^2,
%\end{align*}
%and by (\ref{e3.7}),
%\begin{align*}
%&\int_0^t ds E\left(\int_{\IR^3} d\xi\, \vert{\cal{F}}(\vert\sigma(u(t-s,\cdot))1_{K^D(t-s)}\vert^p)(\xi)\vert^2\right)
% \int_{\IR^3} \mu(d\eta) \vert {\cal{F}} G(s)(\xi- \eta)\vert^2\\
% &\qquad \le C\int_0^t ds\, \left(1+ E\left(\Vert u(t-s)1_{K^D(t-s)}\Vert^q_{L^q(\IR^3)}\right)\right)\\
% &\qquad\qquad\times\sup_{\xi \in \IR^3}
%\int_{\IR^3} \mu(d\eta) \vert {\cal{F}}(G(s)(\xi- \eta)\vert^2\\
% &\qquad < \infty.
%\end{align*}
%&\qquad + \int_0^t ds\, (1+ E(\Vert u(s)1_{K^D(s)}\Vert^q_{L^q(\IR^3)})\\
%&\qquad\quad\times  \sup_{\xi \in \IR^3}
%\int_{\IR^3} \mu(d\eta) \vert {\cal{F}}(G(s)-G_n(s))(\xi- \eta)\vert^2.
%\end{align*}
%By (\ref{e3.7}) and dominated convergence CAREFUL \ref{?}: this has to be used on the isometry (\ref{isometry1}),
%before taking sup??, the last integral tends to zero as $n$ tends to infinity, uniformly in $t$. 
%Therefore (\ref{e3.12}) follows from Gronwall's lemma.
%\hfill $\Box$
%\vskip 16pt

   The next lemma completes the study of the initial conditions needed in this section.

   For a given real function $v$ defined on a bounded domain $\cO\subset\IR^3$ and for $\gamma\in]0,1[$, we denote by $\Vert v\Vert_{\mathcal{C}^\gamma(\cO)}$ the $\gamma$-H\"older semi-norm on $\cO$, that is
$$
   \Vert v\Vert_{\mathcal{C}^\gamma(\cO)} = \sup_{x,y\in\cO; x\neq y}\frac{|v(x)-v(y)|}{|x-y|^\gamma}.
$$

\begin{lemma} Let $v_0, \tilde v_0$ be real-valued functions. Assume that $v_0 \in C^2(\IR^3)$, $\Delta v_0$ is $\gamma_1$-H\"{o}lder continuous and $\tilde v_0$ is $\gamma_2$-H\"{o}lder continuous for some fixed $\gamma_1, \gamma_2 \in \,]0,1[.$ Then, for any $q \in [2, \infty[$, $t\in[0,T]$ and any bounded domain $\cO$, there is a positive constant $C$ depending on $q$, $t$ and $\cO$ such that for every $\rho_1 \in\, ]0, \gamma_1[$ and $\rho_2 \in\, ]0,
\gamma_2[$,
$$
 \left\Vert \frac{d}{dt}G(t) \ast v_0 \right\Vert_{\rho_1,q,\cO} \leq C\left(%\sup_{z\in\cO^t}|\nabla v_0(z)\vert^q
 \Vert v_0\Vert^q_{\rho_1,q,\cO^t}
 +\Vert \Delta v_0\Vert^q_{{\mathcal{C}^{\gamma_1}(\cO^t)}}\right)
$$
and
$$
 \left\Vert G(t) \ast \tilde v_0 \right\Vert_{\rho_2,q,\cO} \leq C \Vert \tilde v_0\Vert^q_{\mathcal{C}^{\gamma_2}(\cO^t)}.
$$
\label{lem3.4}
\end{lemma}

%\begin{remark} \label{rem3.1} Lemma \ref{lem3.2}, the previous lemma along with condition (S) yields
%$$
%\Vert \frac{d}{dt} G(t)*v_0+G(t)*\tilde v_0\Vert_{W^{\gamma,q}(K^D(t))} < \infty,
%$$
%for any $q\in[2,\infty[$, $\gamma\in]0,\gamma_1\wedge\gamma_2]$.
%\end{remark}

\noindent{\it Proof of Lemma \ref{lem3.4}}:  We start by checking the first inequality using the basic identity (\ref{e3.8}). Set $A(x) = \frac{1}{t}(v_0
\ast G(t))(x)$.  H\"{o}lder's inequality with respect to the measure $G(t,du)$ and the properties of $v_0$
 yield
\begin{align*}
\Vert A(\cdot)\Vert^q_{\rho,q,\cO}
%\displaystyle\int_{\IR^3} dx  \int_{\IR^3} dy\, \frac{\vert A(x)-A(y)\vert^q}{\vert x-y \vert^{\rho q}}\\
&\leq \frac{1}{t^q} \left(\int_{\IR^3} G(t,du)\right)^{q-1}
\int_{\IR^3} G(t, du)\, \Vert v_0(\cdot - u)\Vert^q_{\rho,q,\cO}\\
%\int_{\IR^3} dx \int_{\IR^3} dy\, \frac{\vert v_0(x-u)-v_0(y-u)\vert^q}{\vert x-y \vert^{\rho q} }\\
&\leq \frac{1}{t^q}\left(\int_{\IR^3} G(t, du)\right)^q \Vert v_0\Vert^q_{\rho_1,q,\cO^t}
%  \sup_{z \in \cO^t} \vert \nabla v_0(z) \vert^q 
\\
%&\qquad \times \int_{\cO} dx \int_{\cO} dy\,  \vert x-y \vert^{q(1-\rho-\frac{3}{q})} \\
&= C \Vert v_0\Vert^q_{\rho_1,q,\cO^t},
%\sup_{z \in \cO^t} \vert \nabla v_0(z) \vert^q,
\end{align*}
for any $\rho \in\, ]0, 1[.$

   Let us now consider $B(x) = \int_{\vert z\vert \leq 1} \Delta v_0(x+tz)dz.$ The H\"{o}lder-continuity property of $\Delta v_0$ implies that
\begin{align*}
   \Vert B(\cdot) \Vert^q_{\rho,q,\cO}
 & \leq \Vert \Delta v_0\Vert^q_{\mathcal{C}^{\gamma_1}(\cO^t)} 
 \int_{\vert z \vert \le 1} dz \int_{\cO} dx \int_{\cO} dy\, \vert x-y \vert^{(\gamma_1-\rho-\frac{3}{q})q}\\
  &\leq C\, \Vert \Delta v_0\Vert^q_{\mathcal{C}^{\gamma_1}(\cO^t)},
\end{align*}
if $\rho \in\,]0, \gamma_1[.$ This ends the proof of the first inequality in the lemma.

   We now prove the second inequality. Set $R(x) = (G(t) \ast \tilde v_0)(x).$ As for $A(x)$, H\"older's inequality along with the H\"{o}lder-continuity of $\tilde
v_0$ yield
\begin{align*}
   \Vert R(\cdot) \Vert^q_{\rho,q,\cO}
%   \int_{\IR^3} dx \int_{\IR^3} dy\, \frac{\vert C(x) - C(y) \vert^q}{\vert x-y \vert^{\rho q}} 
  & \leq \Vert \tilde v_0\Vert^q_{\mathcal{C}^{\gamma_2}(\cO^t)}  \int_{\cO} dx \int_{\cO} dy\, \vert x-y \vert^{(\gamma_2-\rho-\frac{3}{q})q}\\
  & \leq C \Vert \tilde v_0\Vert^q_{\mathcal{C}^{\gamma_2}(\cO^t)}, 
\end{align*}
if $\rho \in\, ]0, \gamma_2[.$ The proof of the lemma is now complete.
 \hfill $\Box$

\begin{remark} \label{rem3.2} 
Assume that the hypotheses of Lemma \ref{lem3.4} are satisfied. Then those of Lemma \ref{lem3.2} hold for any $q \in
[2, \infty[.$ Consequently, the solution of the deterministic wave equation in dimension $d=3$,
$$
\left(\frac{\partial^2}{\partial t^2}-\Delta\right) w(t,x)= 0,$$
 with initial conditions $v_0,
\tilde v_0$, is a function $t \in[0,T]\mapsto w(t)$, such that each $w(t)$ takes values in $W^{\gamma, q}(\cO)$, for any $q \in
[2, \infty[$, $\gamma \in\, ]0, \gamma_1 \wedge \gamma_2[$ and every bounded domain $\cO$ . 
In other words, $w(t)$ belongs to $W^{\gamma,q}_{\text{loc}}(\IR^3)$.

Indeed, 
$$w(t)=\frac{d}{dt} G(t)*v_0+G(t)*\tilde v_0$$
and the above-mentioned two lemmas establish that
$$\sup_{t\in[0,T]}\left\Vert \frac{d}{dt} G(t)*v_0+G(t)*\tilde v_0\right\Vert_{W^{\gamma,q}(\cO)} < \infty.$$

\end{remark}

   We are now in a position to give the sample path behaviour in the $x$-variable of the solution of (\ref{e3.2'}). In the next statement, $a\ge 1$ is a real number.

\begin{thm} Fix $q\in[2,\infty[$ and assume that $\tau(\beta,\delta):=\frac{2-\beta}{2}\wedge\frac{1+\delta}{2} >
\frac{3}{q}$. Suppose also that

   (a) $\sigma$ and $b$ are Lipschitz continuous functions,

   (b) $v_0$, $\tilde v_0$ are real-valued functions, $v_0 \in C^2(\IR^3)$, and
$\Delta v_0$ and $\tilde v_0$ are H\"{o}lder-continuous functions of order $\gamma_1, \gamma_2 \in \, ]0,1[$,
respectively.

Set 
\begin{equation}
\label{picard0}
u^0(t,x) = \left(\frac{d}{dt} G(t) \ast v_0 + G(t) \ast \tilde v_0\right)(x).
\end{equation}
\noindent Then, for any $t\in[0,T]$,
\begin{equation}\label{e3.13}
   \sup_{n \geq 1}\  E\left(\Vert u_n(t)  \Vert^q_{W^{\gamma,q}(K^D_a(t))}\right) \le 
C\left(1+
%\sup_{t\in[0,T]}\Vert u^0(t)1_{K^D(t)}\Vert^q_{W^{\gamma,q}(K^D(t))}
\Vert u^0(t) \Vert^q_{W^{\gamma,q}(K^D_a(t))} \right),
\end{equation}
for any $\gamma \in\, ]0, \gamma_1 \wedge \gamma_2 \wedge (\tau(\beta,\delta) - \frac{3}{q})[$, and
\begin{equation}\label{e3.14}
   E\left(\Vert u(t) \Vert^q_{W^{\gamma,q}(K^D_a(t))}\right) \le C\left(1+\Vert u^0(t) \Vert^q_{W^{\gamma,q}(K^D_a(t))}
%\sup_{t \in [0,T]}\left\Vert u^0(t)1_{K^D(t)}\right\Vert^q _{W^{\gamma,q}(K^D(t))}
\right).
\end{equation}
%with 
%\begin{equation}
%\label{ia}
%   \tilde I^{\gamma,D}_a(t) = \Vert u^0(t) \Vert^q_{W^{\gamma,q}(K^D_a(t))}.
%\end{equation}
As a consequence, 
\begin{equation}\label{e3.14.bis}
 \sup_{t \in [0,T]}  E\left(\Vert u(t) \Vert^q_{W^{\gamma,q}(K^D(t))}\right) \le C\left(1+\tilde J
%\sup_{t \in [0,T]}\left\Vert u^0(t)1_{K^D(t)}\right\Vert^q _{W^{\gamma,q}(K^D(t))}
\right),
\end{equation}
where $K^D(t)=K^D_1(t)$ and
\begin{align*}
 \tilde J
 &= \Vert v_0 \Vert^q_{W^{\gamma,q}(D^T)} + \Vert \Delta v_0 \Vert^q_{L^q(D^T)} + \Vert \tilde v_0\Vert^q_{L^q(D^T)}\\
 &%+ \sup_{z\in D^T} |\nabla v_0(z)|^q 
 \qquad + \Vert \Delta v_0\Vert^q_{\mathcal{C}^{\gamma_1}(D^T)} + \Vert \tilde v_0\Vert^q_{\mathcal{C}^{\gamma_2}(D^T)}. 
 \end{align*}
%\end{description}

%%%%%%%
%$$
 %  \tilde I_0 = \sup_{t\in[0,T]}\Vert u^0(t) \Vert^q_{W^{\gamma,q}(K^D(t))}.
%$$
%\noindent Then
%\begin{equation}\label{e3.13}
 %  \sup_{n \geq 1}\ \sup_{t \in [0,T]} E\left(\Vert u_n(t)  \Vert^q_{W^{\gamma,q}(K^D(t)}\right) \le 
%C\left(1+
%%\sup_{t\in[0,T]}\Vert u^0(t)1_{K^D(t)}\Vert^q_{W^{\gamma,q}(K^D(t))}
%\tilde I_0 \right),
%\end{equation}
%for any $\gamma \in\, ]0, \gamma_1 \wedge \gamma_2 \wedge (\tau(\beta,\delta) - \frac{3}{q})[$, and
%\begin{equation}\label{e3.14}
 % \sup_{t \in [0,T]} E\left(\Vert u(t) \Vert^q_{W^{\gamma,q}(K^D(t)}\right) \le C\left(1+\tilde I_0
%%%\sup_{t \in [0,T]}\left\Vert u^0(t)1_{K^D(t)}\right\Vert^q _{W^{\gamma,q}(K^D(t))}
%\right).
%\end{equation}

%%%%%%
In particular, for every $\gamma \in\, ]0, \gamma_1 \wedge \gamma_2 \wedge\tau(\beta,\delta)[$,
there exists $q\in\, ]1,\infty[$ large enough so that
\begin{equation}\label{e3.15}
 \sup_{t \in[0,T]} E\left(\Vert u(t) \Vert^q_{W^{\gamma, q}(K^D(t))}\right) \le C\left(1+\tilde J
 %\sup_{t \in [0,T]}\left\Vert u^0(t)1_{K^D(t)}\right\Vert^q _{W^{\gamma,q}(K^D(t))}
 \right). 
\end{equation}
\label{thm3.6}
\end{thm}

\noindent{\it Proof}: For $n \geq 1$, consider the Picard iteration scheme corresponding to equation (\ref{e3.10}). That is,
$u^0_n(t,x) = u^0(t,x)$, defined in (\ref{picard0})
and, for $m \geq 1$,
\begin{align*} 
u^m_n(t, x) 1_{K^D_a(t)}(x)
&= u^0_n(t,x)1_{K^D_a(t)}(x)\\
& + 1_{K^D_a(t)}(x) \int^t_0 \int_{\IR^3} G_n(t-s, x-y)\\
&\quad\times \sigma\left(u^{m-1}_n(s,y)\right) 1_{K^D_a(s)} (y) M(ds,dy)\\
&+ 1_{K^D_a(t)}(x) \int_0^t ds\, \left(G(t-s) \ast b\left(u_n^{m-1}(s, \cdot)\right)1_{K^D_a(s)}(\cdot)\right) (x).
\end{align*}

%Notice that, in the above equation, if $x\in K^D(t)$ then $y\in K^D(s)$ for $0<s\le t$.
Set
$$
   R^{m, \gamma,D}_{n} (t) = E\left(\Vert u^m_n(t) \Vert^q_{W^{\gamma,q}(K^D_a(t))}\right).
%\int_{\IR^3} dx \int_{\IR^3} dy\, \frac{\vert u^m_n(t,x)-u^m_n(t,y)\vert^q}{\vert x-y\vert^{\rho q}} 1_{K^D(t)}(x)1_{K^D(t)}(y)
 $$

   Clearly, for $m \geq 1$, $R^{m,\gamma}_{n,D}(t) \leq C(\tilde I^{\gamma, D}_a(t) + T^{m, \gamma,D,1}_{n}(t) + T_{n}^{m,\gamma,D,2}(t)),$
with
\begin{align*}
%T^\gamma_D(t) &= \Vert u^0(t) \Vert^q_{W^{\gamma,q}(K^D_a(t))},\\
T^{m,\gamma,D,1}_{n}(t) &= E\left(\Vert v^t_{G_n,\sigma(u^{m-1}_n)1_{K^D_a}}\Vert^q_{W^{\gamma,q}(K^D_a(t))}\right),\\
T^{m,\gamma,D,2}_{n}(t) &= E\left(\Vert J^t_{G,b(u^{m-1}_n )1_{K^D_a}} \Vert^q_{W^{\gamma,q}(K^D_a(t))}\right)
\end{align*}
and $\tilde I^{\gamma, D}_a(t) = \Vert u^0(t) \Vert^q_{W^{\gamma,q}(K^D_a(t))}$.
It follows from  Remark \ref{rem3.2} that for  $\gamma \in\, ]0, \gamma_1 \wedge \gamma_2[$ and $a\ge 1$,
\begin{equation}\label{e3.17}
\sup_{t \in [0,T]} \tilde I^{\gamma,D}(t) < \infty.
\end{equation}

Let $g$ be a real-valued Lipschitz function. We shall show that the process 
$$\{g(u_n^{m-1}(t, \cdot)) 1_{K^D_a(t)}(\cdot), \ t \in [0,T]\}$$ 
satisfies
the following properties:
\begin{description}
\item{(i)} $\sup_{n,m \geq 1}\, \sup_{t \in [0,T]} E\left(\Vert g(u_n^{m-1}(t)) \Vert^q_{L^q(K^D_a(t))} \right) <
\infty$,
\item{(ii)} for any $0\le s<t\le T$,
\begin{align*}
&E\left(\Vert g(u_n^{m-1}(s))1_{K^D_a(s)}\Vert^q_{\gamma,q,(K^D_a(t))^{(t-s)(1+\frac{1}{n})}}\right)\\
&\qquad \le C 
E\left(\Vert u_n^{m-1}(s) \Vert^q_{\gamma,q,K^D_a(s)}\right),
\end{align*}
\item{(iii)} for any $0\le s<t\le T$, $|y|\le t-s$,
\begin{align*}
&E\left(\Vert g(u_n^{m-1}(s,\cdot-y))1_{K^D_a(s)}(\cdot-y)\Vert^q_{\gamma,q,K^D_a(t)}\right)\\
&\qquad \le C 
E\left(\Vert u_n^{m-1}(s)\Vert^q_{\gamma,q,K^D_a(s)}\right).
\end{align*}
\end{description}
Indeed, $g$ has linear growth and using the same proof as in Theorem \ref{thm3.1} with $G$ replaced by $G_n$, we
see, as in (\ref{e3.7}), that
$$
   \sup_{n,m \geq 1}\, \sup_{t \in [0,T]}E\left (\Vert u^m_n(t) \Vert^q_{L^q(K^D_a(t))}\right) < \infty,
$$
which proves (i).

   To prove (ii), we use (\ref{KDT}) and the Lipschitz property of $g$ to see that
\begin{align*}
&E\left(\Vert g(u_n^{m-1}(s))1_{K^D_a(s)}\Vert^q_{\gamma,q,(K^D_a(t))^{(t-s)(1+\frac{1}{n})}}\right)\\
&\quad = E\left(\Vert g(u_n^{m-1}(s))\Vert^q_{\gamma,q,(K^D_a(t))^{(t-s)(1+\frac{1}{n})}}\right)\\
&\quad \le C E\left(\Vert u_n^{m-1}(s)\Vert^q_{\gamma,q,(K^D_a(t))^{(t-s)(1+\frac{1}{n})}}\right)\\
&\quad = C E\left(\Vert u_n^{m-1}(s)1_{K^D_a(s)}\Vert^q_{\gamma,q,(K^D_a(t))^{(t-s)(1+\frac{1}{n})}}\right)\\
&\quad \le C E\left(\Vert u_n^{m-1}(s)1_{K^D_a(s)}\Vert^q_{\gamma,q,K^D_a(s)}\right).
\end{align*}

   To prove (iii), let $0\le s<t\le T$ and $|y|\le t-s$. A change of variables yields
$$
\Vert g(u_n^{m-1}(s,\cdot-y))1_{K^D_a(s)}(\cdot-y)\Vert^q_{\gamma,q,K^D_a(t)} \le \Vert g(u_n^{m-1}(s))1_{K^D_a(s)}(\cdot)\Vert^q_{\gamma,q,(K^D_a(t))^{t-s}}.
$$
Then (iii) follows from the set inclusion $(K^D_a(t))^{t-s}\subset K^D_a(s)$ and the Lipschitz property of $g$, in a similar manner as (ii).

   Lemma \ref{l2.2}, and more specifically (\ref{222'}), applied to the stochastic process 
$Z(s)=\sigma(u_n^{m-1}(s))1_{K^D_a(s)}$, $s\in[0,T]$, $\cO = K^D_a(t)$ and $\rho:=\gamma$
yields
\begin{equation*}
   T^{m,\gamma,D,1}_{n}(t) \leq C\,  \int_0^t ds\, E\left(\Vert\sigma(u_n^{m-1}(s))1_{K^D_a(s)}\Vert^q_{W^{\gamma,q}((K^D_a(t))^{(t-s)(1+\frac{1}{n})}}\right).
\end{equation*}
By the properties (i), (ii) proved above with $g\equiv\sigma$, we obtain
\begin{equation}
\label{e3.18}
 T^{m,\gamma,1}_{n,D}(t) \leq C_1 + C_2 \int_0^t ds\, R_{n}^{m-1,\gamma,D}(s).
 \end{equation}

 H\"older's inequality and  property (iii) with $g\equiv b$ imply
\begin{align*}
 &E\left(\Vert J^t_{G,b(u_n^{m-1}) 1_{K^D_a}}\Vert^q_{\gamma,q,K^D_a(t)}\right)\\
 &\quad\le C \int_0^t ds \int_{\IR^3} G(t-s,dy) E\left(\Vert b(u_n^{m-1}(s,\cdot-y))\,1_{K^D_a(s)}(\cdot-y)\Vert^q_{\gamma,q,K^D_a(t)}\right)\\
 &\quad\le C \int_0^t ds E\left(\Vert u_n^{m-1}(s)\Vert^q_{\gamma,q,K^D_a(s)}\right).
\end{align*}
Thus, owing to (i), we have
\begin{equation}
 \label{e3.19}
   T^{m,\gamma,D,2}_{n} (t) \leq C_1 + C_2 \int_0^t ds \, R_{n}^{m-1,\gamma,D}(s).
\end{equation} 

   Hence, the sequence of functions $\left(R_{n}^{m, \gamma,D}(t),\ t \in [0,T],\ m \geq 1\right)$, $\gamma \in\, ]0,
\gamma_1 \wedge \gamma_2 \wedge (\tau(\beta,\delta)- \frac{3}{q})[$, satisfies the inequality
$$
  R^{m,\gamma,D}_{n}(t)\leq C \left(1 + \tilde I_a^{\gamma,D}(t) + \int_0^t ds\,  R_{n}^{m-1,\gamma,D}(s) \right).
$$
By Gronwall's lemma, we obtain 
\begin{equation} \label{e3.20}
\sup_{n,m\ge 1} R_{n}^{m,\gamma,D} (t) \leq C\left(1+ \tilde I_a^{\gamma,D}(t)\right)
\end{equation}
and 
\begin{equation*} 
\sup_{n,m\ge 1} \sup_{t\in[0,T]} R_{n}^{m,\gamma,D} (t) \leq C\left(1+\sup_{t\in[0,T]} \tilde I_a^{\gamma,D}(t)\right)
\end{equation*}
   The Picard iterates satisfy
$$
\lim_{n\to\infty}\sup_{t \in [0,T]} E\left( \Vert(u^m_n(t)-u_n(t)) 1_{K^D(t)} \Vert^q_{L^q(\IR^3)}\right) = 0.
$$
Indeed, this can be proved using the same arguments as those of the proof of Theorem \ref{thm3.1}.

   Therefore, (\ref{e3.13}) follows from (\ref{e3.20}) and Fatou's lemma. Similarly, (\ref{e3.14}) follows from
(\ref{e3.13}), Proposition \ref{prop3.3} and Fatou's lemma. 

To establish (\ref{e3.14.bis}), we apply the same arguments as for the proof of (\ref{e3.7.bis}) in Theorem \ref{thm3.1} to check that
$$
\sup_{t\in[0,T]} E\left(\Vert u(t)\Vert^q_{\gamma, q,K^D(t)}\right)\le C\left( %\sup_{z\in D^T}|\nabla v_0(z)|^q 
\Vert v_0 \Vert_{\gamma,q,D^T}+ \Vert\Delta v_0\Vert^q_{\mathcal{C}^{\gamma_1}(D^T)}+
\Vert \tilde v_0\Vert^q_{\mathcal{C}^{\gamma_2}(D^T)}\right).
$$
This fact together with (\ref{e3.7.bis}) finishes the proof.

%%%%%%%%%%%
Finally, (\ref{e3.15}) is a consequence of  (\ref{e3.14.bis}),
since  $q$ in (\ref{e3.14.bis})  can be arbitrarily large.
 \hfill $\Box$

\medskip

   An important consequence of Theorem \ref{thm3.6} and the Sobolev embeddings is the following.

\begin{cor}
\label{spacecontinuity} Suppose that the hypotheses of Theorem \ref{thm3.6} are satisfied. Then for any fixed $t > 0$, a.s., the sample
paths of $\left(u(t, x)1_{K^D(t)}(x),\ x \in \IR^3\right)$ are $\rho$-H\"{o}lder continuous
 with $\rho \in\, ]0, \gamma_1 \wedge\gamma_2 \wedge \tau(\beta,\delta)[.$

We also notice that for any $q \in [2, \infty[,$
\begin{equation}\label{e3.25}
\sup_{t \in [0,T]}\sup_{x,y \in K^D(t)}E\left(\left\vert u(t,x)-u(t,y)\right\vert^q\right) \leq C \vert x-y \vert^{\rho q},
\end{equation}
for any $\rho \in\, ]0, \gamma_1 \wedge \gamma_2 \wedge \tau(\beta,\delta)[$, and
\begin{equation}\label{e3.26}
   \sup_{t \in[0,T]}\, \sup_{x \in K^D(t)} E\left(\vert u(t,x)\vert^q\right) < \infty.
\end{equation}
\end{cor}

\begin{remark} 
\label{heatcomparison}Notice that H\"older regularity in the space variable is the same as that for the solution
of the stochastic heat equation in any dimension $d\geq 1$ (see \cite{sssarra}, Theorem 2.1).
%cite  Anja Sturm's talk in Vancouver?
\end{remark}

%%%%%% REGULARITY IN TIME

\subsection{Regularity in time}
\label{ss2.3}
   Our next aim is to analyze the behavior in time of the solution of (\ref{e3.2'}). We begin by studying the properties of the
term corresponding to the contribution of the initial conditions.

Throughout the section $D$ is an arbitrary bounded domain of $\IR^3$.

\begin{lemma} Let $v_0$, $\tilde v_0$ be real-valued functions.
\begin{description}
\item{(a)} Let $v_0 \in C^2(\IR^3)$ be such that $\Delta v_0$ is $\gamma_1$-H\"{o}lder continuous for some $\gamma_1 \in\, ]0,1]$. 
Then, the mapping $t \mapsto (\frac{d}{dt} G(t) \ast v_0)(x)$ ($t\leq T$) is also $\gamma_1$-H\"{o}lder continuous  
and
\begin{equation}
\label{refc1}
\sup_{x\in D}\left\Vert \frac{d}{dt} G(\cdot) \ast v_0)(x)\right\Vert_{\mathcal{C}^{\gamma_1}([0,T])}\le C \left(\sup_{y\in D^T}|\nabla v_0(y)| + \Vert \Delta v_0\Vert_{\mathcal{C}^{\gamma_1}([0,T])}\right).
\end{equation}
Therefore, its H\"older semi-norm is uniformly bounded in $x\in D$.
\item{(b)} Assume that $\tilde v_0$ is  $\gamma_2$-H\"{o}lder continuous, for some $\gamma_2 \in\, ]0,1].$ Then the mapping $t \mapsto (G(t) \ast \tilde v_0)(x)$ ($t\leq T$) is also  $\gamma_2$-H\"{o}lder continuous and 
\begin{equation}
\label{refc2}
\sup_{x\in D}\left\Vert (G(\cdot) \ast \tilde v_0)(x)\right\Vert_{\mathcal{C}^{\gamma_2}([0,T])}\le C \left(\sup_{x\in D^T}|\tilde v_0(x)|+\Vert \tilde v_0\Vert_{\mathcal{C}^{\gamma_2}(D^T)}\right).
\end{equation}
Hence, its H\"older semi-norm is uniformly bounded in $x\in D$.
\end{description}
\label{lem3.9}
\end{lemma}

\noindent{\it Proof}: We check (a) by studying each of the terms on the right-hand side of (\ref{e3.8}). Fix $t, \tbar \in [0,T]$ and $x \in D$. Then
\begin{align*}
I_1(t, \tbar,x)&:= \left\vert\frac{1}{t}(G(t) \ast v_0)(x) - \frac{1}{\tbar} (G(\tbar) \ast v_0)(x)\right\vert\\
\\
%&=& \vert \frac{1}{t} \int_{\IR^3} G(t, du) v_0(x-u) - \frac{1}{t} \int_{\IR^3} %G(\tbar,dv)v_0(x-v)\vert\\
&= \frac{1}{t} \left\vert \int_{\IR^3}G(t,du) \left(v_0(x-u)-v_0\left(x- \frac{\tbar}{t} u\right)\right)\right\vert,
\end{align*}
where we have used (\ref{A.2a}). The mean value theorem yields
\begin{align}
 I_1(t, \tbar, x) \vert &\leq \frac{1}{t} \int_{\IR^3}G(t, du) \sup_{y \in D^T} \left(\vert \nabla v_0(y)\vert\right)\, \vert u
\vert\, \frac{\vert t- \tbar \vert}{t}\nonumber
\\
&\leq C\sup_{y\in D^T}\vert \nabla v_0(y)\vert\, \vert t-\tbar \vert,
\label{e3.27}
\end{align}
since $G(t,\cdot)$ is concentrated on $B_t(0)$ and has total mass $t$. By the H\"{o}lder-continuity property of $\Delta v_0$,
\begin{align}
I_2(t, \tbar, x) &:= \left\vert \int_{\vert y \vert < 1} (\Delta v_0(x+ty)-\Delta v_0(x+ \tbar y))dy \right\vert \nonumber
\\
&\leq  \Vert \Delta v_0\Vert_{\mathcal{C}^{\gamma_1}(D^T)} \vert t-\tbar \vert^{\gamma_1},
\label{e3.28}
\end{align}
%for any $\alpha \in\, ]0, \gamma_1[.$
The estimates  (\ref{e3.27}), (\ref{e3.28}) and the identity (\ref{e3.8}) yields the result stated in (a).

%$$
%\sup_{x \in K^D(T)} \left\vert \left(\frac{d}{dt} G(t) \ast v_0\right) (x) - \left(\frac{d}{dt} G(\tbar) \ast v_0\right) (x)\right \vert \leq C \vert t-\tbar
%\vert^{\gamma_1},
%$$
%for any $\alpha \in\, ]0, \gamma_1[.$
 %This yields the result stated in (a).

   Let us now prove (b). For any $t, \tbar \in [0,T]$,
\begin{align*}
   I_3(t, \tbar, x) &:= \vert (G(t) \ast \tilde v_0)(x) - (G(\tbar) \ast \tilde v_0)(x) \vert\\
    \\
%    &=& \vert \int_{\IR^3}(G(t, du) - G(\tbar, du)) \tilde v_0(x-u)\vert \\
&= \left\vert \int_{\IR^3} G(t,du)\left(\tilde v_0(x-u)-\tilde v_0\left (x-\frac{\tbar}{t} u\right)\, \frac{\tbar}{t}\, \right)\right\vert,
\end{align*}
where, in the last equality, we have applied  (\ref{A.2a}).

   Hence
$$
  I_3(t, \tbar,x) \leq I^1_3(t, \tbar,x) + I^2_3(t, \tbar, x),
$$
where
\begin{align*}
   I^1_3(t, \tbar, x) &= \int_{\IR^3}G(t, du) \,\frac{\tbar}{t}\, \left\vert \tilde v_0(x-u)-\tilde v_0\left(x-u \frac{\tbar}{t}\right)
\right\vert, \\
\\
I^2_3 (t, \tbar, x) &= \int_{\IR^3} G(t, du) \,\vert \tilde v_0(x-u) \vert\, \left\vert 1- \frac{\tbar}{t} \right\vert.
\end{align*}
Using the H\"{o}lder-continuity property of $\tilde v_0$, we obtain 
$$\sup_{x \in D} I^1_3(t, \tbar, x) \leq C \Vert \tilde v_0\Vert_{\mathcal{C}^{\gamma_2}(D^T)}\vert
t-\tbar \vert^{\gamma_2}.$$
% with $\alpha \in\, ]0, \gamma_2[.$
Moreover, since $\tilde v_0$ is continuous and has compact
support, 
$$\sup_{x \in D} I^2_3(t, \tbar, x) \leq C \sup_{x\in D^T} |\tilde v_0(x)|\vert t-\tbar \vert,$$
 and statement (b) is established.
\qed
\bigskip

The next theorem provides upper bounds on $L^q$--moments of time increments of the solution of
equation (\ref{e3.2'}), uniformly over $x$ in bounded sets. We recall the notation $\tau(\beta,\delta)=\frac{2-\beta}{2}\wedge\frac{1+\delta}{2}$.

\begin{thm} 
\label{tcontinuity}Assume that
\begin{description}
\item{(a)} $\sigma, b$ are Lipschitz functions,
\item{(b)} $v_0$, $\tilde v_0$ are real-valued functions, $v_0 \in C^2(\IR^3)$, and $\Delta v_0$ and $\tilde v_0$ are H\"{o}lder-continuous functions of order $\gamma_1, \gamma_2 \in\, ]0,1]$, respectively.
\end{description}
   Then, for any $t, \tbar \in [0,T]$ and each $q \in [2, \infty[,$
\begin{equation}\label{e3.29}
\sup_{x \in D} E\left(\vert u(t,x) - u(\tbar,x) \vert^q\right) \leq C \vert t-\tbar\vert^{\alpha q},
\end{equation}
where $C$ is a positive constant and $\alpha \in\, ]0, \gamma_1\wedge \gamma_2 \wedge\tau(\beta,\delta)[$.
\end{thm}

\noindent{\it Proof}:  Fix $x \in D$ and $q \in [2, \infty[.$ Then $x \in K^D(t)$, for all $t \in [0,T]$, so by (\ref{e3.2'}) with $a=1$,
$$
    E\left(\vert u(t,x)-u(\tbar,x)\vert^q\right) \leq C \sum^q_{i=1} T_i(x, t, \tbar),
$$
where
\begin{align*}
T_1(x,t,\tbar) &= \left\vert \left(\frac{d}{dt} G(t) \ast v_0\right)(x) - \left(\frac{d}{dt}G(\tbar) \ast v_0\right)(x)\right\vert^q,\\
\\
T_2(x,t,\tbar) &= \left\vert (G(t) \ast \tilde v_0) (x) - (G(\tbar) \ast \tilde v_0)(x)\right\vert^q,\\
\\
T_3(\cdot,t,\tbar) &= E\Big(\big\vert \int_0^t \int_{\IR^3} G(t-s, \cdot-y) \sigma(u(s,y)) 1_{K^D(s)}(y) M(ds,dy)\\
\\
&\qquad- \int^{\tbar}_0 \int_{\IR^3} G(\tbar-s, \cdot-y) \sigma(u(s,y)) 1_{K^D(s)}(y) M(ds,dy)\big\vert^q\Big)\\
\\
T_4(x,t,\tbar) &= E\Big(\big\vert \int_0^t ds \int_{\IR^3} G(t-s, dy)\, b(u(s, x-y))1_{K^D(s)}(x-y)\\
\\
&\qquad- \int_0^{\tbar} ds \int_{\IR^3} G(\tbar-s, dy)\, b(u(s,x-y))1_{K^D(s)}(x-y)\big\vert^q\Big)
\end{align*}
Lemma \ref{lem3.9} yields
\begin{equation}\label{e3.30}
   \sup_{x \in K^D(T)} \left(T_1(x,t,\tbar) + T_2(x,t,\tbar)\right) \leq C \vert t-\tbar \vert^{q\alpha_1},
\end{equation}
with $\alpha_1 \in\, ]0, \gamma_1 \wedge \gamma_2[.$

%   Set 
%$$
%   Z(s,y) := \sigma(u(s,y))1_{K^D(s)}(y),\qquad (s,y) \in [0,T] \times \IR^3
%$$ 
Let $g$ be a real-valued Lipschitz function and set $\mathcal{I}:=\, ]0, \gamma_1 \wedge \gamma_2 \wedge \tau(\beta,\delta)[$. We next show that
\begin{description}
\item{(1)} $\sup_{t\in[0,T]} E\left(\Vert g(u(t))\Vert^q_{L^q(K^D(t))}\right) <\infty$;
\item{(2)} for any $\gamma\in \mathcal{I}$, $t \in \, ]0,T]$,
\begin{equation*}
   \sup_{s\in [0,t]} E\left(\Vert g(u(s))1_{K^D(s)}\Vert^q_{\gamma,q,(K^D(t))^{t-s}}\right) <\infty.
\end{equation*}
\end{description}
Indeed, (1) follows from the linear growth property of $g$ and Theorem \ref{thm3.1}.
To prove (2), we notice that since $(K^D(t))^{t-s}\subset K^D(s)$, the arguments in the proof of property (ii) in Theorem \ref{thm3.6}
give
\begin{align*}
&\sup_{s\in[0,t]} E\left(\Vert g(u(s))1_{K^D(s)}\Vert^q_{\gamma,q,(K^D(t))^{t-s}}\right)\\
&\quad\le C \sup_{s\in[0,t]} E\left(\Vert u(s) \Vert^q_{\gamma,q,K^D(s)}\right),
\end{align*}
and the conclusion follows from (\ref{e3.15}).

By properties (1) and (2) above, we see that the stochastic process $Z(s,y) := g(u(s,y))1_{K^D(s)}(y),$ $(s,y) \in [0,T] \times \IR^3$, 
satisfies Assumption \ref{assumpI} with $\cO:=D$ and arbitrarily large $q$.
By applying Theorem \ref{thm2.3}, we conclude that for any $q\in[1,\infty[$,
\begin{equation}\label{e3.31}
   \sup_{x \in K^D(T)} T_3(x,t,\tbar) \leq C \vert t-\tbar\vert^{q\alpha_2 },
\end{equation}
with $\alpha_2 \in \mathcal{I}$.

   Finally, we study the contribution of $T_4 (x,t,\tbar)$. Assume $t \leq \tbar$ and consider the decomposition
$$
   T_4(x,t,\tbar) = C\left(T_{4,1}(x,t,\tbar) + T_{4,2}(x,t,\tbar)\right),
$$
with
\begin{align*}
T_{4,1}(x,t,\tbar) &= E\Big(\big\vert \int^{\tbar}_t ds \int_{\IR^3}G(\tbar-s, dy)\, b(u(s,x-y))1_{K^D(s)}(x-y) \big\vert^q\Big),\\
T_{4,2}(x,t,\tbar) &= E\Big(\big\vert \int_0^t ds \int_{\IR^3}\big(G(t-s,dy)-G(\tbar-s,dy)\big)b(u(s,x-y))\\
&\qquad\times 1_{K^D(s)}(x-y)\big\vert^q\Big).
\end{align*}
H\"{o}lder's inequality, the restriction on the growth of $b$ and (\ref{e3.26}) imply that
\begin{align}
\sup_{x \in D} T_{4,1}(x,t,\tbar) &\leq C \vert t-\tbar \vert^{q-1}\nonumber\\
&\quad\times \int^{\tbar}_t ds \int_{\IR^3}G(\tbar-s, dy)\,
\left(1 + \sup_{t \in [0, T]} \sup_{y \in K^D(t)} E(\vert u(t,y)\vert^q)\right) \nonumber\\
&\leq C \vert t-\tbar \vert^q.
\label{e3.32}
\end{align}

   Notice that for $x \in D$ and $y \in B_{t-s}(0)$, $x - y \in K^D(s)$, so $1_{K^D(s)}(x-y)$ can be removed from the
expression for $T_{4,2}(x,t,\tbar)$ when $x \in D$. We split the integral in the definition of $T_{4,2}(x,t,\tbar)$ into a difference of integrals and then we apply the transformations $y \mapsto \frac{y}{t-s}$ and $y \mapsto
\frac{y}{\tbar-s}$, respectively. We obtain
\begin{align*}
T_{4,2}(x,t,\tbar) &= E\Big(\big\vert \int_0^t ds \int_{B_1(0)} G(1, dy)\, b(u(s,x-(t-s)y))\, (t-s)\\
 & \quad - \int_0^t ds \int_{B_1(0)} G(1,dy)\, b(u(s,x-(\tbar-s)y))\, (\tbar-s)\big\vert^q\Big).
\end{align*}
Hence $T_{4,2}(x,t,\tbar) \leq C (T_{4,2,1}(x,t,\tbar) + T_{4,2,2}(x,t,\tbar))$, where
\begin{align*}
T_{4,2,1}(x,t,\tbar) &= \vert t-\tbar\vert^q E\Big(\big\vert \int_0^t ds \int_{B_1(0)}G(1,dy)\,
b(u(s,x-(\tbar-s)y))\big\vert^q\Big),\\
T_{4,2,2}(x,t,\tbar) &= E\Big(\big\vert \int_0^t ds\, (t-s) \int_{B_1(0)}G(1,dy)\, \big(b(u(s,x-(t-s)y))\\
&\quad  - b(u(s,x-(\tbar-s)y))\big)\big\vert^q\Big).
\end{align*}
Clearly, from (\ref{e3.26}) and the linear growth property of $b$, it follows that 
\begin{equation}\label{e3.33}
\sup_{x \in K^D(t)} T_{4,2,1}(x,t,\tbar) \leq C \vert t-\tbar \vert^q.
\end{equation}
Moreover, by (\ref{e3.25}) and the Lipschitz property of $b$,
\begin{align}
\sup_{x \in K^D(t)} T_{4,2,2}(x,t,\tbar) &\leq C \vert t-\tbar \vert^{\rho q} \int_0^t ds \int_{B_1(0)} G(1,dy)\, \vert y \vert^{\rho q}\nonumber\\ 
&\leq C \vert t-\tbar \vert^{\rho q},
\label{e3.34}
\end{align}
for any $q\in[1,\infty[$, $\rho \in\,\mathcal{I}$.

   Putting together (\ref{e3.33}) and (\ref{e3.34}), we obtain
\begin{equation}
   \sup_{x \in K^D(t)} T_{4,2} (x, t, \tbar) \leq C \vert t-\tbar \vert^{\rho q},
\end{equation}
for any $q\in[1,\infty[$, $\rho \in \mathcal{I}.$
Together with (\ref{e3.32}), this yields
\begin{equation}\label{e3.36}
   \sup_{x \in K^D(t)} T_4(x,t,\tbar) \leq C \vert t-\tbar \vert^{\rho q},
\end{equation}
for each $q\in[1,\infty[$, $\rho \in \mathcal{I}.$
\smallskip

With the estimates (\ref{e3.30}), (\ref{e3.31}) and (\ref{e3.36}), we obtain (\ref{e3.29}).
\hfill$\Box$
\vskip 16pt

We summarize the results of this section (Corollary \ref{spacecontinuity} and Theorem \ref{tcontinuity}) as follows.
\begin{thm}
\label{holdercont} Assume that:
\begin{description}
\item{(a)} The covariance of $ F$ is of the form given in Assumption \ref{assumpC};
\item{(b)} the initial values $v_0$, $\tilde{v}_0$ are such that $v_0 \in C^2(\IR^3)$, and $\Delta v_0$ and $\tilde v_0$ are H\"{o}lder continuous with orders $\gamma_1,\gamma_2 \in\, ]0,1]$, respectively;
\item{(c)} the coefficients $\sigma$ and $b$ are Lipschitz.
\end{description}
Then for any $q\in[2,\infty[$, $\alpha \in\, ]0, \gamma_1 \wedge \gamma_2 \wedge \tau(\beta,\delta)[$,
there is $C >0$ such that for $(t,x), (\tbar,y) \in [0,T]\times D$, 
\begin{equation}\label{moment_incr}
   E\left(\vert u(t,x) - u(\tbar,y)\vert^q\right) \leq C\left(\vert t - \tbar\vert + \vert x - y \vert\right)^{\alpha q}.
\end{equation}
In particular, a.s., the stochastic process $\left(u(t,x),\ (t,x) \in[0,T] \times D\right)$ solution of (\ref{e3.2'}) has
$\alpha$-H\"{o}lder continuous sample paths, jointly in $(t,x)$. 
\end{thm}

\begin{remark} 
\label{detwave}
Assume that the initial values $v_0$, $\tilde{v}_0$ satisfy the conditions given in (b) above. With Remark \ref{rem3.2} and Lemma \ref{lem3.9} we conclude that the solution of the deterministic inhomogeneous wave equation (take $\sigma \equiv 0$ in (\ref{e3.2'})) in dimension $d=3$ is $\alpha$--H\"older continuous, jointly in $(t,x)$, for any $\alpha\in\, ]0,\gamma_1\wedge\gamma_2[$.
\end{remark}

%%%%%%%%%%%%%%%%BEGINNING SECTION 3

\section{Sharpness of the results}
\label{s3}

We devote this section to showing that the results of  Theorem \ref{holdercont}
on H\"older continuity in space and time are optimal. To do this, we consider the special case of equation (\ref{e3.2'}) with vanishing initial conditions $v_0$, $\tilde v_0$,
coefficients $\sigma\equiv 1$, $b\equiv 0$, and covariance function of the noise 
given by $f(x)=k_\beta(x)$, with $\beta\in\, ]0,2[$. In this case, the solution of equation (\ref{e3.2'}) defines a stationary Gaussian process.

\begin{thm}
\label{sharp} 
Under the above assumptions on $v_0$, $\tilde v_0$, $\sigma$ and $b$, the solution 
$u(t,x)$ to the stochastic partial differential equation (\ref{e3.2'}) has the following properties:
\begin{description}
\item{(a)} Fix $t\in\,]0,T]$ and a compact set $K \subset \IR^3$. There is a constant $c_1>0$ such that for all $x, y \in K$,
\begin{equation}
\label{sharpx}
E\left(\left\vert u(t,x) - u(t,y)\right\vert^2\right)\ge c_1\vert x-y\vert^{2-\beta}.
\end{equation}
Consequently, a.s. the mapping $x\mapsto u(t,x)$ is {\em not} $\gamma$--H\"older continuous for $\gamma >\frac{2-\beta}{2}$, though it is for
$\gamma<\frac{2-\beta}{2}$.
\item{(b)} Fix $x\in\IR^3$ and $t_0>0$. There is a constant $c_2>0$ such that for all $t, \tbar\in [t_0,T]$ with $|\tbar-t|$ sufficiently small,
\begin{equation}
\label{sharpt}
E\left(\left\vert u(t,x) - u(\tbar-x)\right\vert^2\right)\ge c_2\vert t-\tbar\vert^{2-\beta}.
\end{equation}
Hence, a.s., the mapping $t\mapsto u(t,x)$ is {\em not} $\gamma$--H\"older continuous for $\gamma >\frac{2-\beta}{2}$, while it is for 
$\gamma<\frac{2-\beta}{2}$.
\end{description}
\end{thm}

\noindent{\it Proof.}\quad For any $x\in\IR^3$, set $R(x) = E\left(u(t,x) u(t,0)\right)$. Since
\begin{equation*}
E\left(\left\vert u(t,x) - u(t,y)\right\vert^2\right) = 2\left(R(0)-R(x-y)\right),
\end{equation*}
it suffices to show that for any $x\in\IR^3$,
\begin{equation*}
R(0)-R(x)\ge C_1\vert x\vert^{2-\beta},
 \end{equation*}
 for some positive constant $C_1$. Without loss of generality, we may assume that $t=1$. 

We remark that $R(0)-R(x)$ is a real number. Taking (\ref{fourierG}) into account and integrating with respect to $s$  yields
\begin{align}
 R(0)-R(x)& = \int_0^1 ds \int_{\IR^3}\,\mu(d\xi)(1-e^{i\xi\cdot x})\,\left\vert\mathcal{F}G(1-s)(\xi)\right\vert^2 \nonumber\\
 &= \frac{1}{2}\int_{\IR^3}\,\frac{d\xi}{\vert\xi\vert^{5-\beta}}\left(1-\cos(\xi\cdot x)\right)\left(1-\frac{\sin(2|\xi|)}{2|\xi|}\right),
 \label{cov}
\end{align}
where $\xi\cdot x$ denotes the Euclidean inner product of the vectors $\xi$ and $x$.

If $|\xi|>1$, then $1-\frac{\sin(2|\xi|)}{2|\xi|}\ge\frac{1}{2}$. Thus, using the change of variables
$w=|x|\xi$ and setting $e=\frac{x}{|x|}$ we obtain
 \begin{equation*}
  R(0)-R(x)\ge \frac{1}{4} \,|x|^{2-\beta}\,\int_{|w|>|x|} \frac{dw}{|w|^{5-\beta}}\left(1-\cos(w\cdot e)\right).
  \end{equation*}
Because $x \in K$ and $K$ is compact, the last integral is bounded below by a positive constant, hence (\ref {sharpx}) is proved.

  We now prove (\ref{sharpt}), assuming that $t_0 \leq t < \tbar \le T$. In this situation,
\begin{equation*}
  E\left(\left\vert u(t,x) - u(\tbar,x)\right\vert^2\right) = T_1(t,\tbar;x) + T_2(t,\tbar;x),
\end{equation*}
where
\begin{align*}
  T_1(t,\tbar;x)&= E\left(\int_t^{\tbar}\int_{\IR^3} G(\tbar-s,x-y)\,M(ds,dy)\right)^2,\\
  T_2(t,\tbar;x)&= E\left(\int_0^t\left(G(\tbar-s,x-y )-G(t-s,x-y)\right)\,M(ds,dy)\right)^2.
\end{align*}
Using again the explicit formula (\ref{fourierG}) and integrating with respect to the variable $s$ yields
\begin{align*}
  T_1(t,\tbar;x)&=\int_t^{\tbar} ds\,\int_{\IR^3}\frac{d\xi}{|\xi|^{3-\beta}}\,\left\vert\mathcal{F}G(\tbar-s)(\xi)\right\vert^2\\
  \\
  &= \frac{1}{2}\int_{\IR^3}\,\frac{d\xi}{|\xi|^{5-\beta}}\left((\tbar-t)-\frac{\sin(2(\tbar-t)|\xi|)}{2|\xi|}\right).
\end{align*}
With the change of variables $w=(\tbar-t)\xi$, the last integral becomes 
\begin{equation*}
\vert\tbar-t\vert^{3-\beta}\,\int_{\IR^3}\,\frac{dw}{|w|^{5-\beta}}\left(1-\frac{\sin(2|w|)}{2|w|}\right),
\end{equation*}
which clearly yields
\begin{equation}
  \label{sharpt1}
  T_1(t,\tbar;x)\ge C \vert\tbar-t\vert^{3-\beta}.
\end{equation}
A direct integration in the  variable $s$ yields
\begin{align*}
  T_2(t,\tbar;x)&=\int_0^t ds\,\int_{\IR^3}\frac{d\xi}{|\xi|^{3-\beta}}\left\vert\left(\mathcal{F}G(\tbar-s,\cdot)-\mathcal{F}G(t-s,\cdot)\right)(\xi)\right\vert^2\\
  \\
  &\ge \int_0^t ds\,\int_{|\xi|\ge(\tbar-t)^{-1}}\frac{d\xi}{|\xi|^{5-\beta}}\left(\sin\left((\tbar-s)|\xi|\right)-\sin\left((t-s)|\xi|\right)\right)^2\\
  \\
  &=\int_{|\xi|\ge(\tbar-t)^{-1}}\frac{d\xi}{|\xi|^{5-\beta}}\left(A(t,\tbar,\xi)+B(t,\tbar,\xi)\right),
\end{align*}
where
\begin{align*}
 A(t,\tbar,\xi)&=t\left(1-\cos\left((\tbar-t)|\xi|\right)\right)+\frac{\sin\left((t+\tbar)|\xi|\right)}{2|\xi|}\,\left(1-\cos\left((\tbar-t)|\xi|\right)\right),\\
   \\
  B(t,\tbar,\xi)&=\frac{\sin\left(2(\tbar-t)|\xi|\right)}{4|\xi|} - \frac{\sin\left((t-\tbar)|\xi|\right)}{2|\xi|}.
\end{align*}
Changing the variable $\xi$ into $w=(\tbar-t)\xi$ and bounding below the second term in $A(t,\tbar,\xi)$ yields
  %I DO NOT UNDERSTAND THE FIRST LOWER BOUND
\begin{align*}
  \int_{|\xi|\ge(\tbar-t)^{-1}}\frac{d\xi}{|\xi|^{5-\beta}}A(t,\tbar,\xi)&\ge %\int_{|\xi|\ge(\tbar-t)^{-1}}\frac{d\xi}{|\xi|^{5-\beta}}\\
%  &\quad \times
%  \left(t\left(1-\cos\left((\tbar-t)|\xi|\right)\right)-(\tbar-t)\right)\\
%  \\
  |\tbar-t|^{2-\beta} \int_{|w|>1}\,\frac{dw}{|w|^{5-\beta}}\left(t\left(1-\cos |w|\right)-(\tbar-t)\right)\\
  \\
  &\geq k_1|\tbar-t|^{2-\beta}-k_2|\tbar-t|^{3-\beta}%\ge \frac{k_1}{2}\,|\tbar-t|^{2-\beta},
\end{align*}
(the last inequality uses that fact that $t \geq t_0 >0$). 

  Similarly,
\begin{align*}
   \int_{|\xi|\ge(\tbar-t)^{-1}}\frac{d\xi}{|\xi|^{5-\beta}}B(t,\tbar,\xi)&\geq \int_{|\xi|\ge(\tbar-t)^{-1}}\frac{d\xi}{|\xi|^{5-\beta}} \left(-\frac{1}{4|\xi|} - \frac{1}{2|\xi|} \right)\\
   &\geq -\frac{3}{4} |\tbar-t|^{3-\beta} \int_{|w|>1} \frac{dw}{|w|^{5-\beta}}.
\end{align*}
Therefore,
$$
   T_2(t,\tbar;x) \geq k_1 |\tbar-t|^{2-\beta} - k_2' |\tbar-t|^{3-\beta} 
      \geq \frac{k_1}{2} |\tbar-t|^{2-\beta},
$$
where the last lower bound holds whenever $|\tbar-t|\le \frac{k_1}{2k_2'}$. This finishes the proof of (\ref{sharpt}).
 \smallskip %&=|\tbar-t|^{3-\beta}\int_{|w|>1}\frac{dw}{|w|^{5-\beta}}\left(\frac{\sin(2|w|)}{4|w|}-\frac{\sin(|w|)}{2|w|}\right)\\
%   \\
%   &=-k_3|\tbar-t|^{3-\beta},  \end{align*}
%  for some $k_3>0$. Thus,
%   $$T_2(t,\tbar;x)\ge k_1|\tbar-t|^{2-\beta}-(k_2+k_3)|\tbar-t|^{3-\beta}\ge \frac{k_1}{2}\,|\tbar-t|^{2-\beta},$$
%where the last lower bound holds whenever $|\tbar-t|\le \frac{k_1}{2(k_2+k_3)}$.

   The statements concerning absence of H\"older continuity follow from a well-known result on Gaussian stationary process (see for instance \cite{adler}, Theorem 3.2). We only give
 some details for the space variable, since the arguments for the time variable are the same.

 Assume by contradiction that for a fixed $t\in\, ]0,T]$, the sample paths $x\mapsto u(t,x)$ are $\gamma$--H\"older continuous for some $\gamma>\frac{2-\beta}{2}$, so that
 for any compact set $K\subset\IR^3$, there is $0<C(\omega)<\infty$ with
 \begin{equation*}
 \sup_{x,y\in K,x\ne y}\frac{u(t,x)-u(t,y)}{|x-y|^{\gamma}} < C(\omega).
 \end{equation*}
 This implies that the Gaussian stochastic process 
 \begin{equation*}
 \left(\frac{u(t,x)-u(t,y)}{|x-y|^{\gamma}},\ x,y\in K,\ x\ne y\right)
 \end{equation*}
 is finite a.s., and even, by Theorem 3.2 of \cite{adler}, that
 \begin{equation*}
 E\left( \sup_{x,y\in K,x\ne y}\left( \frac{u(t,x)-u(t,y)}{|x-y|^{\gamma}}\right)^2\right) < \infty.
 \end{equation*}
 In particular, there would exist a finite $K>0$ such that
 \begin{equation*}
  E\left( \vert u(t,x) - u(t,y)\vert^2\right) \le K|x-y|^{2\gamma}\,.
  \end{equation*}
  This clearly contradicts (\ref{sharpx}).
  \hfill $\Box$
  \medskip

\begin{remark}
  \label{optimal}
Under the same assumptions as in Theorem \ref{sharp}, for any fixed $t\in\, ]0,T]$, a.s.  the upper bounds on moments of increments of the process
$\left(u(t,x),\ x\in\IR^3\right)$ remains valid even for $\alpha = \frac{2-\beta}{2}$ (this, however, does not yield an improvement in the H\"older continuity of sample paths of this process). 

Indeed, by (\ref{cov}), we can consider the decomposition
\begin{equation*}
R(0)-R(x)=r_1(x)+r_2(x)+r_3(x),
\end{equation*}
with
\begin{equation*}
r_1(x) = \frac{1}{2}\int_{|\xi|\le 1}\,\frac{d\xi}{\vert\xi\vert^{5-\beta}}\left(1-\cos(\xi\cdot x)\right)\left(1-\frac{\sin(2|\xi|)}{2|\xi|}\right),%\\
%\\
%r_2(x) &= \frac{1}{2}\int_{|\xi||x|\ge 1}\,\frac{d\xi}{\vert\xi\vert^{5-\beta}}\left(1-cos(\xi\cdot x)\right)\left(1-\frac{\sin(2|\xi|)}{2|\xi|}\right)\\
%\\
%r_3(x) &= \frac{1}{2}\int_{1\le|\xi|\le|x|^{-1}}\,\frac{d\xi}{\vert\xi\vert^{5-\beta}}\left(1-cos(\xi\cdot x)\right)\left(1-\frac{\sin(2|\xi|)}{2|\xi|}\right)\\
\end{equation*}
and $r_2(x)$ (resp.~$r_3(x)$) defined by the same expression, but with the domain of integration replaced by $|\xi||x| \geq 1$ (resp.~$1\le|\xi|\le|x|^{-1}$).
Since $1-\cos(\xi\cdot x)\le \frac{1}{2}|\xi|^2 |x|^2$, we clearly have $r_1(x)\le C|x|^2$. Moreover, bounding above the products of trigonometric functions, we obtain
\begin{equation*}
r_2(x) \le C \int_{|\xi||x|\ge 1} \frac{d\xi}{|\xi|^{5-\beta}} = C |x|^{2-\beta}.
\end{equation*}
Finally, for $|x|$ small,
\begin{align*}
r_3(x) &\le C\,|x|^2 \int_{1\le|\xi|\le|x|^{-1}} \frac{d\xi}{|\xi|^{3-\beta}} = C\,|x|^2\left(|x|^{-\beta}-1\right)\\
&\le C\,|x|^{2-\beta}.
\end{align*}
Consequently,
\begin{equation*}
R(0) - R(x) \le C |x|^{2-\beta}.
\end{equation*}
\end{remark}

%%%%%%%%%%%%%%%%%%%%%%BEGINNING OF SECTION 6
\section{Integrated increments of the covariance function}
\label{B}

In this section, we prove technical results on integrals involving regularisations of the fundamental solution $G$ and
one and two-dimensional increments of $f$. Each of these results plays a crucial role in the proofs of section \ref{s1}.  

%We consider next increments of the covariance density $f$.

\begin{lemma} For any $\alpha \in\, ]0, (2-\beta) \wedge 1[$, there is $C < \infty$ such that for all $x,y \in \IR^3$,
\begin{equation}\label{B.2}
   \sup_{n\geq 1} \sup_{s \in [0,T]} I_n(s,x,y) \leq C \vert x-y\vert^\alpha,
\end{equation}
where
$$
   I_n (s, x, y) = \int_{\IR^3} du \int_{\IR^3} dv\, G_n(s,u)G_n(s,v) \vert Df(v-u,x-y) \vert . 
$$
\label{lemB.2}
\end{lemma}

\noindent{\it Proof.} The structure of the function $f$ suggests the decomposition 
$$I_n(s,x,y) \leq I^1_n(s,x,y) + I^2_n(s,x,y),$$
where $I^1_n(s,x,y)$ (resp. $I^2_n(s,x,y)$) denotes the same expression but with the factor in absolute
values replaced by 
$$
   \varphi(x-y+v-u) \ \vert Dk_\beta(v-u,x-y) \vert,
   \qquad \mbox{(resp.}\quad
   \vert D\varphi(v-u,x-y)\vert  k_\beta (v-u)).
$$
Since the function $\varphi$ is bounded, we can apply Lemma \ref{lemB.1}(a) with $d=3$, $b:=\alpha$, $a:= 3-(\alpha+\beta)$, $u:=v-u$, $c:=|x-y|$, $x:=\frac{x-y}{|x-y|}$, to
obtain
\begin{align*}
I^1_n(s,x,y) &\leq C \vert x-y \vert^\alpha  \int_{\IR^3} du \int_{\IR^3} dv\, G_n(s,u)G_n(s,v) \\
&\quad\times \int_{\IR^3} dw\, k_{\alpha+\beta}(\vert x - y\vert w +u-v) \vert Dk_{3-\alpha}(w,\bar e)
\vert,
\end{align*}
where $\bar e = \frac{x-y}{\vert x-y \vert}$, $\alpha + \beta \in\, ]0,2[$.

   Observe that
\begin{align*}
   &  \int_{\IR^3} du \int_{\IR^3} dv\, G_n(s,u)G_n(s,v) k_{\alpha+\beta}(\vert x - y\vert w +u-v)\\
& \quad =  \int_{\IR^3} d\xi\, e^{i\langle \xi, \vert y-x \vert w\rangle}\, \vert {\cal{F}}
G_n(s)(\xi)\vert^2 k_ {3-(\alpha+\beta)}(\xi)\\
& \quad \leq \int_{\IR^3}d\xi\, \vert {\cal{F}} G(s)(\xi)\vert^2
k_{3-(\alpha+\beta)}(\xi),
\end{align*}
and this is uniformly bounded over $s \in [0,T]$ if $\alpha+\beta \in\, ]0,2[$ (see \ref{A.3}). On the other hand,
$$
   \int_{\IR^3} dw\, \vert Dk_{3-\alpha}(w,e) \vert < \infty,
$$
if $\alpha \in\, ]0,1[$, as it is proved in Lemma \ref{lemB.1}(b). Consequently,
\begin{equation}\label{B.3}
   \sup_{n \geq 1} \sup_{t\in [0,T]} I^1_n (t, x, y) \leq C \vert x-y \vert^\alpha,
\end{equation}
with $\alpha \in\, ]0,(2-\beta) \wedge 1[$.

   The mean value theorem,  the properties of $\varphi$ and (\ref{A.3}) imply that
\begin{align}
 I^2_n(t,x,y) &\leq C \vert x-y \vert \displaystyle \int_{\IR^3} du \int_{\IR^3} dv\, G_n (s,u) G_n(s,v) k_\beta(v-u) \nonumber\\
 &\leq C \vert x-y\vert \int_{\IR^3} \vert {\cal{F}} G(s) (\xi)\vert^2 k_{3-\beta} (\xi) d\xi\nonumber \\ 
 &\leq C \vert x-y\vert.\label{B.4}
\end{align}
The estimates (\ref{B.3}) and (\ref{B.4}) establish (\ref{B.2}).
\hfill $\Box$
\vskip 16pt

\begin{lemma} For each $\alpha \in\, ]0, (2-\beta) \wedge (1+ \delta)[$, there is $C < \infty$ such that for all $x,y\in\IR^3$,
\begin{equation}\label{B.5}
   \sup_{n \geq 1} \sup_{s \in [0,T]} J_n(s,x,y) \leq C \vert x-y \vert^\alpha,
\end{equation}
where
\begin{equation}\label{3.36b}
   J_n(s,x,y) = \displaystyle \int_{\IR^3} du \int_{\IR^3} dv\, G_n(s,u)G_n(s,v)\, \vert D^2 f(v-u,y-x) \vert.
\end{equation}
\label{lemB3}
\end{lemma}

\noindent {\it Proof.} Due to the structure of the covariance function $f$, it is easy to check that
$$
   D^2 f(v-u,y-x) = \sum^3_{1=1} \Delta^{2,i}f(v-u,y-x),
$$
where

\begin{align*}
 \Delta^{2,1}f(v-u,y-x)&= \varphi(v-u) D^2 k_\beta(v-u,y-x), \\
 \Delta^{2,2}f(v-u,y-x) &= D\varphi(v-u,y-x)\,(Dk_\beta(v-u,y-x)\\
 &\quad  - Dk_\beta(v-u,x-y)),\\
\Delta^{2,3}f(v-u,y-x)& = k_\beta(x-y+v-u) D^2 \varphi(v-u,y-x).
\end{align*}

Then, by the preceding decomposition, 
$$J_n(s,x,y) \leq \sum^3_{i=1} J_n^i(s,x,y),$$
 where $J^i_n(s,x,y)$ is defined in the
same way as $J_n(s,x,y)$ but with $D^2 f$ replaced by $\Delta^{2,i} f$, $i=1,2,3$.
% = \displaystyle \int^t_0 ds \int_{\IR^3} du \int_{\IR^3} dv G_n(s,n)G_n(s,v) 
%\vert \Delta^{2,i}f(u,v,x,y) \vert, i = 1,2,3.
We now estimate the contribution of each of these terms.

   Since $\varphi$ is bounded, applying Lemma \ref{lemB.1}(c) and (d) with $b:=\alpha$, $a:=3-(\alpha+\beta)$, $\alpha +
\beta \in\, ]0,2[$, $u:=v-u$, $x:=y-x$,  and (\ref{A.3}), yields
\begin{align} 
J^1_n(s,x,y) &\leq \vert x-y \vert^\alpha \int_{\IR^3}du\, G_n(s,u) \int_{\IR^3} dv\, G_n(s,v)\nonumber\\
&\quad \times \int_{\IR^3}dw\, k_{\alpha+\beta}(\vert y-x\vert w + u-v) \vert D^2k_{3-\alpha}(w,e) \vert \nonumber\\
&\leq \vert x-y\vert^\alpha \Big(\sup_{x,y,w} \int_{\IR^3} du \int_{\IR^3} dv\, G_n(s,u)G_n(s,v)\nonumber\\
&\quad\times k_{\alpha+\beta}(\vert y-x\vert w + u-v)\Big)
\int_{\IR^3}dw\, \vert D^2 k_{3-\alpha}(w,e)\vert \nonumber\\
&\leq  C\vert x-y\vert^\alpha. 
\label{B.6}
\end{align}
The properties of the function $\varphi$ in Assumption \ref{assumpC} imply that
\begin{eqnarray*}
   J^2_n(s,x,y) &\leq& C \vert x-y\vert \displaystyle \int_{\IR^3} du \int_{\IR^3} dv\, G_n(s,u) G_n(s,v)\\
  \\
   && \qquad\qquad \times \vert k_\beta (y-x+v-u)-k_\beta (x-y+v-u)\vert\\
  \\
  &\leq& C \vert x-y\vert (J_n^{2,-}(s,x,y) + J_n^{2,+} (s,x,y)),
\end{eqnarray*}
where
$$
%\begin{array}{lll}
  J^{2,\pm}_n (s,x,y) = \displaystyle \int_{\IR^3} du \int_{\IR^3} dv\, G_n(s,u)G_n(s,v)\, \vert Dk_\beta(v-u,\pm (y-x))
\vert,\\
%\\
%  J^{2,2}_n (t,x,y) &=& \displaystyle\int^t_0 ds \int_{\IR^3} du \int_{\IR^3}dv\, G_n(s,u)G_n(s,v)\vert
%k_\beta(y-x+v-u)-k_\beta(v-u)\vert.
%\end{array}
$$
Notice that $J^{2,-}_n(s,x,y)$ coincides with the term $I^1_n(s,x,y)$ in Lemma \ref{lemB.2} when $\varphi \equiv 1$; similarly, $J^{2,+}_n(s,x,y) =
I^1_n(s,y,x)$ with $\varphi \equiv 1.$ Consequently, (\ref{B.3}) yields
\begin{equation}\label{B.7}
   J^2_n(s,x,y) \leq C \vert x-y \vert^{1+\alpha},
\end{equation}
with $\alpha \in\, ]0, (2-\beta) \wedge 1[$.

  The mean value theorem and the properties of $\varphi$ give
\begin{eqnarray*}
  && \vert D^2\varphi(v-u,x-y) \vert\\
  &&\qquad \leq \vert y-x\vert \int^1_0 d\lambda\, \vert \nabla \varphi(v-u+\lambda(x-y)) - \nabla
\varphi(v-u-\lambda(x-y))\vert \\
  && \qquad \leq C\, \vert y-x \vert^{1+\delta}.
\end{eqnarray*}
Therefore, by (\ref{A.3}),
\begin{eqnarray}\nonumber
   J^3_n (s,x,y) &\leq& C \vert y-x\vert^{1+\delta} \int_{\IR^3} du \int_{\IR^3}dv\, G_n(s,u) G_n(s,v) k_\beta(x-y+v-u)\\
\nonumber
   &\leq& C \vert y-x\vert^{1+\delta} \int_{\IR^3}d\xi\, k_{3-\beta}(\xi) \, \vert {\cal{F}} G(s)(\xi)\vert^2  \\
   &\leq& C \vert y-x\vert^{1+\delta}. 
\label{B.8}
\end{eqnarray}
Putting together the estimates (\ref{B.6})-(\ref{B.8}), we obtain (\ref{B.5}).
\hfill $\Box$
\vskip 16pt
%\qed

\begin{lemma} For any $0\leq s \leq t \leq \tbar \leq T$, set
$$
   \nu^n_1(s,t,\tbar) = \int_{\IR^3} du \int_{\IR^3} dv\, G_n(t-s,u)G_n(t-s,v) f\left(\frac{\tbar-s}{t-s}v - u\right)
\frac{\tbar-s}{t-s}.
$$
Then 
$$
   \sup_{n\geq 1}\sup_{0\leq s \leq t \leq \tbar \leq T} \nu^n_1(s,t,\tbar)< \infty.
$$
\label{lemB.4}
\end{lemma}

\noindent{\it Proof.} Consider the change of variables $v \mapsto \frac{\tbar-s}{t-s} v$. Then by Lemma \ref{lemA.1},
\begin{align*}
 \nu^n_1(s,t,\tbar)&= \int_{\IR^3} du \int_{\IR^3} dv\, G_n(t-s,u) G_n\left(t-s, \frac{t-s}{\tbar-s}v\right)\\
 \\
&\qquad\times f(v-u) \left(\frac{t-s}{\tbar-s}\right)^2\\
\\
&=\int_{\IR^3} du \int_{\IR^3} dv\, G_n(t-s,u) G_n(\tbar-s,v)f(v-u).
\end{align*}
Using Assumption \ref{assumpC}, the last integral is bounded by 
$$
\Vert\varphi\Vert_{\infty}\, \int_{\IR^3}du \int_{\IR^3}dv\, G_n(t-s,u) G_n(\tbar-s,v)k_\beta(v-u).
$$ 
Consequently,
\begin{align*}
   \nu^n_1(s,t,\tbar) &\leq C\, \int_{\IR^3} d\xi\, \frac{\vert {\cal{F}} G_n(t-s)(\xi)\vert \, \vert
{\cal{F}}G_n(\tbar - s)(\xi)\vert}{\vert \xi \vert^{3-\beta}}\\
\\
&= C\, \int_{\IR^3} d\xi\, \frac{\vert {\cal{F}} G(t-s)(\xi)\vert}{\vert \xi \vert^{(3-\beta)/2}}\,
 \frac {\vert  {\cal{F}} G(\tbar -s)(\xi))}{\vert \xi \vert^{(3-\beta)/2}}.
\end{align*}
Apply the Cauchy-Schwarz inequality and (\ref{A.3}) to see that this integral is bounded by a finite constant, uniformly in $s, t, \tbar \in[0,T].$ 
\hfill $\Box$
\vskip 16pt

   The next two lemmas deal with time-scaled increments of the covariance density.

\begin{lemma} For each $t, \tbar \in[0,T]$, with $t\le\tbar$; set
\begin{align*}
\nu^n_2(t, \tbar) &= \int^t_0 ds \int_{\IR^3} du \int_{\IR^3}dv\, G_n(t-s,u)G_n(t-s,v)\\
\\
&\quad\times \left\vert \left(\frac{\tbar-s}{t-s}\right)^2
f\left(\frac{\tbar-s}{t-s}(v-u)\right) - \frac{\tbar-s}{t-s} f \left(\frac{\tbar-s}{t-s} v-u\right) \right\vert.
\end{align*}
Then 
\begin{equation}\label{B.9}
\sup_{n \geq 1} \nu^n_2 (t,\tbar) \leq C \vert \tbar-t\vert^\alpha,
\end{equation}
for any $\alpha \in\, ]0,1[$ with $\alpha + \beta \in\, ]0,2[$.
\label{lemB.5}
\end{lemma}

\noindent{\it Proof.}  Consider the decomposition 
$$\nu^n_2(t, \tbar) \leq \nu^{n,1}_2(t, \tbar) + \nu^{n,2}_2(t, \tbar),$$
where $\nu^{n,1}_2(t, \tbar)$ (resp. $\nu^{n,2}_2(t,\tbar)$) is defined in the same way as $\nu^n_2(t,\tbar)$, but with
the expression in absolute values replaced by
$$
   \left\vert \left(\frac{\tbar-s}{t-s}\right)^2 - \frac{\tbar-s}{t-s}\right\vert\, f\left(\frac{\tbar-s}{t-s}(v-u)\right)
$$
and 
$$
   \frac{\tbar-s}{t-s}\,\left\vert f\left(\frac{\tbar-s}{t-s}(v-u)\right)-f\left(\frac{\tbar-s}{t-s}v-u\right)\right\vert,
$$
respectively.

The change of variables $(u,v) \mapsto \frac{\tbar-s}{t-s}(u,v)$ and Lemma \ref{lemA.1} yield
\begin{align}
   \nu^{n,1}_2(t,\tbar) &= \vert \tbar-t\vert \int_0^t \frac{ds}{\tbar-s} \int_{\IR^3}du \int_{\IR^3} dv\,
G_n(\tbar-s,u) G_n(\tbar-s,v)f(v-u)\nonumber\\ 
&\nonumber\\
 &\leq C \vert \tbar-t\vert \int_0^t \frac{ds}{\tbar-s} \int_{\IR^3} d\xi\, \frac{\vert{\cal{F}}
G(\tbar-s)(\xi)\vert^2}{\vert \xi \vert^{3-\beta}} \nonumber\\ 
&\nonumber\\
&\leq C \vert \tbar-t \vert,
\label{B.10}
\end{align}
if $\beta \in\, ]0,2[$, where in the last inequality we have applied Lemma \ref{lemA.3} with $b:=1$.

   Consider the transformation $(u,v) \mapsto (u, \frac{\tbar-s}{t-s} v).$ Then 
\begin{equation}\label{B.10'}
   \nu^{n,2}_2 (t, \tbar) = \int^t_0 ds \int_{\IR^3} du \int_{\IR^3}dv\, G_n(t-s,u)G_n(\tbar-s,v)\left\vert Df\left(v-u, -
\frac{\tbar-t}{t-s} u\right)\right\vert.
\end{equation}
We are going to prove that
\begin{equation}\label{B.11}
   \sup_{n\geq 1} \nu^{n,2}_2 (t, \tbar) \leq C \vert t-\tbar \vert^\alpha,
\end{equation}
with $\alpha \in\, ]0,1[$, $\alpha + \beta \in\, ]0,2[$.

   Notice that the structure of this term is similar to the integral analyzed in Lemma \ref{lemB.2}. Consider the
decomposition 
$$\nu^{n,2}_2 (t, \tbar) \leq \nu^{n,2,1}_2(t,\tbar)+\nu_2^{n,2,2}(t,\tbar),$$
 where
$\nu^{n,2,1}_2(t,\tbar)$ (resp. $\nu^{n,2,2}_2(t,\tbar)$) are as on the right-hand side of (\ref{B.10'}), but with the
absolute value replaced by
$$
\left(\varphi(v-u- \frac{\tbar-t}{t-s} u\right) \,\left\vert Dk_\beta\left (v-u - \frac{\tbar-t}{t-s} u\right)\right\vert,
$$
and 
$$
    \left\vert D\varphi \left(v-u, -\frac{\tbar-t}{t-s}u\right)\right\vert\, k_\beta(v-u),
$$
respectively.

Since $\varphi$ is bounded, Lemma \ref{lemB.1} (a) with $b = \alpha$, $a = 3-(\alpha+\beta)$ $c=\tbar-t$, 
$u:=v-u$, $x=-\frac{u}{t-s}$, yields
\begin{eqnarray*}
   \nu^{n,2,1}_2 (t, \tbar) &\leq& C \vert t-\tbar\vert^\alpha \displaystyle\int_0^t ds \int_{\IR^3} du \int_{\IR^3}dv\,
G_n(t-s,u)G_n(\tbar-s,v)\\
   \\
   && \qquad \times \displaystyle \int_{\IR^3} dw\, k_{\alpha+\beta}((\tbar-t)w-(v-u)) 
    \ \left\vert  Dk_{3-\alpha}(w, -\frac{u}{t-s})\right\vert.
\end{eqnarray*}

   We next prove that if $\alpha + \beta \in \, ]0,2[$, then the integral above is bounded, uniformly in $n$
   and $0\le t\le\tbar\le T$. For this, let
$\nu^{n,2,1,1}_2$ (resp.~$\nu^{n,2,1,2}_2$) be the same expression, but with $\IR^3$ in the $dw$-integral replaced by
$B_3(0)$ (resp. $B_3(0)^c$). For $\nu^{n,2,1,1}_2$, break the absolute value in the $dw$-integral into the sum of two
terms. Applying successively the changes of variables $w \mapsto w - \frac{u}{t-s}$ and $u \mapsto
\frac{\tbar-s}{t-s}u$ and (\ref{A.2}), the first term is bounded, by (\ref{A.3}). Indeed
%, by Lemma \ref{lemB.1},
\begin{align*}
&\int_0^t ds \int_{\IR^3} du \int_{\IR^3}dv\, G_n(t-s,u)G_n(\tbar-s,v)\\
\\
&\qquad\times \int_{\vert w\vert\leq3} dw\,
k_{\alpha+\beta}((\tbar-t)w-(v-u)) k_{3-\alpha}\left(\frac{u}{t-s}-w\right)\\
\\
&\quad \leq  \int_0^t ds \int_{\IR^3} du \int_{\IR^3}dv\, G_n(t-s,u)G_n(\tbar-s,v)\\
\\
&\qquad\times \int_{\vert w \vert\leq 5} 
dw\, k_{\alpha+\beta}\left((\tbar-t)w-v + \frac{\tbar-s}{t-s}u\right) k_{3-\alpha}(w) \\
\\
& \quad = \int^t_0 ds \int_{\IR^3} du \int_{\IR^3}dv\, G_n(\tbar-s,u) G_n(\tbar-s,v) \frac{t-s}{\tbar-s}\\
\\
&\qquad\times\int_{\vert w \vert \leq 5} dw\, k_{\alpha+\beta}((\tbar-t)w-(v-u)) k_{3-\alpha}(w)\\
\\
&\quad \leq \int_{\vert w \vert \leq 5} \frac{dw}{\vert w \vert^{3-\alpha}} \int_0^t ds
\int_{\IR^3}d\xi\, \frac{\vert {\cal{F}} G(\tbar-s)(\xi)\vert^2}{\vert\xi\vert^{3-(\alpha+\beta)}},
\end{align*}
which is bounded uniformly over $0\le t\le \tbar\le T$.
if $\alpha + \beta \in\, ]0,2[$.

   Similarly, the second term is
\begin{align*}
& \int_0^t ds \int_{\IR^3} du \int_{\IR^3}dv\, G_n(t-s,u)G_n(\tbar-s,v)\\
\\
&\qquad\times \int_{\vert w \vert \leq 3} dw\,
k_{\alpha+\beta}((\tbar-t)w-(v-u)) k_{3-\alpha}(w)\\
\\
& \quad \leq \int_{\vert w \vert \leq 3} \frac{dw}{\vert w \vert^{3-\alpha}} \int^t_0 ds \int_{\IR^3}
d\xi\, \frac{\vert{\cal{F}} G(t-s)(\xi)\vert\,  \vert{\cal{F}} G(\tbar-s)(\xi)\vert} {\vert \xi \vert^{3-(\alpha+\beta)}},
\end{align*}
which is bounded uniformly over $0\le t\le \tbar\le T$, if $\alpha + \beta \in\, ]0,2[$.

   Turning to $\nu_2^{n,2,1,2}$, let 
$$
   \psi(\lambda) = k_{3-\alpha}\left(w - \lambda \frac{u}{t-s}\right), \qquad \lambda \in [0,1].
$$ 
Notice that since $\vert u \vert \leq (t-s)(1+\frac{1}{n})$ and $\vert w \vert > 3$, it follows that 
$$\vert w-\lambda
\frac{u}{t-s} \vert \geq \vert \vert w \vert - \lambda (1+ \frac{1}{n})\vert \geq \vert w \vert -2 > \frac{1}{3} \vert w
\vert.$$
 Further, $\vert \psi'(\lambda)\vert \leq C \vert w - \lambda \frac{u}{t-s} \vert^{\alpha-4}$. Thus,
\begin{align*}
\nu_2^{n,2,1,2}(t,\tbar) &\leq  \int_0^t ds \int_{\IR^3} du \int_{\IR^3}dv\, G_n(t-s,u)G_n(\tbar-s,v) \\
\\
 & \quad \times \int_{\vert w \vert > 3} dw\,\vert k_{\alpha + \beta}(\tbar-t)w+u-v)\vert  \int^1_0
d\lambda\, k_{4-\alpha}\left(w-\lambda \frac{u}{t-s}\right)\\
\\
&\leq C\int_{\vert w \vert > 3} \frac{dw}{\vert w \vert^{4-\alpha}} \int_0^t ds \int_{\IR^3}d\xi\,
\frac{\vert {\cal{F}} G(t-s)(\xi)\vert\, \vert {\cal{F}} G(\tbar-s)(\xi)\vert}{\vert \xi\vert^{3-(\alpha+\beta)}} 
\end{align*}
which is bounded uniformly over $0\le t\le \tbar\le T$, if $\alpha \in\, ]0,1[$ and $\alpha + \beta \in\,]0,2[$.

   Therefore, we have proved that
\begin{equation}\label{B.12}
   \sup_{n \geq 1} \nu_2^{n,2,1}(t, \tbar) \leq C \vert t-\tbar\vert^\alpha,
\end{equation}
if $\alpha \in\,]0,1[$ and $\alpha + \beta \in\,]0,2[$.

   By the properties of $\varphi$, we have
$$
   \nu_2^{n,2,2}(t, \tbar) \leq C \vert\tbar-t\vert \int_0^t ds \int_{\IR^3} d\xi\, \frac{\vert {\cal{F}}
G_n(t-s)(\xi)\vert\, \vert{\cal{F}} G_n(\tbar-s)(\xi)\vert}{\vert \xi \vert^{3-\beta}}
$$
and therefore
\begin{equation}\label{B.13}
\sup_{n \geq 1} \nu^{n,2,2}_2 (t, \tbar) \leq C \vert \tbar-t \vert.
\end{equation}
The estimates (\ref{B.12}) and (\ref{B.13}) prove (\ref{B.11}). Finally, (\ref{B.9}) follows from (\ref{B.10}) and
(\ref{B.11}).
\hfill $\Box$
\vskip 16pt
%\qed

   The last lemma of this section gives information on the second-order increments of the time-scaled covariance function. Its proof is somewhat lengthy because we seek the best possible result: indeed, the case $\beta \in\, ]0,1[$ requires considerable effort, while the case $\beta \in [1,2[$ is not so complicated. 

\begin{lemma} For any $s,t,\tbar \in [0,T]$, $s \leq t \leq \tbar$, $u,v \in \IR^3$, we set
\begin{eqnarray}
   \Delta^2 f(s,t,\tbar, u,v) &=& \left(\frac{\tbar-s}{t-s}\right)^2 f \left(\frac{\tbar-s}{t-s} (v-u)\right) - \frac{\tbar-s}{t-s}
f\left(\frac{\tbar-s}{t-s} v-u\right)\nonumber\\
&\nonumber\\
&& \qquad\qquad - \, \frac{\tbar-s}{t-s} f\left(v-\frac{\tbar-s}{t-s}u\right) + f(v-u) \label{B.13.1}
\end{eqnarray}
and
$$
\nu^n_3(t,\tbar) = \displaystyle \int^t_0 ds \int_{\IR^3} du \int_{\IR^3} dv\, G_n(t-s,u)G_n(t-s,v)\, \vert \Delta^2
f(s,t,\tbar,u,v)\vert.
$$
Then
\begin{equation}\label{B.14}
\sup_{n \geq 1} \nu^n_3(t, \tbar) \leq C \vert t-\tbar\vert^\alpha,
\end{equation}
for any $\alpha \in\, ]0, (2-\beta) \wedge (1+\delta)[$.
\label{lemB.6}
\end{lemma}

\noindent{\it Proof.}  Assume first that $\beta \in [1,2[$. In this case, we consider a decomposition of $\Delta^2
f(s,t,\tbar,u,v)$ into first-order increments, which leads to the following decomposition:
$$
\nu^n_3(t, \tbar) = \nu^{n,1}_3 (t, \tbar) + \nu^{n,2}_3 (t, \tbar) + \nu^{n,3}_3(t, \tbar),
$$
where $\nu^{n,i}_3 (t,\tbar)$, $i=1,2,3$, is defined as $\nu^n_3 (t,\tbar)$, but with $\Delta^2 f$ replaced by
$\Delta^{2,i} f$, and
\begin{align*}
   \Delta^{2,1} f(s,t, \tbar,u,v) &= 
   \left\vert \left(\frac{\tbar-s}{t-s}\right)^2 f\left(\frac{\tbar-s}{t-s}(v-u)\right)
    - \frac{\tbar-s}{t-s} f\left(\frac{\tbar-s}{t-s}v-u\right)\right\vert,\\
     \\
   \Delta^{2,2} f(s,t, \tbar,u,v) &= 
   \left\vert 1- \frac{\tbar-s}{t-s} \right\vert\, f(v-u),\\
     \\
   \Delta^{2,3} f(s,t, \tbar,u,v) &=
   \frac{\tbar-s}{t-s}\, \left\vert f(v-u)-f\left(v-\frac{\tbar-s}{t-s} u\right)\right\vert.
\end{align*}
Notice that $\nu^{n,1}_3(t,\tbar)$ is equal to $\nu^n_2(t,\tbar)$ defined in Lemma \ref{lemB.5}. Thus, by (\ref{B.9}),
\begin{equation}\label{B.15}
   \sup_{n\geq 1} \nu^{n,1}_3(t, \tbar) \leq C \vert t-\tbar \vert^\alpha,
\end{equation}
for any $\alpha \in\, ]0,(2-\beta) \wedge 1[$.

   For the term $\nu^{n,2}_3(t, \tbar)$, we have
$$
   \nu^{n,2}_3(t, \tbar) \leq C \vert t-\tbar \vert \int_0^t \frac{ds}{t-s} \int_{\IR^3}d\xi\, \frac{\vert {\cal{F}}
G(t-s)(\xi)\vert^2}{\vert \xi \vert^{3-\beta}}.
$$
Then, by Lemma \ref{lemA.3} with $b =1$,
\begin{equation}\label{B.16}
   \nu^{n,2}_3 (t, \tbar) \leq C \vert t-\tbar \vert,
\end{equation}
for any $\beta \in\, ]0,2[$.

   The analysis of $\nu^{n,3}_3 (t, \tbar)$ does not differ very much from that of the term $\nu^{n,2}_2(t, \tbar)$ in the
proof of Lemma \ref{lemB.5}. Indeed, the change of variables $(u,v) \mapsto (\frac{\tbar-s}{t-s} u,v)$ and Lemma
\ref{lemA.1} yield
\begin{align*}
 \nu^{n,3}_3 (t, \tbar) &= \int^t_0 ds \int_{\IR^3} du \int_{\IR^3} dv\, G_n(\tbar-s,u) G_n(t-s,v)\, \\
 \\
&\quad\times \left\vert Df\left(v-u,\frac{\tbar-t}{\tbar - s}u\right) \right\vert\, .\\
\end{align*}
Comparing this expression with (\ref{B.10'}), we observe that in the integrands, the roles of $t$ and $\tbar$ are exchanged. 
However, carrying on calculations
similar to that for $\nu^{n,2}_2(t,\tbar)$ in the proof of Lemma \ref{lemB.5}, we encounter similar expressions, with $ds$ replaced
by $\frac{ds}{t-s}$. 
Consequently,  using Lemma \ref{lemA.3} instead of (\ref{A.3}), we obtain
\begin{equation}\label{B.17}
\sup_{n \geq 1} \nu^{n,3}_3 \leq C \vert t-\tbar \vert^\alpha,
\end{equation}
for each $\alpha \in\, ]0,(2-\beta) \wedge 1[$. Notice that for $\beta \in [1,2]$, $(2-\beta) \wedge 1 = 2-\beta$.
Consequently, the estimate (\ref{B.14}) follows in this case from (\ref{B.15})-(\ref{B.17}).
\smallskip

   Consider now the case $\beta \in\, ]0,1[$. We decompose $\Delta^2f(s,t,\tbar,u,v)$ into the sum $\sum^3_{i=1}
\bar\Delta^{2,i}f(s,t,\tbar,u,v)$, with the  definitions
\begin{align*}
\bar\Delta^{2,1}f(s,t,\tbar,u,v) &=\left(\left(\frac{\tbar-s}{t-s}\right)^2-2 \frac{\tbar-s}{t-s}+1\right)f(v-u)\\
\\
& = \left(\frac{\tbar-
t}{t-s} \right)^2 f(v-u),\\
\\
\bar\Delta^{2,2}f(s,t,\tbar,u,v) &= \left(\left(\frac{\tbar-s}{t-s}\right)^2 -
\frac{\tbar-s}{t-s}\right)\left(f\left(\frac{\tbar-s}{t-s}(v-u)\right)-f(v-u)\right),\\
\\
\bar\Delta^{2,3}f(s,t,\tbar,u,v) &=\frac{\tbar-s}{t-s}\tilde \Delta^{2,3}f(s,t,\tbar,u,v),
\end{align*}
where 
\begin{align*}
\tilde \Delta^{2,3}f(s,t,\tbar,u,v) & =
\Big(f\big(\frac{\tbar-s}{t-s}(v-u)\big) -f\big(\frac{\tbar-s}{t-s}v-u\big)-f\big(v-\frac{\tbar-s}{t-s}u\big)\\
&\quad + f(v-u)\Big).
\end{align*}
Then
$$
\nu^{n}_3(t,\tbar) \le \sum^3_{i=1} \overline{\nu}^{n,i}_3(t, \tbar),
$$
with $\overline{\nu}^{n,i}_3(t,\tbar)$ defined as $\nu^n_3(t,\tbar)$, but with $\Delta^2 f$ replaced by $\bar\Delta^{2,i} f$,
$i = 1,2,3.$

   We start by estimating $\overline{\nu}^{n,1}_3(t,\tbar)$. We have
$$
\overline{\nu}^{n,1}_3 (t, \tbar) \leq C \vert t-\tbar\vert^2 \int_0^t \frac{ds}{(t-s)^2} \int_{\IR^3}d\xi\,
\frac{\vert{\cal{F}} G(t-s)(\xi) \vert^2}{\vert \xi \vert^{3-\beta}}.
$$
Then, using Lemma \ref{lemA.3} with $b=2$, we conclude that
\begin{equation}\label{B.18}
   \sup_{n\geq 1} \overline{\nu}^{n,1}_3 (t, \tbar) \leq C \vert t-\tbar \vert^2,
\end{equation}
if $\beta \in\, ]0,1[$.

   We proceed now with $\overline{\nu}^{n,2}_3(t,\tbar)$, which we split into two terms:
\begin{align*}
\overline{\nu}^{n,2,1}_3 &\leq \vert t-\tbar \vert \int_0^t ds \int_{\IR^3} du \int_{\IR^3} dv\,
G_n(t-s,u)G_n(t-s,v) \frac{\tbar-s}{(t-s)^2}\\
\\
&\quad\times \left\vert f\left(\frac{\tbar-s}{t-s}(v-u)\right) - f\left(\frac{\tbar-s}{t-s}v-u\right)\right\vert,\\
\\
\overline{\nu}^{n,2,2}_3 &\leq \vert t-\tbar \vert  \int_0^t ds  \int_{\IR^3} du \int_{\IR^3} dv\,
G_n(t-s,u)G_n(t-s,v) \frac{\tbar-s}{(t-s)^2}\\
\\
 &\quad\times \left\vert f \left(\frac{\tbar-s}{t-s} v-u\right) - f(v-u) \right\vert.
\end{align*}
The change of variables $(u,v) \mapsto \left(u, \frac{\tbar-s}{t-s} v\right)$ and Lemma \ref{lemA.1} yield
\begin{align*}
 \overline{\nu}^{n,2,1}_3 &\leq \vert t-\tbar \vert \int_0^t \frac{ds}{t-s} \int_{\IR^3} du \int_{\IR^3}
dv\, G_n(t-s,u) G_n(\tbar-s,v)\\
\\
&\quad\times \left\vert Df\left(v-u, -\frac{\tbar-t}{t-s}u\right)\right\vert.
\end{align*}
Moreover,
\begin{align*}
\overline{\nu}^{n,2,2}_3 &\leq C \vert t-\tbar \vert \int_0^t \frac{ds}{(t-s)^2} \int_{\IR^3} du
\int_{\IR^3} dv\, G_n(t-s,u)G_n(t-s,v)\\
\\
&\quad\times \left\vert Df\left(v-u,\frac{\tbar-t}{t-s}v\right)\right\vert.
\end{align*}
Proceeding as for the analysis of the term $\nu^{n,3}_3$, or $\nu^{n,2}_2$ in Lemma \ref{lemB.5} (observe the similarity of $\nu^{n,2}_2$
in (\ref{B.10'}) with $\overline{\nu}^{n,2,1}_3, \overline{\nu}^{n,2,2}_3)$, and using Lemma \ref{lemA.3}, we obtain
\begin{equation}\label{B.19}
   \sup_{n \geq 1} \overline{\nu}^{n,2}_3 (t, \tbar) \leq C \vert t-\tbar\vert^{1+\rho_1},
\end{equation}
for any $\rho_1 \in\, ]0,1[$ (we recall that $\beta \in\,]0,1[$).
\smallskip

   The study of $\overline{\nu}^{n,3}_3(t,\tbar)$ is more intricate. Due to the product structure of the covariance
function, we have the following decomposition
$$
\bar \Delta^{2,3} f(s,t,\tbar,u,v) = \frac{\tbar-s}{t-s}\,\sum^4_{i=1} \Delta^{2,3,i} f(s,t,\tbar,u,v),
$$
where 
\begin{align*}
 \Delta^{2,3,1} f(s,t,\tbar,u,v) &= \tilde\Delta^{2,3} k_{\beta}(s,t,\tbar,u,v)\,\varphi \left(\frac{\tbar-s}{t-s}(v-u)\right),\\
\\
 \Delta^{2,3,2} f(s,t,\tbar,u,v) &= D\varphi\left(\frac{\tbar-s}{t-s}v-u,-\frac{\tbar-t}{t-s}u\right)\,Dk_{\beta}\left(v-u,\frac{\tbar-t}{t-s}v\right),\\
\\
  \Delta^{2,3,3}f(s,t,\tbar,u,v) &= D\varphi\left(v-\frac{\tbar-s}{t-s}u, \frac{\tbar-t}{t-s}v\right)\,Dk_{\beta}\left(v-u,-\frac{\tbar-t}{t-s}u\right),\\
\\
  \Delta^{2,3,4}f(s,t,\tbar,u,v) &=\tilde\Delta^{2,3} \varphi(s,t,\tbar,u,v)\,k_{\beta}(v-u).
\\
\end{align*}

Then $\overline{\nu}^{n,3}_3(t, \tbar) \leq \sum^4_{i=1} \overline{\nu}^{n,3,i}_3(t,\tbar),$ where
\begin{align*}
\overline{\nu}^{n,3,i}_3(t,\tbar)& = \int_0^t ds \int_{\IR^3}du \int_{\IR^3}dv\, G_n(t-s,u)G_n(t-s,v)\\
\\
&\quad\times\frac{\tbar-s}{t-s}\, \left\vert \Delta^{2,3,i} f\left(s,t,\tbar,u,v\right)\right\vert.
\end{align*}
We shall now analyze the contribution of each of these terms.

By Lemma \ref{lemB.1} (e), applied with $d=3$, $b:=\alpha$, $a:=3-(\alpha+\beta)$, $c:=\tbar-t$, $u:=v-u$,
$x:=-\frac{u}{t-s}$, $y:=\frac{v}{t-s}$,
\begin{equation*}
\bar \nu_3^{n,3,1}\le |\tbar-t|^{\alpha} I(t,\tbar),
\end{equation*}
where
\begin{align*}
I(t,\tbar) &= \int_0^t ds\, \frac{\tbar-s}{t-s} \int_{\IR^3} du \int_{\IR^3} dv\, G_n(t-s,u) G_n(t-s,v)\\
\\
&\quad \times \int_{\IR^3} dw\, k_{\alpha+\beta}(v-u-(\tbar-t)w)\left\vert \bar D^2 k_{3-\alpha}\left(w, -\frac{u}{t-s},\frac{v}{t-s}\right)\right\vert.
\end{align*}

 Our next objective is to check that $\sup_{0\le t\le \tbar\le T}I(t,\tbar) < \infty$. This will be carried out by taking into account first the small values of
$w$ ($\vert w \vert \leq 5$), and then the large ones.

Set $\tbar-t = h$ and split the last absolute value above into four terms, so that the above integral when $w$ is integrated
on the set $\vert w \vert \leq 5$ equals to $\sum_{j=1}^4 I_j(t,\tbar)$.

Clearly,
\begin{align*}
I_1(t,\tbar) &:= \int_0^t ds\, \frac{\tbar-s}{t-s} \int_{\IR^3} du \int_{\IR^3} dv\, G_n(t-s,u) G_n(t-s,v)\\
\\
&\quad \times \int_{\vert w \vert < 5} dw\, k_{\alpha+\beta}(v-u-(\tbar-t)w) k_{3-\alpha}(w)\\
\\
&\le 
 \left(\int_{\vert w \vert < 5} \frac{dw}{\vert w \vert^{3-\alpha}}\right) \int_0^t ds \, \frac{\tbar-s}{t-s}\, \int_{\IR^3} d\xi\,
\frac{{\cal{F}} G(t-s)(\xi)\vert^2}{\vert \xi \vert^{3-(\alpha+ \beta)}} < \infty,
\end{align*}
if $\alpha + \beta \in\, ]0,2[$, $\alpha > 0.$

The analysis of $I_j(t,\tbar)$, $j = 2,3,4$, are all similar, so we only give the details for 
\begin{align*}
I_4(t,\tbar) & =  \int_0^t ds \frac{\tbar-s}{t-s} \int_{\IR^3} du \int_{\IR^3} dv\, G_n(t-s,u)G_n(t-s,v)\\
&\quad\times\int_{\vert w \vert < 5} dw\, k_{\alpha+\beta}( v-u - h w ) k_{3-\alpha}\left(w+\frac{v-u}{t-s}\right).
\end{align*}
Consider the change of variables $\wbar = w + \frac{v-u}{t-s}$. Then
\begin{align*}
I_4(t,\tbar) & \leq \int_{\vert w \vert < 7} \frac{dw}{\vert w \vert^{3-\alpha}} \int_0^t ds  \frac{\tbar-s}{t-s} \int_{\IR^3} du
\int_{\IR^3} dv\, G_n(t-s,u) G_n(t-s,v)\\
&\quad\times k_{\alpha+\beta}\left(hw+\frac{\tbar-s}{t-s}(v-u)\right).
\end{align*}
Apply the change of variables $(u,v) \mapsto \frac{\tbar-s}{t-s}(u,v).$ Then, by Lemma \ref{lemA.1},
\begin{align*}
 I_4(t,\tbar) &\leq \int_{\vert w \vert < 7} \frac{dw}{\vert w \vert^{3-\alpha}} \int_0^t ds\, \frac{t-s}{\tbar-s}   \int_{\IR^3} du
\int_{\IR^3} dv\, G(\tbar-s,u) G(\tbar-s,v)\\
\\
 &\quad\times k_{\alpha+\beta}\left(hw+(v-u)\right)\\
 \\
 &\leq C \int_0^{\tbar} ds\,  \int_{\IR^3}d\xi\, \frac{\vert{\cal{F}} G(t-s)(\xi)\vert^2}{\vert \xi
\vert^{3-(\alpha+\beta)}},
\end{align*}
which is bounded, uniformly over $0\le t\le\tbar\le T$, if $\alpha + \beta \in\, ]0,2[$, $\alpha > 0.$

Therefore, we have proved that the contribution to $I(t,\tbar)$ of the term in which the $dw$-integral is restricted to $\{\vert w
\vert < 5 \}$ is finite. 
\smallskip

   We end the proof by checking that $I_5(t,\tbar) < \infty$, where $I_5(t,\tbar)$ is defined in the same way as $I(t,\tbar)$ but with the
$dw$-integral restricted to the domain $\{\vert w\vert \geq 5\}$.

   Define 
$$
\Psi(\lambda,\mu)= 
k_{3-\alpha}( w - \lambda\frac{u}{t-s} + \mu \frac{v}{t-s}),\, \lambda, \,\mu \in[0,1],
$$ 
so that
\begin{align*}
\bar D^2\,k_{3-\alpha}(w,-\frac{u}{t-s},\frac{v}{t-s}) &= \Psi(1,1)-\Psi(0,1)-\Psi(1,0)+\Psi(0,0)\\
& = \int^1_0 d\lambda  \int_0^1 d \mu\, \frac{\partial^2\Psi}{\partial\lambda\partial\mu}(\lambda,\mu).
\end{align*}
A simple computation shows that
$$
  \left\vert \frac{\partial^2\Psi}{\partial\lambda\partial\mu}(\lambda,\mu)\right\vert \leq C  k_{5-\alpha}\left(w- \lambda \frac{u}{t-s}+ \mu \frac{v}{t-s}\right).
$$
Consequently 
\begin{align*}
I_5(t,\tbar) &:= \int_0^t ds\, \frac{\tbar-s}{t-s} \int_{\IR^3} du \int_{\IR^3} dv\, G_n(t-s,u) G_n(t- s,v)\\
\\
&\quad\times \int_{\vert w \vert \geq 5} dw\, k_{\alpha+\beta}(v-u-hw)\,\tilde D^2 k_{3-\alpha}\left(w,-\frac{u}{t-s},
\frac{v}{t-s}\right)\\
\\
&\leq  C \int_0^t ds\, \frac{\tbar-s}{t-s} \int_{\IR^3} du \int_{\IR^3} dv\, G_n(t-s,u) G_n(t- s,v)\\
\\
&\quad\times \int_{\vert w \vert \geq 5} dw\, k_{\alpha+\beta}(v-u-hw)\\ 
\\
&\quad\times\int^1_0 d \lambda \int^1_0 d\mu\, k_{5-\alpha}\left(w-\lambda \frac{u}{t-s} + \mu
\frac{v}{t-s}\right).
\end{align*}
We can give a lower bound of $\vert w - \lambda \frac{u}{t-s} + \mu \frac{v}{t-s} \vert,$ independent of $u,v$ on
the set $\{\vert w \vert \geq 5\}$. Indeed, by the triangle inequality 
\begin{align*}
\left \vert w-\lambda \frac{u}{t-s} + \mu\, \frac{v}{t-s}\right \vert &\ge \left\vert\, \vert w \vert - \left\vert \lambda
\frac{u}{t-s} - \mu\, \frac{v}{t-s} \right\vert\, \right\vert\\
& \geq \vert\, \vert w \vert -2(\lambda+\mu) \vert\\
&\geq \vert w \vert - 4 \geq \frac{\vert w \vert}{5}.
\end{align*}
Hence,
\begin{align*}
 I_5 (t,\tbar)&\leq C  \int_{\vert w \vert \geq 5}  \frac{dw}{\vert w \vert^{5-\alpha}} \int_0^t ds\,
\frac{\tbar-s}{t-s} \int_{\IR^3} du \int_{\IR^3} dv\, G_n(t-s,u) G_n(t-s,v)  \\
\\
& \quad\times  k_{\alpha+\beta}(v-u-hw) \\
\\
&\leq C \int_{\vert w \vert \geq 5}  \frac{dw}{\vert w \vert^{5-\alpha}} \int_0^t ds \frac{\tbar-s}{t-s}\int_{\IR^3} d\xi\,
\frac{\vert{\cal{F}} G(t-s)(\xi)\vert^2}{\vert \xi \vert^{3-(\alpha+\beta)}}.
\end{align*}
The $dw$-integral is finite for any $\alpha \in\,]0,2[$ and, by Lemma \ref{lemA.3},  the $ds\,d\xi$-integral is also finite for $\alpha + \beta \in\,
]0,2[$. Consequently, $I_5(t,\tbar)$ is finite whenever $\alpha+\beta \in\,]0,2[$; therefore
\begin{equation}\label{B.20}
\sup_{n \geq 1} \overline{\nu}^{n,3,1}_3 (t, \tbar) \leq C \vert t-\tbar \vert^{\alpha}
\end{equation}
for any $\alpha \in\,]0, 2-\beta[$.

 By the properties of $\varphi$, we have
\begin{align*}
\overline{\nu}^{n,3,2}_3(t,\tbar) &\leq C \int_0^t ds\, \frac{\tbar-s}{t-s} \int_{\IR^3}du
\int_{\IR^3}dv\, G_n(t-s,u)G_n(t-s,v)\\
\\
&\quad\times \frac{\vert t-\tbar\vert}{t-s}\, \vert u \vert\, \left\vert D k_\beta \left(v-u, \frac{\tbar-t}{t-s}v\right)\right\vert\\
\\
&\leq C \vert t-\tbar \vert \int_0^t ds \frac{1}{t-s}\int_{\IR^3}du \int_{\IR^3}dv\, G_n(t-s,u)G_n(t-s,v)\\
\\
&\quad \times \left\vert D k_\beta \left(v-u, \frac{\tbar-t}{t-s}v\right)\right\vert .
\end{align*}
In the last inequality, we have used the fact that $\vert u \vert \leq (\tbar-s)(1+\frac{1}{n}).$

We can study the contribution of this last integral with similar arguments as those used in the analysis of
the term $\nu^{n,2,1}_2(t, \tbar)$ in
the proof of Lemma \ref{lemB.5},  concluding that 
\begin{equation}\label{B.21}
   \sup_{n \geq 1} \bar\nu_3^{n,3,2}\leq C \vert t-\tbar \vert^{1+\rho_3},
\end{equation}
for any $\rho_3 \in\, ]0, (2-\beta) \wedge 1[$ (see (\ref{B.12})). Clearly, the same bound holds for the term
$\bar\nu_3^{n,3,3}$.

By the assumptions on $\varphi,$ we have
\begin{align*}
& \left\vert \varphi \left(\frac{\tbar-s}{t-s} (v-u)\right)-\varphi
 \left(\frac{\tbar-s}{t-s}v-u\right)-\varphi \left(v-\frac{\tbar-s}{t-s}u\right)+\varphi(v-u)\right\vert\\
 \\
&\qquad \leq \int_0^1 d \lambda  \left\vert \nabla \varphi \left(\frac{\tbar-s}{t-s}v-u+\lambda
\frac{t-\tbar}{t-s}u\right)-\nabla \varphi \left(v-u+\lambda \frac{t-\tbar}{t-s}u\right) \right\vert\\
\\ 
&\qquad \qquad \times\frac{\vert t-\tbar \vert}{t-s}\, \vert u\vert
\leq C \frac{\vert t-\tbar \vert^{1+\delta}}{(t-s)^{1+\delta}} \,\vert u \vert \,\vert v
\vert^\delta.
\end{align*}
It follows that
\begin{align*}
\overline{\nu}^{n,3,4}_3(t,\tbar) &\leq C \vert t-\tbar \vert^{1+\delta}  \displaystyle\int_0^t ds \frac{\tbar-s}{t-s}
\int_{\IR^3}du \int_{\IR^3}dv\, G_n(t-s,u)G_n(t-s,v)\\ 
\\
&\quad\times\frac{\vert u\vert\, \vert v\vert^\delta}{(t-s)^{1+\delta}}
 k_\beta \displaystyle(v- u).
\end{align*}
The fraction with $u$ and $v$ is bounded by $2^{1+\delta}$. Therefore
\begin{equation*}
\overline{\nu}^{n,3,4}_3(t,\tbar) \leq   C \vert t-\tbar \vert^{1+\delta} \int_0^t ds \frac{\tbar-s}{t-s}\int_{\IR^3}d\xi\, \frac{\vert {\cal{F}}
G(t-s) (\xi)\vert \vert \overline{{\cal{F}} G(t-s)(\xi)}}{\vert \xi \vert^{3-\beta}},
\end{equation*}
and by Lemma \ref{lemA.3},
\begin{equation}\label{B.22}
   \sup_{n\geq 1} \overline{\nu}_3^{n,3,4} (t, \tbar) \leq C \vert t-\tbar\vert^{1+\delta}
\end{equation}
if $\beta \in\,]0,2[$.
\smallskip

Summarising the results obtained in (\ref{B.20})-(\ref{B.22}) yields
\begin{equation}\label{B.23}
   \sup_{n \geq 1} \overline{\nu}^{n,3}_3 (t, \tbar) \leq C \vert t-\tbar \vert^\alpha,
\end{equation}
with $\alpha \in\, ]0,(2-\beta) \wedge (1+\delta)[$.

   Finally, by (\ref{B.18}), (\ref{B.19}) and (\ref{B.23}), we conclude that
$$
   \sup_{n \geq 1}\nu^n_3 (t, \tbar) \leq C \vert t-\tbar \vert^\alpha,
$$
wiht $\alpha \in\, ]0,(2-\beta) \wedge (1+\delta)[$, when $\beta \in\,]0,1[.$
The proof of the lemma is complete.
\hfill$\Box$
\vskip 16pt
\noindent{\bf Acknowledgment:} 
The second named author is pleased to thank the
{\it Institut de Math\'ematiques} of the {\it Ecole Polytechnique F\'ed\'erale de Lausanne} for 
its hospitality and financial support during a visit where part of this work was carried out.

\end{document}